\documentclass[11pt]{article}

\usepackage[margin=1in]{geometry}

\usepackage{amsmath}
\usepackage{amsthm}
\usepackage{graphicx}
\usepackage{xcolor}
\usepackage{hyperref}  % load before cleveref
\usepackage{cleveref}

\theoremstyle{plain}
\newtheorem{theorem}{Theorem}[section]

\theoremstyle{definition}
\newtheorem{definition}[theorem]{Definition}

\newcommand{\newsiamremark}[2]{\theoremstyle{definition}\newtheorem{#1}[theorem]{#2}}

\newcommand{\email}[1]{\href{mailto:#1}{\texttt{#1}}}
\usepackage{sam}

\title{Fast Construction of Efficient Cut Cell Quadratures}
\author{Samuel F. Potter\thanks{Robert McNeel \& Associates, Seattle, WA (\email{sam@mcneel.com}).}}
\date{}

\begin{document}

\maketitle

\begin{abstract}
  Accurate quadrature on cut cells is a recurring difficulty in high-order immersed and interface methods, where meshing can be both expensive and fragile. Optimization-based formulations provide an alternative: they are fast, avoid meshing entirely by requiring only samples inside the cut cell, and apply to a wide class of domains, including but not limited to level set representations. This survey develops the moment-matching viewpoint underlying these formulations and organizes a family of methods into a unified ``quadrature toolbox''. These include approaches based on least squares, Steinitz elimination, nonnegative least squares, Xiao--Gimbutas elimination, Newton iteration, basis pursuit, and linear programming. Within this framework, several refinements of existing methods and a few new heuristics appear naturally. These include a new proof of the positivity of least squares quadrature rules for a fine enough grid, a simplified form of Steinitz elimination with improved asymptotic behavior, observations about stabilized nonnegative least squares, and a clarification of the relationship between basis pursuit and Ryu--Boyd quadrature. We also describe two practical warm starts for Newton’s method based on the first primal or dual iterates of a primal–dual interior-point method applied to the Ryu--Boyd LP. In one dimension, the dual strategy reduces to locating critical points of a related Legendre polynomial, suggesting how structural descriptions of this sort might extend to higher dimensions. Overall, the optimization viewpoint clarifies how to construct fast, robust, high-order quadrature rules without relying on meshing, and provides a common language for comparing and combining a broad set of techniques.
\end{abstract}

\section{Introduction}

Numerical integration over complicated and unstructured domains shows up in many settings: cut cells in finite element methods, trimmed surfaces in CAD, implicitly defined regions in level set methods, and so on.
In each, the task of ``computing a quadrature rule'' is usually treated as a single, monolithic problem which is tightly coupled to the underlying geometry.
Requiring that quadrature rules have positive weights and have nodes contained within the domain of integration, this task naturally decomposes into three conceptually distinct modular pieces:
\begin{enumerate}
\item an algorithm to decide whether a point lies inside or outside the domain (the \emph{inside-outside test}),
\item a means of integrating a family of test functions over each cell (\emph{moment integration}),
\item and an algorithm for using these two pieces of information to produce a quadrature rule.
\end{enumerate}
The first piece is fundamentally geometric and has been studied at length in many different contexts.
The second piece can be viewed as a conceptually separate subproblem in its own right.
In many cases it can even be formulated recursively: moment integration over a cut cell typically involves integrating fluxes over the cut cell boundary, which is itself a quadrature problem, albeit one that is easier or better structured than the original.
The third---what we refer to as ``quadrature optimization''---can be formulated independently of any particular geometry under consideration once the first two pieces are in place.

This paper is \emph{only} about quadrature optimization, although we discuss the inside-outside tests and moment integration in passing.
We are also primarily interested in integration on \emph{cut cells}: regions obtained as the result of intersecting a larger domain with a regular background cell, such as a square or cube.
We assume that, for each cut cell, we have at our disposal an inside-outside test and a means of integrating all moments in a finite-dimensional polynomial space.
Together, these provide a means of discretizing the cut cell into a grid of points contained within the cell and a moment vector associated with that cell.
The central question is, then: given this grid and moment vector, how should one choose nodes and weights so that the resulting quadrature rule has good accuracy and stability properties, and reflects whatever structural preferences (positivity, locality, symmetry, sparsity) the application demands?
Posed this way, the problem is no longer tied to any particular representation of geometry: different geometric settings simply feed different moment vectors and discretizations into the same optimization machinery.

The goal of this work is to make that optimization machinery explicit and usable.
We develop a small ``quadrature toolbox'' that produces quadrature rules by solving sparse, underdetermined moment-matching problems, with a focus on cut cell applications.
Within this framework, a large family of algorithms can be used to generate quadratures, including simple least squares solves, Steinitz-type elimination, nonnegative least squares, Xiao--Gimbutas elimination, basis pursuit, and Ryu and Boyd's linear programming formulation.
These methods differ in cost, robustness, and the structural properties of the resulting rules, but they all operate on the same underlying algebraic object: the moment-matching system.
The emphasis throughout is that quadrature optimization is modular, geometry-agnostic, and amenable to experimentation: once an inside-outside test and moment integrals are supplied, the optimization layer can run independently to produce high quality rules.

A common way to obtain moments on complicated domains is to mesh or otherwise decompose the cut cell and then apply standard quadratures elementwise.
This sidesteps the need for an inside-outside test and provides a straightforward means of evaluating moments, but it does so at the cost of additional machinery that is substantially more complex than what is required otherwise.
Meshing can be intricate, fragile, and slow; and, even when a mesh is available, mapping reference Gauss rules onto mesh elements is likely to overintegrate relative to the polynomial space of interest.
By contrast, in the absence of a mesh, moments can typically be computed directly---often by pushing integrals to the boundary via the divergence theorem---without constructing intermediate subcells.
From the perspective of quadrature optimization, then, meshing is unnecessary overhead: once an inside–outside test and moment integrals are available, the construction of a quadrature rule should be viewed as an independent numerical task.
This is the viewpoint adopted in this paper.

\subsection{Contributions}

This paper is mainly expository, in that it systematizes the moment-matching perspective and organizes the landscape of optimization-based methods into a unified framework (the ``quadrature toolbox'').
However, within this structure, we make a number of new contributions to the individual methods discussed:
\begin{itemize}
\item \emph{Least squares quadrature.}
  We extend an existence result of Huybrechs to multiple dimensions, ensuring the existence of Tchakaloff-efficient quadratures under mild conditions.
\item \emph{Steinitz elimination.}
  We present a simplified approach to Steinitz elimination with an improved asymptotic complexity.
\item \emph{Nonnegative least squares quadrature.}
  We revisit the Lawson--Hanson algorithm and provide a careful study of a stabilized version of the algorithm, originally due to Lawson and Hanson but, to our knowledge, not discussed in the literature.
\item \emph{Xiao--Gimbutas elimination.}
  We do a careful comparison of different node significances and clarify their relation to the moment-matching residual.
\item \emph{Basis pursuit quadrature.}
  We clarify the relationship between basis pursuit quadrature and Ryu--Boyd quadrature, showing that it is essentially a degenerate form of Ryu--Boyd quadrature in which the linear program reduces to feasibility testing.
\item \emph{Ryu--Boyd quadrature---a primal heuristic.}
  We carefully consider the behavior of a standard primal-dual interior point method as applied to the Ryu--Boyd quadrature optimization problem and show that it is possible to take only a single primal-dual step before enough information is available to produce a warm start for Newton's method applied to the moment-matching equations.
  Since this guess is based on an inspection of the primal variable of the iteration, and since we do not have a theory of convergence, we refer to it as a ``primal heuristic'' for the Ryu--Boyd problem.
\item \emph{Ryu--Boyd quadrature---a dual heuristic.}
  By studying the dual variables of the same interior point method, we are able to circumvent the linear programming formalism entirely, producing an even simpler heuristic for generating a warm start for the moment-matching equations.
  We provide a theoretical justification for this in 1D using a simple argument based on the theory of orthogonal polynomials.
\end{itemize}

\section{Preliminaries}

\subsection{Notation and basic definitions}

Throughout, $d > 0$ denotes the dimension of the ambient space $\mathbb{R}^d$.
A multi-index is written $\malpha = (\alpha_1, \hdots, \alpha_d)$, where $\alpha_i \geq 0$ for each $i = 1, \hdots, d$.
If $\mx = (x_1, \hdots, x_d) \in \bbR^d$, a $d$-variate monomial is written:
\begin{equation}
  \mx^{\malpha} = x_1^{\alpha_1} \cdots x_d^{\alpha_d}.
\end{equation}
We write $|\malpha| = \alpha_1 + \cdots + \alpha_d \geq 0$ for the degree of a multi-index.

Let $I$ be an index set.
For each $i \in I$, let $\malpha_i$ be a multi-index and let $c_i \in \bbR$ be a coefficient.
The total degree of the polynomial $p(\mx) = \sum_{i \in I} c_i \mx^{\malpha_i}$ is:
\begin{equation}
  \deg p = \max_{i\in I} |\malpha_i|.
\end{equation}
The so-called partial degree of $p$ is $\max_i \max_j (\malpha_i)_j$, which coincides with spaces of tensor product polynomials.
We will restrict our attention to spaces of polynomials up to a given total degree, reflecting their importance in the theory of approximation.
The set of $d$-variate polynomials of total degree at most $N$ is denoted:
\begin{equation}
  \bbP_N^d = \Span\Set{\mx^{\malpha} : \dim\malpha = d, |\malpha| \leq N},
\end{equation}
where $\dim\bbP_N^d = {N + d \choose d}$.
Note that the dimensions of $\bbP_N^2$ and $\bbP_N^3$ are the triangular and tetrahedral numbers.
Throughout, we let $m = \dim\bbP_N^d$.

\subsection{Background on quadrature}

Throughout, we will refer to a \emph{quadrature rule} (or, simply, a \emph{quadrature}) as a set of nodes $\mx_1, \hdots, \mx_n \in \bbR^d$ and weights $w_1, \hdots, w_n \in \bbR$.
The focus of this work is on computing quadratures:
\begin{enumerate}
\item with nodes contained in an integration domain of interest, $\Omega \subseteq \mathbb{R}^d$, which may be nonconvex (hence, we require $\mx_1, \hdots, \mx_n \in \Omega$),
\item and with positive weights: $w_i > 0$ for $i = 1, \hdots, n$.
\end{enumerate}
The number of nodes $n$ is also referred to as the \emph{order} of the quadrature.
We define the operator:
\begin{equation}
  \calQ\squareb{f} = \sum_{i=1}^n w_i f(\mx_i),
\end{equation}
to represent the application of the quadrature to a function $f$ defined on $\bbR^d$.
When there is no risk of confusion, we may refer to the operator $\calQ$ as the quadrature.
The number of \emph{degrees of freedom} of the quadrature is $\dof(\calQ) = (d + 1)n$.
Although the methods presented here can be adapted to a nonuniform density (or more general measure) $\mu$ supported on $\Omega$, we restrict our attention to $\mu \equiv 1$ for simplicity and because of its centrality in applications.

For $d = 2, 3$, we will frequently have a cut cell in mind as a representative example of $\Omega$.
To this end, we frequently assume that $\Omega$ is a subset of ${[-1, 1]}^{d}$ such that $\partial\Omega$ and $\partial[-1, 1]^d$ overlap nontrivially ($\Omega$ is cut from $[-1, 1]^d$).

The $\malpha$th monomial moment on $\Omega$ is:
\begin{equation}
  I_{\malpha} = \int_\Omega \mx^{\malpha}d\mx,
\end{equation}
where $\mx^{\malpha} = \mx_1^{\alpha_1}\cdots\mx_d^{\alpha_d}$.
If $\phi_1, \hdots, \phi_m$ form a basis for $\bbP_N^d$, we define moments in terms of each $\phi_i$ analogously. E.g., the $\phi_i$ moment is the integral $\int_\Omega \phi_i(\mx)d\mx$.
If $\calQ\squareb{p} = \int_{\Omega}p(\mx)d\mx$ for all $p \in \bbP_N^d$, we say that $\calQ$ \emph{integrates} $\bbP_N^d$.
Moments will be assumed known and computed by external means---typically via boundary integration using the divergence theorem.

Following Xiao and Gimbutas~\cite{xiaoNumericalAlgorithmConstruction2010}, the \emph{efficiency} of an order $n$ quadrature that integrates $\bbP_N^d$ is defined to be:
\begin{equation}
  \eff(\calQ) = \frac{\dim{\bbP_N^d}}{\dof(\calQ)} = \frac{m}{(d + 1)n}. \label{eq:efficiency}
\end{equation}
This definition can be extended to other finite-dimensional function spaces straightforwardly.
Our focus is developing fast algorithms for quickly computing reasonably efficient quadratures.
We expect $0 \leq \eff(\calQ) \leq 1$ to hold, which is indeed the case, although the upper bound is not obvious---see \Cref{thm:stroud}.
Gaussian quadrature in 1D obtains $\eff(\calQ) = 1$, but little about the maximum efficiency of multidimensional quadratures is known.

To illustrate, let $n > 0$, let $x_1, \hdots, x_n$ be the zeros of the $n$th Legendre polynomial, and let $w_i = \int_{-1}^1 \ell_i(x)dx$, where $\ell_i(x)$ is the $i$th Lagrange basis polynomial determined by the nodes $x_1, \hdots, x_n$. The corresponding quadrature $\calQ_1$ is the order $n$ Gauss--Legendre rule. This quadrature integrates all polynomials in $\bbP_{2n-1}^1$, hence $\eff(\calQ_1) = 1$. Quadratures for rectangles in dimensions greater than one can be obtained by taking tensor products of Gauss--Legendre rules. 
\begin{example}[Efficiency of tensor product Gauss--Legendre quadrature]
   Let $d > 1$ and define the order $n$ tensor-product Gauss--Legendre quadrature $\calQ_d$ by the formula:
  \begin{equation}
    \calQ_d\squareb{f} = \sum_{i_1=1}^n \cdots \sum_{i_d=1}^n w_{i_1} \cdots w_{i_d} f(x_{i_1}, \hdots, x_{i_d}).
  \end{equation}
  This quadrature integrates all monomials $\mx^{\malpha}$ where $\max_i \alpha_i \leq 2n - 1$; hence, it integrates all monomials $\mx^{\malpha}$ for which $|\malpha| \leq 2n - 1$. To reiterate: the former set of monomials is the set of monomials of partial degree at most $2n - 1$, the latter set is $\bbP_{2n-1}^d$. The latter is a subset of the former, and, over all possible $N$, the largest $\bbP_N^d$ so contained. Consequently, with respect to the total degree, we can see that $\calQ_d$ ``over-integrates''.

  Since $\dof(\calQ_d) = (2n)^d$, and since $\dim(\bbP_{2n-1}^d) = {2n - 1 + d \choose d}$, we have:
  \begin{equation}
    \begin{split}
      \eff(\calQ_d) &= \frac{{2n - 1 + d \choose d}}{(d+1)n^d} = \frac{2n(2n+1)\cdots(2n+d-1)}{(d+1)d!n^d} \\
                    &= \frac{(2n)^d + \cdots}{(d+1)!n^d} \sim \frac{2^d}{(d+1)!} \mbox{ as } n \to \infty.
    \end{split}
  \end{equation}
  This recovers $\eff(\calQ_1) = 1$, and $\eff(\calQ_2) \sim 2/3$ and $\eff(\calQ_3) \sim 1/3$. Clearly, $\eff(\calQ_d) \sim 0$ if $d \to \infty$ and $n \to \infty$ simultaneously.
\end{example}

We generically refer to a quadrature for a known, simple geometry (a canonical domain) with an efficiency higher than what could otherwise be achieved by ``conventional means'' as an efficient quadrature. 
For $d > 1$, very little is known about the maximum achievable efficiency of quadratures for canonical domains, let alone more exotic ones.
If $\Omega$ isn't canonical, a standard procedure for numerical integration is to decompose $\Omega$ into a disjoint union of subdomains, and define each subdomain as the image of a canonical domain.
By using a known quadrature for the preimage of each subdomain and suitably adjusting its weights we get a mapped quadrature.
Even if decomposition can be avoided, mapping can have a deleterious effect on efficiency.

\begin{example}[Efficiency of mapped quadratures]\label{ex:mapped-quadrature-efficiency}
  This example illustrates the inefficiency of mapped quadratures and motivates the native, moment-based approach used throughout this paper.
  Let $\Omega$ be a deformation of the unit square with one side given by the polynomial $h(x) \geq 0$:
  \begin{equation}
    \Omega = \Set{(x, y) \in \bbR^2 : 0 \leq x \leq 1 \mbox{ and } 0 \leq y \leq h(x)}.
  \end{equation}
  We integrate the monomial $x^p y^q$ over $\Omega$ by a change of variables:
  \begin{equation}
    \int_{\Omega} x^p y^q dA(x, y) = \int_0^1 \int_0^1 x^p (y h(x))^q h(x) dy dx.
  \end{equation}
  We denote the polynomial degree with respect to $x$ of the integrand by $r(p, q) = p + (q + 1) \deg h$; the degree with respect to $y$ is $q$.
  Consequently, to integrate $x^p y^q$ we can first apply integration by substitution as above and then integrate over the unit square using a Gauss--Legendre rule.
  This requires an $\lceil r(p, q)/2 \rceil$-point rule to integrate over $x$ and a $\lceil q/2 \rceil$-point rule for $y$; consequently, to integrate each monomial in $\mathbb{P}_n^2$, we need a rule with $\lceil n/2 \rceil$ points to integrate in $y$ and $\lceil r^* / 2 \rceil$ points to integrate in $x$, where:
  \begin{equation}
    r^* = \max_{0 \leq p \leq n} r(p, n - p).
  \end{equation}
Denote this rule by $\mathcal{Q}_{n,h}$.
Its efficiency satisfies:
\begin{equation}
  \eff(\mathcal{Q}_{n,h}) \sim \left. \left( \frac{n(n+1)}{2} \right) \middle/ \left(\frac{3}{4} r^* n \right) \right. = \frac{2(n+1)/3}{\max_{0 \leq p \leq n} \Big(p + (n - p + 1) \deg h\Big)}.
\end{equation}
Since the maximum in the denominator is bounded below by $(n + 1) \deg h$, we have:
\begin{equation}
  \eff(\mathcal{Q}_{n,h}) \lesssim \frac{2}{3} \Big( \deg h \Big)^{-1}.
\end{equation}
This recovers $\eff(\calQ_2) \sim 2/3$ if $h$ is linear, as expected.
But we can now see that if $h$ is of even moderate degree, the efficiency of the rule decays rapidly, resulting in a large number of quadrature points needed to integrate $\mathbb{P}_n^2$ exactly.
\end{example}

It can be argued that if $h$ in \Cref{ex:mapped-quadrature-efficiency} is smooth, then a reasonably high-order Gauss--Legendre rule will be accurate enough to avoid compromising convergence rates in a simulation.
Indeed, numerically integrating the moments with near minimum error is often unnecessary.
For example, the effect of quadrature on the order of convergence of finite element methods has been studied extensively, giving clear guidelines for how to systematically underintegrate~\cite{ciarletCombinedEffectCurved1972,babuskaEffectNumericalIntegration2011}.
However, when $h(x)$ is only piecewise smooth, the need to decompose the cut cell becomes much more acute: mapping breaks down, and the efficiency of the resulting quadrature deteriorates rapidly (although, see \Cref{sec:pruning-inefficient-rules} for a prescription).
A native quadrature rule defined directly on the cut cell avoids this issue entirely.
For simplicity, we elect to integrate moments as accurately as we can to avoid confounding our results with this additional source of error.
Systematic underintegration of the moments does not meaningfully affect the qualitative results presented in this work.

\subsection{The need for native cut cell quadratures}

The particular focus of this work is designing quadrature rules for cut cells, primarily for use in immersed methods~\cite{burmanCutFEMDiscretizingGeometry2015,wunschEnrichedImmersedFinite2025}.
Basic approximation theory tells us that if the continuity class of the boundary is low, mapped quadrature is doomed.
If $h(x)$ in \Cref{ex:mapped-quadrature-efficiency} is a polynomial approximating the function $y(x) = |x - 1/2| + 1/2$, then it is well known that tens of thousands (if not hundreds of thousands) of coefficients (e.g., computed in the Chebyshev basis) are needed for a decent approximation.
The obvious remedy is to further subdivide $\Omega$---that is, to build a mesh on it.
But avoiding meshing is the motivation for choosing an immersed method in the first place.
By punting and moving the possibility of mesh failure elsewhere, we have not solved our original problem.
Consequently, there is a genuine need to construct native cut cell quadrature rules without recourse to meshing.
As we will see, such rules can not only be constructed without meshing, they are typically dramatically more efficient.
In particular, it is easy to construct rules that integrate $\mathbb{P}_n^d$ with an efficiency of $1/(d+1)$.
With a little work, we can push that efficiency even higher.

\subsection{Moment-matched quadratures}\label{sec:moment-matched-quadratures}

Quadratures for non-canonical domains can be generated using numerical optimization to ensure that the quadrature integrates each moment correctly.
These moment-matched quadratures are our focus.
The approach is quite general and can be used to design quadratures for a wide variety of function spaces.
For simplicity, and since it is the case of the greatest practical importance, we restrict our attention to $\bbP_N^d$.
\begin{definition}
  Let $w_1, \hdots, w_n \in \bbR$ and let $\mx_1, \hdots, \mx_n \in \bbR^d$, and write $\calQ\squareb{f} = \sum_{i=1}^n w_i f(\mx_i)$.
  Let $\bbP_N^d = \Span\Set{\phi_1, \hdots, \phi_m}$.
  The moment-matching equations for the quadrature $\calQ$ with respect to $\bbP_N^d$ are:
  \begin{equation}
    I_i = \int_\Omega \phi_i(\mx)d\mx = \calQ\squareb{\phi_i} = \sum_{j=1}^n w_j \phi_i(\mx_j), \qquad i = 1, \hdots, m. \label{eq:moment-matching-equations}
  \end{equation}
  Letting $\mw = (w_1, \hdots, w_n)$ and $\mX = (\mx_1, \hdots, \mx_n) \in {(\bbR^d)}^n$, the residual function $\mF : \bbR^n \times \bbR^{dn} \to \bbR^m$ for the moment-matching equations is $\mF_i(\mw, \mX) = I_i - \calQ\squareb{\phi_i}$.
\end{definition}

A quadrature integrates $\bbP_N^d$ if the moment-matching equations are satisfied.
If we fix $\mX$ and find $\mw$ that satisfies \eqref{eq:moment-matching-equations}, then the moment-matching equations are linear in $\mw$.
Defining the Vandermonde type matrix
\begin{equation}
  \mV = \begin{bmatrix}
    \phi_1(\mx_1) & \cdots & \phi_m(\mx_1) \\
    \vdots & \ddots & \vdots \\
    \phi_1(\mx_n) & \cdots & \phi_m(\mx_n)
  \end{bmatrix} \in \bbR^{n \times m},
\end{equation}
this linear system reads
\begin{equation}
  \mV^\top \mw = \m{I}.
\end{equation}
If we seek a pair $(\mw, \mX)$ that simultaneously satisfies \eqref{eq:moment-matching-equations}, then the system is nonlinear, and we turn to the evaluation of its residual.
We write:
\begin{equation}\label{eq:moment-matching-jacobian}
  \sfD\mF = \begin{bmatrix}
    \sfD_{\mw}\mF & \sfD_{\mX}\mF
  \end{bmatrix},
\end{equation}
where:
\begin{equation*}
  \sfD_{\mw}\mF = \begin{bmatrix}
    \phi_1(\mx_1) & \cdots & \phi_1(\mx_n) \\
    \vdots & \ddots & \vdots \\
    \phi_m(\mx_1) & \cdots & \phi_m(\mx_n) 
  \end{bmatrix}, \;
  \sfD_{\mX}\mF = \begin{bmatrix}
    w_1 \sfD\phi_1(\mx_1) & \cdots & w_n \sfD\phi_1(\mx_n) \\
    \vdots & \ddots & \vdots \\
    w_1 \sfD\phi_m(\mx_1) & \cdots & w_n \sfD\phi_m(\mx_n)
  \end{bmatrix},
\end{equation*}
and $\sfD\phi_i = \begin{bmatrix} \partial\phi_i/\partial x_1 & \cdots & \partial\phi_i/\partial x_d \end{bmatrix} \in \bbR^{1 \times d}$.

As we discuss in \Cref{sec:gauss-newton}, we will often generate a quadrature that nearly satisfies \Cref{eq:moment-matching-equations}.
Iteratively improving such a rule with Newton's method is highly effective.
Because the moment-matching equations are typically underdetermined, Newton's method is used in its generalized (pseudoinverse) form.
This update coincides with the standard Gauss--Newton step for the residual, so the distinction between ``Newton'' and ``Gauss--Newton'' is entirely notational; we use the terms interchangeably.

\subsection{Existence of positive quadratures}\label{sec:existence-of-quadratures}

Although the task of finding a stable, positive quadrature appears daunting at first sight, Martinsson et al.~\cite{martinssonInterpolationIntegrationFinitedimensional2007} showed that a stable but not necessarily positive quadrature exists for a set of bounded functions where
\begin{enumerate}[nosep]
\item the integration region is taken to be essentially any set,
\item the quadrature has the same number of nodes as the number of functions being integrated,
\item and the quadrature can be computed from a suitable discretization of the set using a rank-revealing QR decomposition~\cite{guEfficientAlgorithmsComputing1996}.
\end{enumerate}
Such a quadrature is often a reasonable choice.

Nevertheless, we seek a quadrature with positive weights since they are expected by practitioners.
There are several reasons for this:
\begin{itemize}
\item They are (arguably) more stable than quadratures with negative weights~\cite{huybrechsStableHighorderQuadrature2009,keshavarzzadehNumericalIntegrationMultiple2018}.
  For instance, following Huybrechs~\cite{huybrechsStableHighorderQuadrature2009}, if $\kappa(\mathcal{Q}) = \sum_i |w_i|$, and if $f$ and $\tilde{f}$ are such that $|f(\mx) - \tilde{f}(\mx)| \leq \epsilon$ for $\mx \in \Omega$ and for $\epsilon > 0$, then:
  \begin{equation}
    \big|\mathcal{Q}[f] - \mathcal{Q}[\tilde{f}]\big| \leq \epsilon \kappa(\mathcal{Q}).
  \end{equation}
  Clearly, $\kappa(\mathcal{Q})$ is minimized if $\mathcal{Q}$ is a positive quadrature.
  Likewise, a quadrature with negative weights could amplify the effect of round-off error.
\item In many industrial finite element simulations (e.g., in nonlinear elasticity simulations with large deformations and contact, such as the simulation of a car crash~\cite{wriggersComputationalContactMechanics2006}), ``mass lumping'' is used extensively to enable efficient explicit time-stepping.
  One popular approach to mass lumping is the ``nodal quadrature method'', in which quadrature nodes are used for finite element collocation nodes---consideration of the Lagrangre basis implies that the mass matrix must be diagonal~\cite{voetMathematicalTheoryMass2023}.
  Stability of the scheme is ensured by using a positive quadrature~\cite{cohenGeometricModelingSplines2001}.
\item There have been claims and some numerical evidence presented that Newton's method applied to nonlinear problems is more stable when using nonnegative moment-matched quadratures than ones with negative weights~\cite{garhuomNonnegativeMomentFitting2022b,hubrichNumericalIntegrationDiscontinuous2017,hubrichNumericalIntegrationNonlinear2019}.
\end{itemize}
The idea that quadrature nodes may additional be used as collocation nodes is an important consideration.
Huybrechs makes it clear that nodes which are unstable for polynomial interpolation can be used without any trouble as quadrature nodes~\cite{huybrechsStableHighorderQuadrature2009}.
On the other hand, using optimization to explicitly seek interpolation nodes with a good Lebesgue constant was explored by Hesthaven for the simplex~\cite{hesthavenElectrostaticsAlmostOptimal1998}.
We do not consider this issue further.

Positive quadratures with $\eff(\mathcal{Q}) = 1/(d+1)$ exist under very general conditions.
\begin{theorem}\label{thm:tchakaloff}
  Let $\Omega \subseteq \mathbb{R}^d$ be a compact set, let $m > 0$, let $\phi_1, \hdots, \phi_m : \Omega \to \mathbb{R}$ be a set of $m$ linearly independent functions, and let $\bbP = \Span\,\{\phi_1, \hdots, \phi_m\}$.
  Then, there exists a positive quadrature with order not greater than $m$, with nodes in $\Omega$, and which integrates $\bbP$.
\end{theorem}

\begin{proof}
  This is a specialization of Tchakaloff's theorem to the present setting~\cite{bayerProofTchakaloffsTheorem2006}.
  Davis provided a constructive proof of Tchakaloff's theorem~\cite{davisConstructionNonnegativeApproximate1967}.
  We provide our own constructive proof adapted to the setting of this paper in \Cref{sec:least-squares-quadrature,sec:steinitz-elimination}.
\end{proof}

\noindent In the context of this paper, we will refer to \Cref{thm:tchakaloff} as Tchakaloff's theorem, bearing in mind that the actual theorem of Tchakaloff is more general.

For a generic quadrature, $\eff(\mathcal{Q}) \leq 1$ is necessary generically.
Consider the map $\mF : \bbR^n \to \bbR^{dn} \to \bbR^m$.
If $\eff(\mathcal{Q}) > 1$, then $m > (d+1)n$---that is, the dimension of the range of $\mF$ is greater than the dimension of its domain.
In this case, since $\mF$ is a smooth mapping, its image is a codimensional subset of $\bbR^m$.
There could be a structural relationship between $\mF$ and $\mI$ which allows the moment-matching equations to be satisfied regardless in special cases---but, this is not the case in practice, especially for unstructured domains such as cut cells.

Lower bounds on the number of quadrature nodes are available.
A classical lower bound is due to Stroud~\cite{stroudApproximateCalculationMultiple1971} (see also Mysovskikh~\cite{mysovskikhApproximationMultipleIntegrals1980}).
\begin{theorem}\label{thm:stroud}
  Let $\mathcal{Q}$ be a positive quadrature exact on $\bbP_N^d$ and let $n$ be the number of nodes in $\mathcal{Q}$.
  Then:
  \begin{equation}
    n \geq {\kappa + d \choose d} \quad \mbox{where} \quad \kappa = \left\lfloor\frac{N}{2}\right\rfloor.
  \end{equation}
\end{theorem}
This bound can be sharpened using additional structural information about the commutators arising in the recursions defining systems of multivariate orthogonal polynomials, but we will not pursue this here~\cite{xuCommonZerosPolynomials2020}.

In the foregoing discussion, we consider quadratures which are positive, minimal (having at most $m$ nodes, matching Tchakaloff's theorem), and exact---i.e., which exactly satisfy the moment-matching equations for the quadrature.
If we relax the last requirement, the existence proof is simpler.

\begin{theorem}\label{thm:caratheodory-quadrature}
  Let $\mX \subseteq \mathbb{R}^d$ be a point set and let $n = |\mX| < \infty$.
  Let $\bbP = \Span\,\{\phi_1, \hdots, \phi_m\}$.
  Let $\mV = \begin{bmatrix} \phi_1(\mX) & \cdots & \phi_m(\mX) \end{bmatrix}$ and let $\mI \in \mathbb{R}^m$ such that $\mI_j = \int_\Omega \phi_j(\mx)d\mx$.
  Assume that $\mw^*$ minimizes:
  \begin{equation}
    \begin{split}
      \mbox{minimize}\qquad & \left\|\mV^\top \mw - \mI\right\|_2 \\
      \mbox{subject to}\qquad & \mw \geq \m{0}.
    \end{split}
  \end{equation}
  Then, there exists $\mw^\star \geq \m{0}$ such that $\|\mV^\top \mw^\star - \mI\|_2 \leq \|\mV^\top \mw^* - \mI\|_2$ and $\operatorname{card}(\mw^\star) \leq m$.
\end{theorem}

\begin{proof}
  First, we recall Caratheodory's theorem, which states that any point in the convex hull of a set in $\mathbb{R}^m$ can be written as a convex combination of $m + 1$ points in the original set~\cite{barvinokCourseConvexity2002}; equivalently, any point in the conic hull of a set in $\mathbb{R}^m$ can be written as a conic combination of $m$ points taken from the generating set.
  
  Now, let $\m{R}^* = \m{I} - \mV^\top \mw^*$.
  Of course, $\m{R}^* \neq \m{0}$ in general.
  Since $\mV^\top \mw^* = \m{I} - \m{R}^*$, we can see that $\m{I} - \m{R}^*$ is a conic combination of the columns of $\mV^\top$; that is, $\m{I} - \m{R}^* \in \operatorname{cone}\,\Set{\mphi(\mx_1), \hdots, \mphi(\mx_n)}$.
  Then, by Caratheodory's theorem, $\mI - \mR^*$ can be written as a conic combination of at most $m$ points in $\Set{\mphi(\mx_1), \hdots, \mphi(\mx_n)}$.
  Let $\mw^\star \in \mathbb{R}^n$ be the nonnegative vector whose $j$th component is given by Caratheodory's theorem and zero otherwise; consequently, $\mV^\top \mw^\star = \mI - \mR^*$.
\end{proof}

\noindent \Cref{thm:caratheodory-quadrature} is \emph{not} constructive.
Caratheodory's theorem merely indicates that an $m$-sparse combination of columns of $\mV^\top$ exists, but not how to compute it.

\section{Related work}\label{sec:related}

We give a high-level overview of existing successful algorithms for quadrature generation, most of which use some form of numerical optimization or numerical linear algebra.
We emphasize those approaches most relevant to the ``quadrature optimization'' viewpoint adopted in this work.

When comparing these works, we have to bear several different performance characteristics in mind:
\begin{itemize}
\item The error in the quadrature rules. This will be due to the approximation of the boundary of the domain and the approximation of the moments.
The latter will typically be affected by the former.
For each of these methods, it is straightforward to achieve either a fixed (high) order of accuracy, or to drive this error as close to numerical roundoff as possible.
Consequently, error alone is not a meaningful axis of comparison.
\item Whether the rules have positive weights and whether the nodes are contained inside the domain (strictly in the interior or with nodes allowed on the boundary---we reject rules with nodes outside $\Omega$).
\item Issues related to numerical stability.
\item How quickly the rules can be computed.
\item The efficiency of the rules (eq. \eqref{eq:efficiency}).
\end{itemize}
The last three points are the axes along which there is the most difference.

\subsection{Algebraic approaches}\label{sec:algebraic-approaches}

Although our focus is on using numerical optimization and methods from the scientific computing literature to compute quadrature rules, there is a substantial body of work based on algebraic techniques.
Since this literature provides a trove of information on the quadrature problem and is historically older and more established, we mention it first.
We restrict ourselves to describing a few significant works which can serve as entry points to this literature:
\begin{itemize}
\item Cools provides a readable survey of the field as of 1997~\cite{coolsConstructingCubatureFormulae1997}. The methods are restricted to simple canonical domains with obvious symmetries (all of $\mathbb{R}^n$, the cube, the ball, the simplex, etc.). The focus is on methods involving either polynomial ideal theory or invariant theory.
\item Collowald and Hubert provide an approach to constructing multidimensional quadratures on unstructured domains using moment theory~\cite{collowaldAlgorithmsComputingCubatures2018}. In particular, they consider eigenvalue problems involving the Hankel operator, a generalization of the one-dimensional moment matrix that underlies the construction of Gaussian quadratures on the interval. In principle, this work generalizes the approach to computing two-dimensional quadratures due to Vioreanu and Rokhlin (\Cref{sec:vioreanu-rokhlin-quadrature}) to multiple dimensions---possibly even to nonconvex domains~\cite{vioreanuSpectraMultiplicationOperators2014}. 
  We stress that the experimental results presented in this work are done symbolically (i.e., using arbitrary precision rational numbers).
  Numerical instability is an active concern, since the underlying multivariate Hankel matrices are highly ill-conditioned, and the connection to multidimensional rootfinding inherits the numerical difficulties identified by Graf et al.~\cite{grafNumericalInstabilityAlgebraic2024}.
\item Riener and Schweighofer are able to apply abstract optimization theory to the moment problem to develop algorithms and derive strong upper bounds on the order of certain quadrature rules in the plane~\cite{rienerOptimizationApproachesQuadrature2018}.
\end{itemize}

\subsection{Huybrechs' accurate equispaced quadratures}\label{sec:huybrechs}

Foundational work on using numerical methods to construct quadrature rules was done by Huybrechs~\cite{huybrechsStableHighorderQuadrature2009}.
In this work, it was shown that stable quadratures with equispaced nodes (or a subset thereof) can be designed for the interval $[a, b] \subseteq \mathbb{R}$ by solving least squares and nonnegative least squares problems.
Let $x_1 < \cdots < x_n$ be an equispaced grid such that $x_1 = a$ and $x_n = b$ and let $\bbP_m^1$ denote the polynomial space being integrated.
This is a signal work for us, as many of the methods developed here can be viewed as natural multidimensional extensions of Huybrechs’ approach.

Huybrechs shows that there exists some $n = O(m^2)$ such that for all $n \geq m$, the least squares quadrature (see \Cref{sec:least-squares-quadrature}) has positive weights, providing a fast algorithm for computing these weights which exploits properties of polynomials that are orthogonal with respect to the discrete uniform measure supported on the equispaced grid (see, e.g., Gautschi~\cite{gautschiOrthogonalPolynomialsComputation2004}).
He also presents some preliminary results showing that the nonnegative least squares quadrature (see \Cref{sec:nonnegative-least-squares-quadrature}) computed on a grid with $O(n^2)$ nodes produces a positive quadrature with $m + 1$ nonzero weights.
Huybrechs uses MATLAB's implementation of nonnegative least squares (\texttt{lsqnonneg}), which uses Lawson and Hanson's active set method~\cite{lawsonSolvingLeastSquares1995}, but which masks the underlying algorithmic structure of the Lawson--Hanson method used and obscures key performance characteristics that matter in higher dimensions.

\subsection{Vioreanu--Rokhlin quadrature}\label{sec:vioreanu-rokhlin-quadrature}

If $\Omega \subseteq \mathbb{C}$ is compact and convex and $V$ is a finite-dimensional subspace of $L^2(\Omega)$, Vioreanu and Rokhlin showed that the eigenvalues of the complex multiplication operator $\mathcal{M}_{x + iy} : V \to L^2(\Omega)$ lie within $\Omega$, providing numerical evidence that these eigenvalues supply an excellent warm start for a nonlinear rootfinder applied to the moment-matching equations~\cite{vioreanuSpectraMultiplicationOperators2014}.
If $\Omega$ possesses any symmetries, the eigenvalues of the multiplication operator have the same symmetries.
Afterwards, quadrature weights are computed using least squares, and the subsequent quadrature is optimized by applying Gauss--Newton to the moment-matching equations (any symmetries that might exist in the eigenvalues are not enforced).
Unfortunately, if $\Omega$ is nonconvex, the eigenvalues may not lie within $\Omega$, somewhat limiting the applicability of the method to canonical two-dimensional domains which are convex or nearly convex.
The method also has trouble generating quadratures for domains with continuous symmetry groups, due to the spectral degeneracy of the multiplication operator.
We note again that the connection between quadrature and the multiplication operator, along with extensions of these ideas to higher dimenions, is explained in the work of Collowald and Hubert~\cite{collowaldAlgorithmsComputingCubatures2018} (\Cref{sec:algebraic-approaches}).
In short, the method provides high quality warm starts for convex planar domains, but the extension to nonconvex geometries or three dimensions remains unclear (but see Collowald and Hubert~\cite{collowaldAlgorithmsComputingCubatures2018}, discussed in \Cref{sec:algebraic-approaches} above).

\subsection{Methods based on recursive moment evaluation}

Several works recursively reduce the problem in dimensionality in a divide and conquer approach.
Moments on cell edges are first evaluated, followed by cell faces, in order to compute easier integrals using standard quadratures.
This is done in different ways:
\begin{itemize}
\item Initial work along these lines assuming a polyhedral cut cell was done in the convex case by Mousavi and Sukumar~\cite{mousaviNumericalIntegrationPolynomials2011} and the nonconvex case by Sudhakar and Wall~\cite{sudhakarQuadratureSchemesArbitrary2013}.
\item M\"{u}ller et al.\ explored this approach in an early work~\cite{mullerHighlyAccurateSurface2013}. The crux of their approach is to use the divergence theorem to recursively reduce volume and surface integrals into the constituent pieces of their respective boundaries. However, nonnegativity and interior containment are not enforced, which limits direct applicability. Also, while the use of the divergence theorem for efficiency is crucial, their treatment relies on simplifying smoothness assumptions that obscure the intrinsic complexity of moment computation for general cut cells. In particular, they assume that the domains are smooth and that the level set functions defining them lie in the basis being integrated.
\item Saye gave an algorithm for recursively constructing mapped tensor product Gauss rules on cut cells defined as the region bounded by the graph of a function, specifically for use in a discontinuous Galerkin simulation~\cite{sayeHighOrderQuadratureMethods2015}.
  He assumes a level set domain, and uses a method reminiscient of the Taylor model to bound the partials of the level set function to determine regions where the zero level set can be parametrized as the graph of a function~\cite{makinoTaylorModelsOther2003,neumaierTaylorFormsUse2003}.
\item Gunderman et al.\ present algorithms which recursively apply numerical antidifferentiation, the generalized Stokes' theorem, and Gaussian quadrature to construct rules for areas and volumes bound by rational curves and trimmed NURBS surfaces~\cite{gundermanSpectralMeshFreeQuadrature2021,gundermanHighAccuracyMeshFreeQuadrature2021}. 
\end{itemize}

\subsection{Ryu--Boyd quadrature}

Ryu and Boyd explored using linear programming to design positive quadratures~\cite{ryuExtensionsGaussQuadrature2015}.
They formulate the problem as a semi-infinite linear program (LP) over a measure space, and produce a finite-dimensional LP by approximating the continuous measures with discrete measures.
They prove that the semi-infinite LP recovers Gaussian quadratures on the interval in one dimension.
Practically speaking, they place a regular grid of quadrature nodes inside $\Omega$ and solve an LP whose primal variables are the quadrature weights at each node.
The LP has nonnegativity constraints enforcing positivity and equality constraints ensuring that the polynomial space is integrated.
The choice of cost function for the LP is open-ended and has a dramatic effect on the type of quadrature produced---they illustrate the effect of different choices with numerical examples.
Once the LP is solved, the quadrature must be extracted and further optimized.
The strength of the method is its adaptability to different integration domains.
They do not consider algorithmic details or carry out performance studies---we consider this method in more detail in \Cref{sec:ryu-boyd-quadrature} and show how it can be made efficient in \Cref{sec:fast-primal-heuristic}.

We note that Mehrotra and Papp explored a simpler version of Ryu and Boyd's semi-infinite LP several years earlier~\cite{mehrotraGeneratingMomentMatching2013}.
Their focus is on developing a method alternative to quasi-Monte Carlo methods in stochastic optimization.
While Ryu and Boyd consider a semi-infinite LP whose cost function is the quadrature applied to an arbitrary residual function (see \Cref{sec:ryu-boyd-quadrature} for details), Mehrotra and Papp use the total mass of the quadrature weights for their LP's cost.
This is equivalent to choosing a constant function for the residual in Ryu and Boyd's formulation.
We show in \Cref{sec:basis-pursuit-quadrature} that this reduces to feasibility testing for the sampled, finite-dimensional version of this LP; hence, it is unclear whether it offers any practical advantage over the nonnegative least squares algorithms presented in \Cref{sec:nonnegative-least-squares-quadrature}.

Mehrotra and Papp also present an algorithm for solving the finite-dimensional dual LP in polynomial time which merits further investigation as a method of solving Ryu and Boyd's semi-infinite LP.
With the choice of a nontrivial residual function, this could offer a compelling method of computing cut cell quadratures.
However, it should be noted that while Ryu and Boyd follow the solution of their LP with a clustering and reoptimization phase, which significantly improves the efficiency of their quadrature rules, Mehrotra and Papp's method appears to only achieve the efficiency guaranteed by Tchakaloff's theorem (\Cref{thm:tchakaloff}).

\subsection{Xiao--Gimbutas quadrature}

The goal of the work by Xiao and Gimbutas~\cite{xiaoNumericalAlgorithmConstruction2010} is to develop an algorithm for computing moment-matched (see \Cref{sec:moment-matched-quadratures}) positive quadratures on canonical domains which have nearly the minimum number of nodes possible.
In summary, they start by creating a quadrature for the domain through some combination of mapping and decomposition, where known positive quadratures are used for each of the resulting subdomains.
Optionally, they explicitly enforce any symmetries of the domain that are desired in the quadrature,  reducing the number of variables in the moment-matching equations.
To create a warm start for a Gauss--Newton iteration, they select a minimal subset of the nodes of the mapped and decomposed quadrature using a rank-revealing QR decomposition, along the lines of what is suggested in Martinsson et al.~\cite{martinssonInterpolationIntegrationFinitedimensional2007}.
They then proceed to solve the moment-matching equations and eliminate the ``least significant'' node (see \Cref{sec:xiao-gimbutas-elimination}) one at a time, using recursive backtracking to search over all sequences of eliminated nodes until a near-optimal quadrature is found.
This potentially gives a worst case complexity which is exponential in the number of candidate nodes.
So, while they are able to compute very nearly optimal quadratures, their algorithm is computationally intensive (timing statistics are not reported).
We note that the goal of their work is very different: they seek nearly optimal quadratures which are usable on standard mesh elements (e.g. triangles, tetrahedra, quadrilaterals, etc.), and are content to pay a high upfront cost.
This is achieved with great effect, with Xiao--Gimbutas quadratures being some of the most efficient known in a variety of cases (see, e.g., the Python library \href{https://github.com/inducer/modepy}{modepy}~\cite{kloecknerModepy2024}).

\subsection{Other quadratures based on numerical optimization}

A variety of other works in a similar vein have been carried out.
Thiagarajan and Shapiro solve the moment matching equations with a fixed grid using a QR decomposition~\cite{thiagarajanAdaptivelyWeightedNumerical2014}.
Glaubitz considers basis pursuit and least squares quadratures for use with experimental data, typically high-dimensional; some results on the existence and positivity of these quadratures are proved~\cite{glaubitzConstructingPositiveInterpolatory2020,glaubitzStableHighorderCubature2021,glaubitzConstructionApplicationProvable2023}.
Other similar studies on generating high-dimensional quadratures using methods from optimization have been carried out by Narayan and collaborators~\cite{keshavarzzadehNumericalIntegrationMultiple2018,jakemanGenerationApplicationMultivariate2018}.
With Monte Carlo applications in mind, Elefante solves a sequence of nonnegative least squares problems on low-discrepancy point sets until a small enough error in the moment-matching equations is achieved~\cite{elefanteCQMCImprovedCode2022}.
Bui et al.\ explore using moment-matched cut cell quadratures for elastoplastic simulations using the finite cell method~\cite{buiEfficientCutcellQuadrature2020}.
Garhuom et al.\ use nonnegative least squares to prune adaptive quadtree rules~\cite{garhuomNonnegativeMomentFitting2022b}.
The diversity of these methods underscores the broad applicability of quadrature optimization, and highlights the need for a unified framework that treats quadrature generation as a standalone numerical problem.

\section{Quadrature toolbox}

\begin{figure}
  \includegraphics[width=0.5\linewidth]{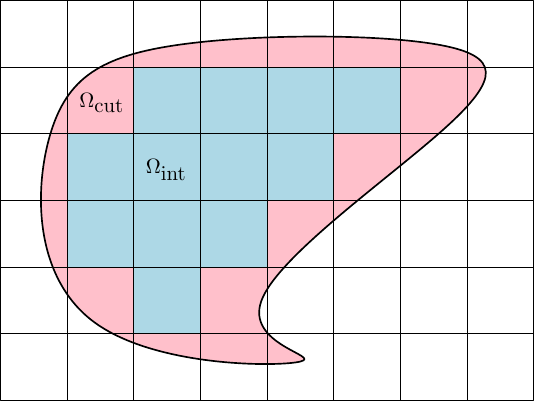}%
  \hspace{0.0399\linewidth}%
  \includegraphics[width=0.46\linewidth]{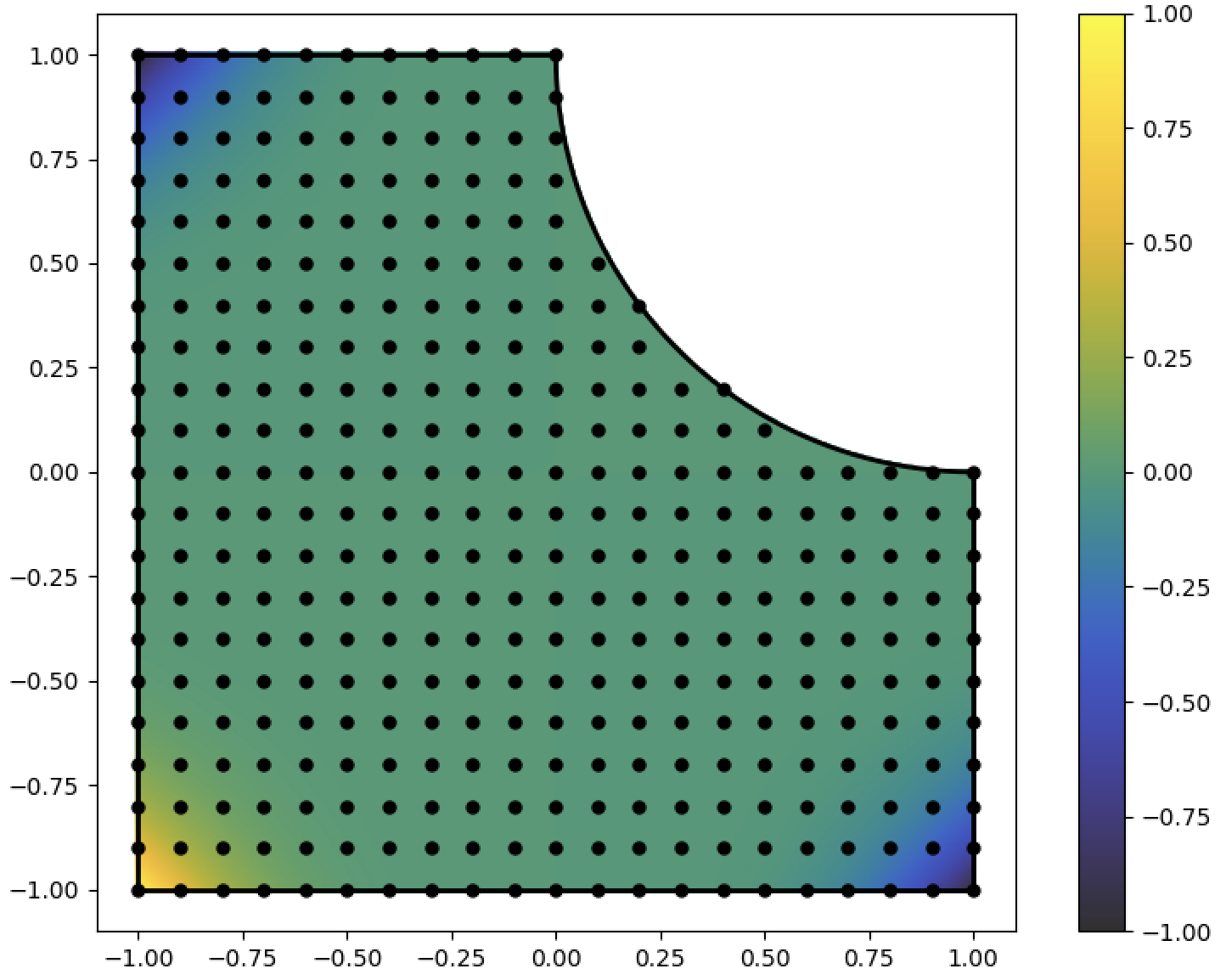}
  \caption{\emph{Left}: A cartoon overview of a domain immersed in a background cutting mesh. The resulting integration domains are either fully interior ($\Omega_{\operatorname{int}}$) or cut ($\Omega_{\operatorname{cut}}$). \emph{Right}: A simple canonical cut cell: the biunit square with a unit disk subtracted from its top right corner. The cut cell fits $[-1, 1] \times [-1, 1]$ tightly. An example of the sort of discretized grid we use to solve to moment-matching equations using sparse optimization is shown: a $21 \times 21$ grid is determined by setting $h = 1/5$; and, after removing nodes which lie strictly outside of the cut cell, 355 nodes remain. Shown in the background is an example of a monomial basis function, $\phi(x, y) = x^5 y^5$, restricted to the cut cell.}\label{fig:cut-cell-diagrams}
\end{figure}

The methods discussed in this section are, at their core, different ways of exploiting the same structure: the feasible moment-matching set is an $m$-dimensional convex cone.
Any quadrature with more than $m$ nodes has internal redundancy, and each method removes that redundancy using a different perspective from convex geometry or optimization.

This section presents these algorithms as a modular toolbox and shows how each fits into this common moment-matching framework:
\begin{enumerate}
\item In \Cref{ssec:discretizing-the-cut-cell}, we explain the spatial discretization used throughout, yielding a uniform grid of candidate nodes contained inside the cut cell.
\item We study in some depth what happens when the quadrature weights are computed using ordinary least squares in \Cref{sec:least-squares-quadrature}. The efficiency of these rules tends to zero for large degree and is of little practical use, but they have important theoretical ramifications.
\item Steinitz elimination is explained in \Cref{sec:steinitz-elimination}. Coupled with the least squares rules presented in \Cref{sec:least-squares-quadrature}, it gives a constructive proof of Tchakaloff's theorem and produces positive rules with efficiency $1/(d+1)$.
\item An alternative proof of Tchakaloff's theorem that furnishes a useful direct algorithm is direct solution of the nonnegative least squares (NNLS) quadrature optimization problem.
  This is discussed in \Cref{sec:nonnegative-least-squares-quadrature}.
  The classic primal active set method due to Lawson and Hanson is both efficient and extremely stable and accurate.
\item A natural model that has been used for quadrature optimization several times in the literature is basis pursuit~\cite{mehrotraGeneratingMomentMatching2013,glaubitzStableHighorderCubature2021}. In \Cref{sec:basis-pursuit-quadrature}, we point out that basis pursuit leads to a degenerate form of Ryu and Boyd's LP and that it likely offers no practical advantage over NNLS.
\item In \Cref{sec:gauss-newton}, we explain how Gauss--Newton can be used to polish nearly optimal quadratures by directly solving the nonlinear moment-matching equations. This arises as a post-processing step in other methods.
\item As discussed in \Cref{sec:ryu-boyd-quadrature}, Ryu and Boyd introduced a semi-infinite LP which yields a generalization of the basis pursuit LP after discretization.
  Their introduction of an arbitrary ``residual function'' allows additional control over the quadrature that is computed.
  Although the choice of this function has no bearing on whether the underlying polynomial space is integrated correctly, selecting it carefully can affect the efficiency of the quadrature, as well as qualitative traits, such as whether the rule is ``Gauss'' (all interior nodes) or ``Lobatto'' (nodes on the boundary allowed).
\item In \Cref{sec:xiao-gimbutas-elimination}, we describe Xiao--Gimbutas elimination, a form of nonlinear node elimination which carefully deletes quadrature nodes without blowing up the residual of the moment-matching equations.
  We explicitly explain how Xiao and Gimbutas's node significances relate to the moment-matching residuals, and demonstrate how it can be used to enhance the efficiency of NNLS quadratures, sometimes quite dramatically.
\item Finally, in \Cref{sec:pruning-inefficient-rules}, we point out that many of these algorithms can be used to ``prune'' known, inefficient composite rules---say, those that have been arrived at by meshing a cut cell. This can be a useful addition to systems that already follow this approach, allowing for a tie-in between the ``native'' approach advocated here and existing alternatives.
\end{enumerate}
Following this, in \Cref{sec:fast-construction}, we synthesize these ideas and present novel fast algorithms for computing cut cell quadratures.

\subsection{Discretizing the cut cell}\label{ssec:discretizing-the-cut-cell}

All quadrature algorithms in this paper begin with the same discretization step: generating a fixed set of candidate nodes inside the cut cell.
This subsection describes how we construct this node set, what geometric information is required, and how we choose the grid resolution.

\paragraph{Computing a tight bounding box and sampling the nodes}
The algorithms presented in this section are unified in the sense that they all start by discretizing $\Omega$ into a uniform grid of $n$ nodes where $n \gg m$, fixing this choice of nodes, and subsequently computing the weights.
The discretization step is intentionally simple: we always use a uniform grid on a bounding box for the cut cell.
See \Cref{fig:cut-cell-diagrams} for the setup.
The geometric challenge lies in supplying two primitives.
We need to be able to compute a tight bounding box of $\Omega$ and to be able to filter nodes that lie outside $\Omega$.
Since these depend entirely on our representation of the geometry, we treat these primitives as black boxes---although, below, we do briefly survey common issues with different geometric formats.
Beyond this, we do not spend any time on the question of how best to discretize $\Omega$, although there are alternatives that may be interesting to investigate in the future~\cite{fornbergFastGeneration2D2015,vandersandeFastVariableDensity2021,duhFastVariableDensity2021,duhDiscretizationNonuniformRational2024}.

\begin{figure}
  \includegraphics[width=\linewidth]{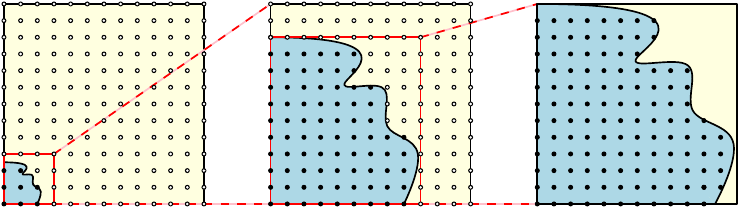}
  \caption{A diagram showing the flow of \Cref{alg:tight-bounding-box}. From left to right, we use the inside-outside test to progressively tighten the fit of an initially coarse bounding box.}\label{fig:bounding-box}
\end{figure}

The first step is computing a tight bounding box.
In the absence of a primitive for computing a tight bounding box, the inside-outside test alone is enough to drive a simple refinement procedure (see \Cref{fig:bounding-box} for a diagram):
\begin{sfpalgorithm}[Compute a tight bounding box using only the inside-outside test.]
  \label{alg:tight-bounding-box}
  \begin{enumerate}[nosep]
  \item Initialize the bounding box to the background mesh element from which $\Omega$ has been cut.
  \item Sample a uniform grid on the current box with spacing $h > 0$.
  \item Evaluate the inside-outside test at each grid node.
  \item Shrink the box to the minimal axis-aligned box containing: 1) all nodes inside or on the boundary of $\Omega$, and 2) all outside nodes that are adjacent to an inside node.
  \item If bounding box changed, go to Step 2; otherwise, terminate.
  \end{enumerate}
\end{sfpalgorithm}
\noindent Even if an algorithm for computing a tight bounding box is available, \Cref{alg:tight-bounding-box} may be worth evaluating for its simplicity and speed.

Once we've fit a bounding box to $\Omega$, we compute our set of candidate nodes and remap the grid to a reference integration domain:
\begin{sfpalgorithm}[Discretize cut cell $\Omega$ with a uniform grid.]
  \label{alg:discretize-cut-cell}
  \begin{enumerate}[nosep]
  \item Fit a tight bounding box to the cut cell, possibly using \Cref{alg:tight-bounding-box} Denote it by $[a_1, b_1] \times \cdots \times [a_d, b_d] \supseteq \Omega$.
  \item Map $[a_1, b_1] \times \cdots \times [a_d, b_d]$ to $[-1, 1]^d$, denote the mapping by $\mF$.
  \item Discretize $[-1, 1]^d$ into a uniform regular grid with spacing $h > 0$.
  \item Let $\mX$ be the set of grid nodes contained in $\mF(\Omega)$.
  \end{enumerate}
\end{sfpalgorithm}
\noindent In the first step, by ``tight'' we mean that we try to make the rectangle $[a_1, b_1] \times \cdots \times [a_d, b_d]$ as small as possible while still containing $\Omega$.
No issue is presented if the boundaries of the bounding box and $\Omega$ overlap.

Following the rescaling done in Step 2 of \Cref{alg:discretize-cut-cell}, we must use the change of variables formula as the starting point for computing moments:
\begin{equation}
  \int_\Omega \phi(\mx)d\mx = \int\limits_{\mF^{-1}(\Omega)} \phi(\mF(\m{u})) |\sfD\mF(\m{u})|d\m{u} = \frac{1}{2^d} \prod_{j=1}^d (b_j - a_j) \int\limits_{\mF^{-1}(\Omega)} \phi(\mF(\m{u}))d\m{u}.
\end{equation}
Rescaling $\Omega$ serves two purposes.
Most obviously, it allows a more even distribution of grid nodes over $\Omega$ which avoids some issues caused by tiny sliver cut cells.
However, it also improves the tendency of our choice of residual functions to create Lobatto\footnote{Colloquially, a ``Lobatto rule'' is a quadrature rule with points on the boundary of the domain.} rules when computing Ryu--Boyd quadratures, as in \Cref{sec:fast-primal-heuristic}.

\paragraph{Implicit vs.\ parametric geometry}
The inside-outside test and possibly a primitive for generating a tight bounding box (alternative to \Cref{alg:tight-bounding-box}) depend entirely on how $\Omega$ is represented:
\begin{itemize}
\item For implicit domains (i.e., defined by level-set functions), inside-outside testing is trivial: just evaluate the level set function.
Computing a tight bounding box is more involved and may require either solving small constrained optimization problems or using an isosurfacing routine to approximate the boundary.
Basic methods like marching cubes~\cite{lorensenMarchingCubesHigh1998} are simple but incur geometric and topological errors unless augmented with higher-order meshing and interval‐arithmetic techniques (see \Cref{sec:comparison})~\cite{plantingaIsotopicApproximationImplicit,snyderGenerativeModelingComputer2014}).
\item For parametric or B-rep geometries~\cite{stroudBoundaryRepresentationModelling2006}, the situation is inverted.
Computing a bounding box is easy, but inside-outside testing is notoriously delicate.
B-reps are rarely watertight, which undermines ray-tracing predicates.
This has motivated robust alternatives based on layer‐potential integrals and fast winding numbers~\cite{jacobsonRobustInsideoutsideSegmentation2013,barillFastWindingNumbers2018,spainhourRobustContainmentQueries2024}.
Simpler closest-point predicates are also possible when normal orientation is reliable.
\end{itemize}

Again, we stress that we assume that the solutions to these problems have been provided as black boxes ahead of time, although they are \emph{not} trivial.
The way these problems should be solved depends completely on user requirements and the underlying geometry.
Different approaches are discussed extensively in the literature, but the details are scattered across many disparate sources.

\paragraph{Choosing the grid resolution}
A natural question is how finely one should sample the uniform grid in order to resolve all potential optimal quadratures on $\Omega$.
Following Huybrechs, a simple uniform discretization is a perfectly reasonable choice on a one-dimensional interval~\cite{huybrechsStableHighorderQuadrature2009}.
How well this works for a cut cell appears to be largely independent of the cut cell's geometry, provided that the grid is fine enough.
Martinsson et al.\ give strong evidence that this should indeed be the case~\cite{martinssonInterpolationIntegrationFinitedimensional2007}.

To choose a sufficiently fine grid, Huybrechs' prescription is that the number of nodes needed to integrate $\bbP_N^1$ is $O(N^2)$ with a small constant (see \Cref{sec:huybrechs}).
By extension, to integrate $\bbP_N^d$, we choose $h = O(N^{-{2d}})$, resulting in a grid with $O(N^{2d})$ nodes, albeit with an extremely small constant.

We now give a heuristic argument that this choice is essentially necessary for a uniform grid.
Recall that the density of quadrature points on $[-1, 1]$ approaches $\rho(x) = \frac{1}{\sqrt{1 - x^2}}$ as $n \to \infty$.
Consider the Chebyshev nodes $x_j = -\cos(\pi j/n)$ for $j = 0, 1, \hdots, n$.
Since the nodes cluster most strongly near the endpoints, it suffices to examine the spacing between $x_0$ and $x_1$.
Thus, the minimal spacing is $\frac{\pi^2}{2n^2} - O(n^{-4})$, and a uniform grid can resolve this only if $h = O(n^{-2})$.

\begin{figure}
  \centering
  \includegraphics[width=0.48\linewidth]{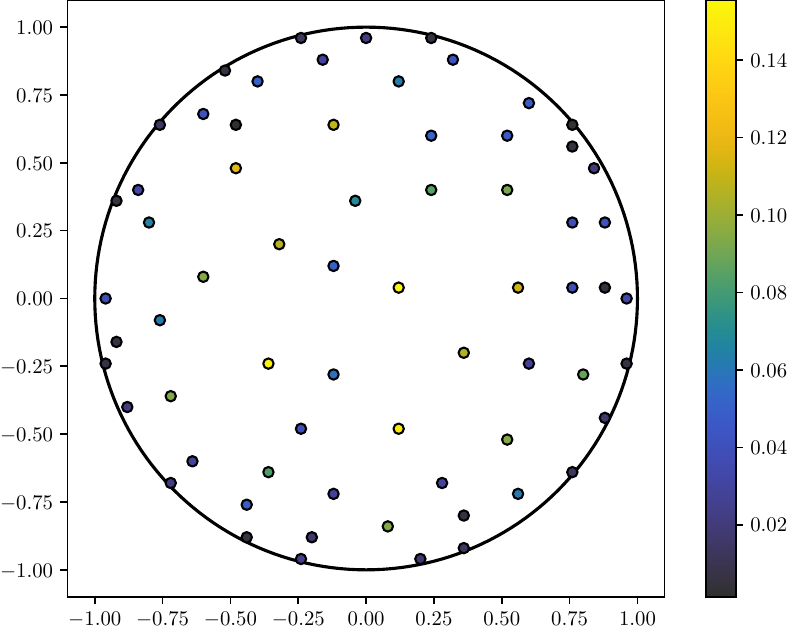}%
  \hspace{0.025\linewidth}%
  \includegraphics[width=0.48\linewidth]{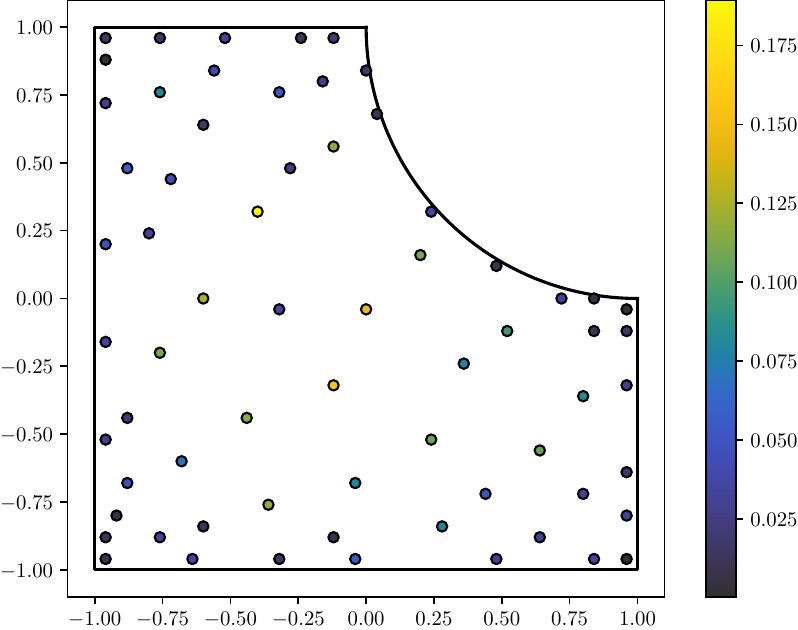}
  \caption{Two simple examples of NNLS quadratures computed using the Lawson--Hanson algorithm, starting from a uniform grid restricted to $\Omega$. \emph{Left}: $\Omega$ is the unit disk. \emph{Right}: $\Omega$ is the biunit square with a unit disk centered at $(1, 1)$ removed. Both quadratures integrate $\mathbb{P}_{10}^2$ (up to numerical roundoff) and have exactly $\dim{\mathbb{P}_{10}^2} = 66$ points, as expected. The moment-matching equations are solved directly in the monomial basis. In both cases, $\mV$ has 66 columns and $\sim$2,000 rows.}\label{fig:simple-nnls-2d}
\end{figure}

How this generalizes to higher dimensions is less clear.
It may sometimes be possible to ``beat'' this requirement on $h$ if nodes do not cluster along a one-dimensional boundary feature (e.g., consider a disk, where this does not happen; see \Cref{fig:simple-nnls-2d}).
For domains with boundary singularities or cusps, the behavior can be worse.
In two and three dimensions, our conclusion is that even if $h = O(n^{-2})$ is necessary in the worst case, the hidden constant is very small.

Since this constant is unknown in advance, a practical strategy is simply to refine the grid until a method succeeds.
As we show in \Cref{thm:steinitz-grid-doubling-complexity}, increasing the grid resolution until a quadrature is successfully computed has no adverse effect on the overall complexity of our algorithms.

\subsection{Least squares quadrature}\label{sec:least-squares-quadrature}

\begin{figure}
  \includegraphics[width=\linewidth]{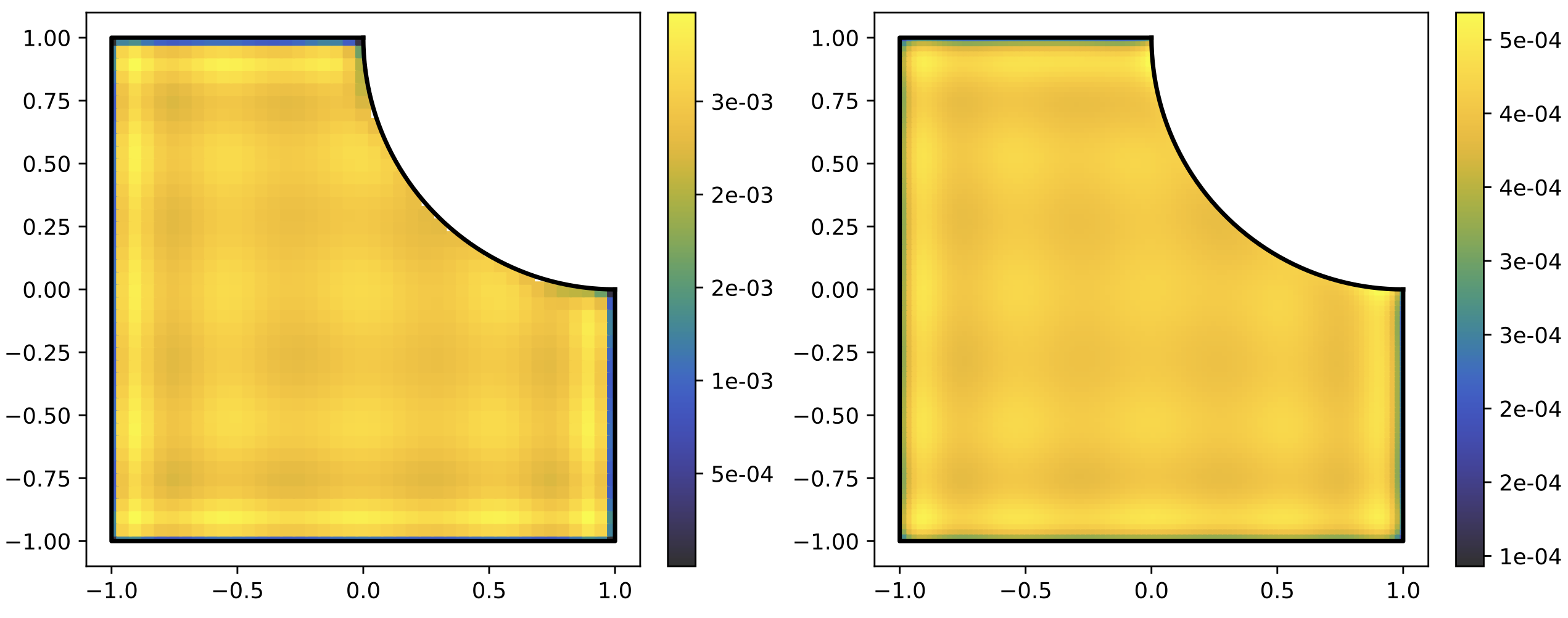}
  \caption{Examples of least squares quadrature rules computed following \Cref{sec:least-squares-quadrature}. Each plot shows a least squares rule integrating $\bbP_{10}^2$. Both rules use uniform grids discretizing the domain $\Omega$ according to \Cref{ssec:discretizing-the-cut-cell}. \emph{Left}: the least squares rule for the coarsest grid yielding a positive quadrature (in this case, the grid has 41 nodes along each axis discretizing the bounding box of $\Omega$). \emph{Right}: another least squares rule, this time with 100 nodes along each axis. Considering these rules, it is intuitively clear that the weights of the least squares rule converge to samples of a continuous density.}\label{fig:least-squares-quadrature}
\end{figure}

Least squares quadratures offer a clean and inexpensive baseline.
They solve the moment-matching equations in closed form and reveal much of the asymptotic structure of the problem.
Unfortunately, they are also extremely dense and therefore inefficient in practice.
Their primary role for us is as a theoretical reference point---useful for understanding asymptotics---and as a start guess for methods discussed later (see \Cref{sec:steinitz-elimination} and \Cref{sec:fast-primal-heuristic}).

As we have seen, Huybrechs~\cite{huybrechsStableHighorderQuadrature2009} and Glaubitz~\cite{glaubitzConstructingPositiveInterpolatory2020} explored the conditions under which solving rectangular moment-matching equations with fixed nodes produces positive quadratures.
With the nodes fixed, the optimization problem to solve is:
\begin{equation}
  \mbox{minimize} \quad \frac{1}{2} \norm{\mV^\top\mw - \mI}_2^2. \label{eq:least-squares}
\end{equation}
Since $n \geq m$, the basis matrix $\mV \in \bbR^{n \times m}$ is tall or square.
Usually, we have $n \gg m$, so we generally think of this as an underdetermined least squares problem.

For an underdetermined problem, the solution to \Cref{eq:least-squares} is nonunique.
The minimum $\ell_2$ norm solution is easily computed using either the QR or LU decomposition~\cite{golubMatrixComputations2013b}.
\begin{sfpalgorithm}[Compute least squares quadrature.]
  \label{alg:least-squares}
  \noindent \emph{Input}: domain $\Omega \subseteq \bbR^d$, polynomial space $\mathbb{P}$ with $\dim\mathbb{P}=m$, moment vector $\mI$, discretization of $\Omega$ into uniform grid $\mX$ containing $O(m^2)$ points. \\
  \emph{Output}: quadrature rule $(\mX, \mw)$ that integrates $\mathbb{P}$ on $\Omega$.
  \begin{enumerate}[nosep]
  \item Compute $\mw = (\mV^\top)^\dagger \mI$ using LU or QR decomposition.
  \end{enumerate}
\end{sfpalgorithm}
\noindent Above, $(\mV^\top)^\dagger$ denotes the pseudoinverse of $\mV^{\top}$.
Since $\mV$ has full column rank, $(\mV^{\top})^\dagger = \mV (\mV^\top \mV)^{-1}$.
To use the LU decomposition to compute the minimum norm solution of \eqref{eq:least-squares}, we compute the LU factorization of $\mV^\top \mV$.
To apply the QR decomposition, we let $\mV = \mQ\mR$ be the reduced QR decomposition of $\mV$.
We briefly remind the reader of the complexity of this procedure.
\begin{theorem}\label{thm:least-squares-complexity}
  Assuming that $n = O(m^{2d})$, \Cref{alg:least-squares} runs in $O(m^2 n) = O(m^{2d+2})$ time.
\end{theorem}
\begin{proof}
  The complexity follows immediately from the cost of forming $\mV^\top \mV$ and applying either LU or QR; see standard complexity estimates~\cite{golubMatrixComputations2013b}.
\end{proof}
Whether or not this can be improved upon is an interesting question.
Huybrechs devised an asymptotically faster method of constructing least squares quadratures on the interval using properties of discrete orthogonal polynomials~\cite{huybrechsStableHighorderQuadrature2009}.
It may be possible to devise an analog of his algorithm for quadratures in two and higher dimensions using known properties of multivariate orthogonal polynomials~\cite{olverFastAlgorithmsUsing2020}.

The drawback of least squares quadrature is that generally $\mw_j \neq 0$ for all $j$.
Although this is far too dense to be of practical use, these rules play a central theoretical role: they provide the first concrete demonstration that positive quadratures exist for a sufficiently fine discretization.
Huybrechs proved a version of the following theorem in one dimension~\cite{huybrechsStableHighorderQuadrature2009}.
The version here is a specialization of a result due to Glaubitz~\cite{glaubitzStableHighorderCubature2021}.
The strategy is to rewrite the least squares solution using discrete orthogonality to expose positivity.

\begin{theorem}\label{thm:positive-least-squares-quadrature}
  Let $\bbP$ be a finite-dimensional polynomial space defined on $\bbR^d$ containing constants.
  Then, there is an integer $n_0 > 0$ such that for all $n \geq n_0$ a set of points $\mx_1, \hdots, \mx_n \subseteq \Omega \subseteq [-1, 1]^d$ exists such that the resulting minimum norm least squares quadrature is positive.
\end{theorem}

\begin{proof}
  First, let $N$ be a positive integer, and let $h = 1/N$.
  Define the scalars $-1 = t_0 < t_1 < \cdots t_N = 1$ so that $t_{k+1} - t_k = h$ for $k = 1, \hdots, N$ and consider the cubes $R_{k_1,\hdots,k_d} = [t_{k_1-1}, t_{k_1}] \times \cdots \times [t_{k_d-1}, \hdots, t_{k_d}]$, where $0 < k_\ell \leq N$ for $\ell = 1, \hdots, d$.
  Defining the multi-index $\m{k} = (k_1, \hdots, k_d)$, we denote the centroid of $R_{\m{k}}$ by $\mx_{\m{k}}$.
  We consider two subsets of these cubes:
  \begin{itemize}[nosep]
  \item \emph{The interior cubes}: $\mathcal{R}_{\operatorname{int}} = \{R_{\m{k}} : R_{\m{k}} \subseteq \Omega \wedge R_{\m{k}} \cap (\partial\Omega \backslash \partial [-1, 1]^d) = \emptyset\}$ (i.e., the cubes which lie inside $\Omega$ but which do not overlap the part of its boundary interior to $[-1, 1]^d$).
  \item \emph{The cut cubes}: $\mathcal{R}_{\operatorname{cut}} = \{R_{\m{k}} : R_{\m{k}} \cap (\partial\Omega \backslash \partial [-1, 1]^d) \neq \emptyset\}$.
  \end{itemize}
  The families $\mathcal{R}_{\operatorname{int}}$ and $\mathcal{R}_{\operatorname{cut}}$ partition $\Omega$.
  To this end, we define:
  \begin{equation}
    \Omega_{\operatorname{int}} = \bigcup_{R_{\m{k}} \in \mathcal{R}_{\operatorname{int}}} R_{\m{k}}, \qquad \Omega_{\operatorname{cut}} = \bigcup_{R_{\m{k}} \in \mathcal{R}_{\operatorname{cut}}} R_{\m{k}},
  \end{equation}
  so that $\Omega = \Omega_{\operatorname{int}} \cup \Omega_{\operatorname{cut}}$.

  Let $f : \Omega \to \bbR$ be a smooth function.
  Recall that if $\|f_{\malpha}\|_{L^\infty(\Omega)} = O(1)$ for $\malpha$ with $|\malpha| = 2$, the error in the midpoint rule applied on a single cube $R_{\m{k}} \in \mathcal{R}_{\operatorname{int}}$ is:
  \begin{equation}
    \left|\int_{R_{\m{k}}} f(\mx)d\mx - h^d f(\mx_{\m{k}})\right| = O(h^{d+2}).
  \end{equation}
  Since $\Omega$ is regular, $|\mathcal{R}_{\operatorname{int}}| = O(N^d) = O(h^{-d})$, so that the composite midpoint rule on $\Omega_{\operatorname{int}}$ satisfies:
  \begin{equation}
    \left|\int_{\Omega_{\operatorname{int}}} f(\mx)d\mx - \sum_{R_{\m{k}} \in \mathcal{R}_{\operatorname{int}}} h^d f(\mx_{\m{k}})\right| = |\mathcal{R}_{\operatorname{int}}| O(h^{d+2}) = O(h^2).
  \end{equation}
  Again, since $\Omega$ is regular, we also have $|\mathcal{R}_{\operatorname{cut}}| = O(N^{d-1}) = O(h^{1-d})$.
  Consequently, we can bound:
  \begin{equation}
    \begin{split}
      &\left|\int_{\Omega_{\operatorname{cut}}} f(\mx)d\mx\right| \leq \sum_{R_{\m{k}}\in\mathcal{R}_{\operatorname{cut}}} \left|\int_{R_{\m{k}}} f(\mx)d\mx\right| \\
      &\hspace{8em} \leq h^d  \sum_{R_{\m{k}}\in\mathcal{R}_{\operatorname{cut}}} \underbrace{\max_{\mx \in R_{\m{k}}} |f(\mx)|}_{O(1)} = O(h^d) O(h^{1-d}) = O(h). \label{eq:tmp1}
    \end{split}
  \end{equation}
  This lets us bound:
  \begin{equation}
    \left|\int_{\Omega} f(\mx)d\mx - \sum_{R_{\m{k}}\in\mathcal{R}_{\operatorname{int}}} h^d f(\mx_{\m{k}})\right| = \underbrace{\left|\int_{\Omega_{\operatorname{cut}}} f(\mx)d\mx\right|}_{O(h)} + O(h^2) = O(h).
  \end{equation}
  In summary, the result of applying the composite midpoint rule for $\Omega_{\operatorname{int}}$ to a smooth enough function $f$ (ignoring the part of $f$ supported on $\Omega_{\operatorname{cut}}$) is $O(h)$ accurate if $\Omega$ is regular.

  Next, with $n_{\operatorname{int}} = |\mathcal{R}_{\operatorname{int}}|$, we let $\mX_h = (\mx_{\m{k}})_{R_{\m{k}}\in\mathcal{R}_{\operatorname{int}}} \in \bbR^{n_{\operatorname{int}}}$ and define the following discrete inner product:
  \begin{equation}
    (f, g)_{\mX_h} = \sum_{R_{\m{k}}\in\mathcal{R}_{\operatorname{int}}} f(\mx_{\m{k}})g(\mx_{\m{k}}).
  \end{equation}
  For our target polynomial space, we let $\{\phi_i\}_{i=1}^m$ be a basis satisfying $(\phi_i, \phi_j)_{\mX_h} = \delta_{ij}$ for each $i$ and $j$, and such that $\phi_1(\mx) \equiv n_{\operatorname{int}}^{-1/2}$.
  Using the notation:
  \begin{equation}
    f(\mX_h) = \big(f(\mx_{\m{k}})\big)_{R_{\m{k}}\in\mathcal{R}_{\operatorname{int}}} \in \bbR^{n_{\operatorname{int}}}
  \end{equation}
  we define the Vandermonde matrix $\mV \in \bbR^{n_{\operatorname{int}} \times m}$ by:
  \begin{equation}
    \mV = \begin{bmatrix} \phi_1(\mX_h) & \cdots & \phi_m(\mX_h) \end{bmatrix}.
  \end{equation}
  We also define the usual moment vector $\mI \in \bbR^m$ by $\mI_i = \int_\Omega \phi_i(\mx)d\mx$.
  The discrete orthogonality of our basis gives $(\mV{}^\top \mV)_{ij} = \delta_{ij}$ for each $i$ and $j$.
  Consequently, from $\mw = \mV (\mV^\top \mV)^{-1} \mI$, we have:
  \begin{equation}
    \mw = \mV \mI = \sum_{i=1}^m \left(\int_\Omega\phi_i(\mx)d\mx\right) \phi_i(\mX_h).
  \end{equation}
  The first term in this sum is given by:
  \begin{equation}
    \left(\int_\Omega\phi_1(\mx)d\mx\right) \phi_1(\mX_h) = \frac{\operatorname{vol}(\Omega)}{n_{\operatorname{int}}^{1/2}} \m{1}_{n_{\operatorname{int}}}.
  \end{equation}
  We also have:
  \begin{equation}
    \int_\Omega\phi_i(\mx)d\mx = \sum_{R_{\m{k}}\in\mathcal{R}_{\m{k}}} h^d \phi_i(\mx_{\m{k}}) + O(h)  = n^{1/2}_{\operatorname{int}} h^d \underbrace{\big(\phi_i, \phi_1\big)_{\mX_h}}_{=0} + O(h) = O(h).
  \end{equation}
  Since $\|\phi_i\|_{\mX_h}^2 = 1$ for each $i$, we have $|\phi_i(\mx_{\m{k}})| \leq 1$ for each $\m{k}$.
  Altogether, we conclude:
  \begin{equation}
    \left\|\sum_{i=2}^m \phi_i(\mX_h) \int_\Omega\phi_i(\mx)d\mx\right\|_\infty \leq \sum_{i=2}^m \underbrace{\left\|\phi_i(\mX_h)\right\|_\infty}_{=O(1)} \underbrace{\left|\int_\Omega \phi_i(\mx)d\mx\right|}_{=O(h)} = O(h).
  \end{equation}
  
\end{proof}

Although least squares quadratures are dense and inefficient, their positivity and limiting behavior make them a natural starting point for sparse methods. In the next section we show how Steinitz elimination can be used to reduce least squares rules to positive and efficient quadratures.

\subsection{Steinitz elimination}\label{sec:steinitz-elimination}

\begin{figure}
  \includegraphics[width=\linewidth]{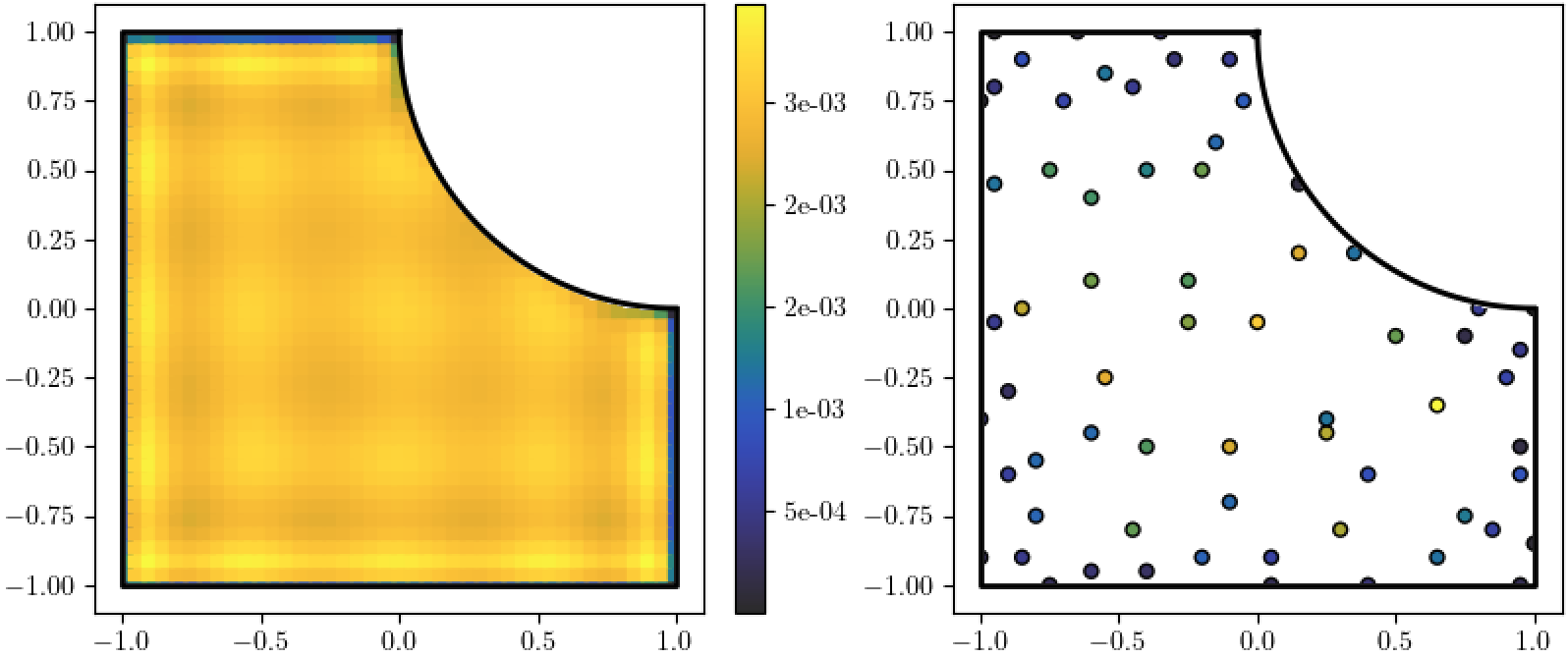}
  \caption{An example of a reduced quadrature computed using \Cref{alg:steinitz-grid-doubling}. \emph{Left}: the first positive least squares quadrature found after adaptively increasing the grid resolution. \emph{Right}: the order $m$ quadrature computed using Steinitz elimination (\Cref{alg:steinitz}).}\label{fig:steinitz-grid-doubling}
\end{figure}

Steinitz elimination is a constructive procedure for reducing an inefficient, positive quadrature to a positive quadrature with a minimal number of nodes, where this number matches what is predicted by Tchakaloff's theorem (\Cref{thm:tchakaloff}).
Specfically, given a quadrature with $n > m$ nodes, the moment-matching constraints define an affine subspace of $\mathbb{R}^n$ and the positivity constraints further restrict the set of possible weights to a convex polytope that lies in the nonnegative orthant in $\bbR^n$.
Steinitz's idea is simple: move the weight vector along a direction which leaves all moments unchanged until one weight hits zero, remove that node, and repeat the process in $\bbR^{n-1}$.
After at most $n - m$ such eliminations, one arrives at an $m$-node positive quadrature, achieving exactly the bound guaranteed by Tchakaloff.

Starting from a suboptimal but positive least squares quadrature, where $w_i > 0$ for each $i$, we can use Steinitz elimination to reduce this quadrature to one which is order $m$.
Pruning a positive least squares quadrature in this way was first investgated by Davis \Cref{thm:caratheodory-quadrature}~\cite{davisConstructionNonnegativeApproximate1967} and was later revisited by Glaubitz~\cite{glaubitzConstructingPositiveInterpolatory2020}.

Glaubitz provides a convenient proof of this result~\cite{glaubitzConstructingPositiveInterpolatory2020}; we also recapitulate this result here.

\begin{theorem}\label{thm:steinitz-elimination}
  Let $d > 0$, $N > 0$, and $m = \dim\bbP_N^d$.
  Let $\mathcal{Q}$ be a positive quadrature with $n \geq m$ nodes that integrates $\bbP_N^d$.
  Then, there exists a positive quadrature supported on exactly $m$ nodes that integrates $\bbP_N^d$.
  The nodes of the reduced quadrature are a subset of the original quadrature's nodes.
\end{theorem}

\begin{proof}
  Let $\mX = (\mx_1, \hdots, \mx_n)$ and $\mw = (w_1, \hdots, w_n)$ denote the nodes and weights of the order $n$ positive quadrature, so that:
  \begin{equation}
    \mV^\top \mw = \mI
  \end{equation}
  holds, where $\mV \in \bbR^{n \times m}$ is a matrix with $\mV_{ij} = \phi_j(\mx_i)$ for some basis $\{\phi_j\}$ of $\bbP_N^d$.

  Now, let $\mathcal{I} = \{\, i : w_i > 0\,\}$ and let $\mV_{\mathcal{I}}$ be the submatrix of $\mV$ containing only those rows.
  If $|\mathcal{I}| = m$, we're done.
  Otherwise, $\mV_{\mathcal{I}}^\top$ has a nontrivial nullspace.
  Without loss of generality, let $\delta\mw\in\bbR^n$ be a vector with at least one negative component, and which satisfies:
  \begin{enumerate}
  \item $\delta\mw_i = 0$ for each $i \notin \mathcal{I}$, and
  \item $\delta\mw_{\mathcal{I}} \in \operatorname{null}(\mV_{\mathcal{I}}^\top)$,
  \end{enumerate}
  By this choice, $\mV^\top\delta\mw = 0$ and shifting $\mw$ by a scalar multiple of $\delta\mw$ preserves all moments.

  For each index with $\delta\mw_i < 0$, define:
  \begin{equation}
    \alpha_i = -\frac{w_i}{\delta\mw_i} > 0,
  \end{equation}
  and let $\alpha_* = \min_{\{i \,:\, \delta\mw_i < 0\}} \alpha_i > 0$.
  By construction, $\mw + \alpha_* \delta\mw \geq 0$ and at least one index---call it $i_*$---satisfies $w_{i_*} + \alpha_* \delta\mw_{i_*} = 0$.
  Thus, setting $\mw \gets \mw + \alpha_* \delta\mw$, at least one positive weight becomes zero without violating the moment constraints, and the number of nonzero components strictly decreases.

  Repeating this argument at most $n - m$ times, we reduce the original order $n$ quadrature to a quadrature that is order $m$.
\end{proof}

Some care is needed if we use \Cref{thm:steinitz-elimination} to reduce a positive least squares quadrature, since the initial order of the quadrature is $n \gg m$.
The most costly part of the algorithm outlined in the proof of \Cref{thm:steinitz-elimination} is the computation of the nullspace.
Rather that compute the nullspace of the matrix $\mV_{\mathcal{I}}^\top$, where $\mathcal{I}$ consists of the indices corresponding to positive weights in our quadrature, it is enough to take a subset of $m + 1$ indices of $\mathcal{I}$.
This is enough to guarantee that there will always be a nontrivial nullspace---the fact that this procedure works is ultimately guaranteed to us by Carath\'{e}odory's theorem (see \Cref{thm:caratheodory-quadrature}).
This can be formalized as an algorithm as follows:
\begin{sfpalgorithm}[Prune a positive quadrature using Steinitz elimination.]
  \label{alg:steinitz}
  \noindent\emph{Input}: nodes $\mX = (\mx_1,\ldots,\mx_n)$ and weights $\mw \ge 0$ defining a quadrature that integrates $\bbP_N^d$, and a moment vector $\mI\in\bbR^m$. \\
  \noindent\emph{Output}: an order $m$ positive quadrature.
  \begin{enumerate}[nosep]
  \item Let $\mathcal{I}$ consist of the first $m + 1$ indices for which $w_i > 0$, so that $|\mathcal{I}| \leq m + 1$.
  \item If $|\mathcal{I}| = m$, stop.
  \item Form $\mV_{\mathcal{I}} \in \mathbb{R}^{(m+1)\times m}$ consisting of the rows of $\mV$ indexed by $\mathcal{I}$.
  \item Compute a nonzero vector $\delta\mw_{\mathcal{I}} \in \operatorname{null}(\mV_{\mathcal{I}}^\top)$ and extend it to $\delta\mw \in \mathbb{R}^n$ by setting $\delta\mw_i = 0$ for $i\notin\mathcal{J}$.
  \item Flip the sign of $\delta\mw$ if necessary so that at least one component is negative.
  \item Set $\mw \gets \mw + \alpha_* \delta\mw$, where  $\alpha_* = \min_{\{i:,\delta\mw_i < 0\}} \alpha_i$.
  \item Return to Step 1.
  \end{enumerate}
\end{sfpalgorithm}

Once we cap the number of rows in $\mV_{\mathcal{I}}^\top$ at each step, we reduce the size of the nullspace computation from one involving potentially $O(n)$ rows to one that only depends on $m+1$ rows, without changing the information available to the algorithm.

\begin{theorem}\label{thm:steinitz-complexity}
  \Cref{alg:steinitz} reduces an order $O(n)$ quadrature to an order $m$ quadrature in $O(m^3 n)$ operations.
\end{theorem}

\begin{proof}
  Initially finding and subsequently updating $\mathcal{I}$ requires $O(n)$ operations amortized over all elimination steps.
  The bulk of the work at each step is in computing a null vector for an $(m+1)\times m$ matrix.
  This costs $O(m^3)$ flops, following standard complexity estimates in numerical linear algebra~\cite{golubMatrixComputations2013b}.
  The remaining work in each step---e.g., forming $\alpha_i$, updating $\mw$---is $O(m)$, and therefore negligible.
  Each Steinitz step eliminates at least one positive weight, so there are at most $O(n)$ total steps.
  This gives a complexity of $O(m^3 n)$ overall.
\end{proof}
\noindent It is unclear whether it is possible to reduce this complexity further.
We will see in \Cref{sec:nonnegative-least-squares-quadrature} that this complexity can be beat by updating and downdating a QR decomposition when computing nonnegative least squares rules.
But it does not seem possible to reuse information between different blocks in order to reduce the complexity of computing null vectors in \Cref{alg:steinitz}.

We now show how to combine an adaptive procedure for computing a positive least squares quadrature on a minimally coarse grid, followed by its subsequent reduction using \Cref{alg:steinitz}.
See \Cref{fig:steinitz-grid-doubling} for an example.

\begin{sfpalgorithm}[Adaptively compute a least squares quadrature and subsequently prune it using Steinitz elimination.]
  \label{alg:steinitz-grid-doubling}
  \noindent\emph{Input}: a basis $\{\phi_j\}$ for the polynomial space $\bbP_N^d$, the moment vector $\mI$, and an initial grid spacing $h = 2^{1 - p}$ with $p = O(1)$. \\
  \noindent\emph{Output}: an order $m$ positive quadrature.
  \begin{enumerate}[nosep]
  \item Compute a tight bounding box of $\Omega$ using \Cref{alg:tight-bounding-box}.
  \item Run \Cref{alg:least-squares} to compute an initial least squares quadrature, with $\mw_{\operatorname{LS}}$ the corresponding weight vector.
  \item If the residual of the least squares system is too large, or if any component of $\mw_{\operatorname{LS}}$ is negative, set $h = h/2$ and go to Step 2.
  \item Run Steinitz elimination following \Cref{alg:steinitz} to reduce the quadrature to order $m$.
  \end{enumerate}
\end{sfpalgorithm}

\begin{theorem}\label{thm:steinitz-grid-doubling-complexity}
  Let $n$ be the number of nodes in $\Omega$ the final time Step 3 of \Cref{alg:steinitz-grid-doubling} is executed.
  Then, the complexity of \Cref{alg:steinitz-grid-doubling} is $O(m^3 n)$.
\end{theorem}

\begin{proof}
  Let $h_1 = 1$ and let $h_j = 2^{1-j}$ be the grid spacing after the $j$th time of Step 3 is run; let $n_j$ be the corresponding number of nodes in $\Omega$.
  Halving $h$ increases the number of nodes in $\Omega$ by a factor of $2^d$, up to a constant, which gives $n_{j+1} = O(2^d n_j)$.
  If $J$ is the final value of $j$, then this also gives $n_j = O(2^{-(J-j)d}n_J)$, where $n_J$ is the number of nodes in $\Omega$ after the final time Step 3 executes.

  At grid level $j$, the least squares solve is $O(m^2 n_j)$ by \Cref{thm:least-squares-complexity}.
  Summing over all grids gives:
  \begin{equation}
    \sum_{j=1}^J O(m^2 n_j) = O\left(m^2 n_J \sum_{j=1}^J 2^{-(J-j)d}\right) = O(m^2 n_J).
  \end{equation}
  Thus, the entire least squares phase costs $O(m^2 n)$, where $n = n_J$.
  By \Cref{thm:steinitz-complexity}, Steinitz elimination costs $O(m^3 n)$, which dominates.
\end{proof}
In summary, grid refinement followed by Steinitz elimination yields an $m$ point positive quadrature without incurring excessive overhead from the refinement phase.
This adaptive scheme therefore provides an efficient and robust baseline for all subsequent sparse methods.

\subsection{Nonnegative least squares quadrature}\label{sec:nonnegative-least-squares-quadrature}

The theorems in \Cref{sec:existence-of-quadratures} show that if $n$ is large enough, then an order $m$ positive quadrature exists.
Since these proofs are constructive, they could be used as the basis of an algorithm for computing an efficient positive quadrature.
However, what this construction suggests is that given a ``dense enough'' point set, we should be able to select a subset of points which together yield an efficient positive quadrature.
In the language of optimization, we could consider the nonnegative least squares problem:
\begin{equation}
  \begin{split}
    \label{eq:nnls}
    \mbox{minimize}\quad & \frac{1}{2} \norm{\mV^\top\mw - \mI}_2^2 \\
    \mbox{subject to}\quad & \mw \geq 0,
  \end{split}
\end{equation}
whose global optima constitute the polytope $\{\mw\in\bbR^n : \mV^\top \mw = \mI \mbox{ and } \mw \geq 0\}$.
Since there is no unique minimizer, the choice of algorithm used to solve \eqref{eq:nnls} will bias the selection.
Again, if $n$ is large enough, ordinary least squares produces the minimum $\ell_2$ norm solution (corresponding to a point somewhere in the interior of the feasible set---more on this later when we turn to linear programming), but \Cref{thm:positive-least-squares-quadrature} indicates that every component will be nonzero, forcing us to turn to the relatively inefficient Steinitz elimination process described in \Cref{sec:steinitz-elimination}.

A suitable algorithm for solving \eqref{eq:nnls} is the active set method due to Lawson and Hanson~\cite{lawsonSolvingLeastSquares1995}.
To derive this algorithm, define the Lagrangian of \eqref{eq:nnls}:
\begin{equation}
  L(\mw, \malpha) = \frac{1}{2} \left\|\mV^\top\mw - \mI\right\|_2^2 + \malpha^\top \mw.
\end{equation}
The Karush--Kuhn--Tucker (KKT) conditions for a constrained local optimum are:
\begin{equation}
  L_{\mw}(\mw, \malpha) = \mV\big(\mV^\top\mw - \mI\big) + \malpha = \m{0}, \quad \malpha^\top \mw = 0, \quad \mw \geq 0, \quad \malpha \leq 0.
\end{equation}
In rough outline, an active set method for solving \eqref{eq:nnls} chooses an initial iterate $\mw_0$ and determines the inactive set $J_0 = \Set{i : 1 \leq i \leq n \mbox{ and } (\mw_0)_i > 0}$.
After computing $\malpha_0 = \mV(\mI - \mV^\top\mw_0)$, we determine the maximum component of $\malpha_0$, remove it from $\calI_0$, and compute $\mw_1$ by updating the new inactive component of $\mw_0$.
Continuing in this way, $\mw_k$ and $\malpha_k$ are determined for $k > 1$, until all Lagrange multipliers are nonpositive.
At each step, a small least squares problem is solved to compute the nonzero components of $\mw_k$.
Occasionally, it is necessary to backtrack to ensure that the constraints are satisfied exactly for each iteration.

In an efficient implementation of the Lawson--Hanson algorithm, since only one component is added or removed from the active set at a time, we can maintain a running QR decomposition, using standard methods for updating and downdating the factorization one column at a time as we go~\cite{golubMatrixComputations2013b}.
These details are summarized in Lawson and Hanson's classic monograph~\cite{lawsonSolvingLeastSquares1995}.
We formalize this in an algorithm.
\begin{sfpalgorithm}[Compute quadrature using Lawson--Hanson NNLS algorithm.]
  \label{alg:lawson-hanson-nnls}
  \noindent\emph{Input}: grid nodes $\mx_1, \hdots, \mx_n \in \Omega$, basis matrix $\mV \in \bbR^{n \times m}$, moment vector $\mI \in \bbR^m$. \\
  \noindent\emph{Output}: order $m$ quadrature which minimizes \eqref{eq:nnls}.
  \begin{enumerate}[nosep]
  \item Initialize $\mw = \m{0}$, compute the active set $\mathcal{I} = \{1, \hdots, n\}$.
  \item Denote the inactive set by $\mathcal{J} = \mathcal{I}^c$.
  \item Set up the running QR decomposition. Letting $\mQ$ be the running $\mQ$ factor, we store $\mQ^\top \mV^\top$ and $\mQ^\top \mI$, and update these as we go. The $\mR$ factor is stored in the columns of $\mQ^\top\mV^\top$ indexed by $\mathcal{J}$.
  \item Compute Lagrange multipliers $\malpha = \mV(\mI - \mV_{\mathcal{J}}^\top\mw_{\mathcal{J}})$. (When $\mathcal{J} = \emptyset$, we interpret ``$\mV_{\mathcal{J}}^\top \mw_{\mathcal{J}}$'' as the zero vector in $\bbR^m$.)
  \item If $\max_i \alpha_i \le 0$, we're done: $\mw$ is optimal.
  \item Set $k = \Arg\max_{i \in \mathcal{I}} \alpha_i$, remove $k$ from $\mathcal{I}$, and update the running QR factorization to include $k$th column of $\mV^\top$.
  \item Define $\hat{\mw}$ by setting $\hat{\mw}_{\mathcal{I}} = \m{0}$ and $\hat{\mw}_{\mathcal{J}} = (\mV_{\mathcal{J}}^\top)^\dagger \mI$. The latter is computed using our running QR decomposition via an upper triangular with $\mR$.
  \item If $\min_i\hat{\mw}_i \geq 0$, accept the step and set $\mw \gets \hat{\mw}$.
    Otherwise, backtrack along the line segment from $\mw$ to $\hat{\mw}$ until a component in $\hat{\mw}_\mathcal{J}$ hits zero---add that index to $\mathcal{I}$, downdate the QR decomposition, and return to Step 6.
  \item Go to Step 3.
  \end{enumerate}
\end{sfpalgorithm}
\noindent Steps 3--8 in \Cref{alg:lawson-hanson-nnls} constitute one ``iteration'' of the Lawson--Hanson algorithm.
In practice, we have found that for the NNLS problems under consideration, \Cref{alg:lawson-hanson-nnls} always terminates in at most $O(m)$ iterations.
Since Lawson--Hanson is essentially a primal active set method for quadratic programming, there is no reason to expect this number of iterations in general; and, indeed, known worst case bounds are exponential~\cite{nocedalNumericalOptimization2006}.
This is consistent with the picture that the NNLS solution lies at an extreme point of the feasible polytope, whose geometry is characterized by Carath\'{e}dory's theorem.
For the structured moment-matching problems considered here, the observed behavior appears to be essentially linear in $\dim\bbP_N^d$.

\begin{theorem}\label{thm:lawson-hanson-nnls-complexity}
  If \Cref{alg:lawson-hanson-nnls} terminates in $O(m)$ iterations so that $|\mathcal{J}| = O(m)$ is maintained throughout, its runs in $O(m^2 n)$ operations.
\end{theorem}

\begin{proof}
  Maintaining the various index sets takes $O(n)$ operations, amortized over the algorithm's iterations.
  The main contributions to the runtime are in computing the Lagrange multipliers, updating and downdating the QR decomposition, and performing an upper triangular solve to compute $\hat{\mw}_{\mathcal{J}}$.
  To compute $\malpha$ in Step 3, we first compute $\mV_{\mathcal{J}}^\top\mw_{\mathcal{J}}$ at a cost of $O(m^2)$ (since we assume $|\mathcal{J}| = O(m)$), followed by a matrix-vector product with $\mV$ for an additional $O(mn)$, making the overall cost $O(mn)$.
  In Steps 4 and 7 we update and downdate the QR factorization---the cost per step is $O(mn)$ due to updating $\mQ^\top\mV^\top$.
  Likewise, the cost the triangular solve is $O(m^2)$ per step.
  Summed over $O(m)$ steps, the overall cost is $O(m^2n)$.
\end{proof}

Lawson and Hanson's Fortran implementation is available in the public domain, and for many years it provided the heavy lifting for \texttt{scipy.optimize.nnls}.
It contains an important feature which provides a significant amount of additional stability and accuracy on ill-conditioned problems: the dual variables are read off the rotated residual carried by the running QR factorization, rather than recomputed from $\mV^\top\mw - \mI$.
In some cases, we found that the residual attained by the na\"{i}ve implementation can be off by many orders of magnitude relative to that attained by the stabilized one.
To the best of our knowledge, the reason for it is given neither in their book nor in a comment in any implementation we have examined.

The subsequent history of \texttt{scipy.optimize.nnls} illustrates how easily such a detail is lost.
SciPy 1.12 retired the Fortran in favor of a concise Python implementation~\cite{github-issue} organized around the normal equations, which recomputes the dual vector directly and so does not retain the safeguard.
A Cython rewrite in SciPy 1.15~\cite{polatRewriteNnlsCython2024}, and the C implementation current at the time of writing~\cite{polatRewriteSlsqpNnlsC2025}, both follow the original Fortran closely and thereby inherit it again.
We stress that none of this reflects an oversight.
The safeguard is invisible in the source, amounting to little more than a choice of which rows to accumulate, and a reimplementation that does not begin from the Fortran has no indication that anything has been given up; conversely, translating numerical code literally---as the later rewrites do---is exactly the right instinct, precisely because it preserves properties the translator has not had to identify.
The episode does suggest, though, that a safeguard which survives only as long as translation remains literal is one worth writing down.
For more details, see \Cref{sec:lawson-hanson}.

\begin{figure}
  \centering
  \includegraphics[width=\linewidth]{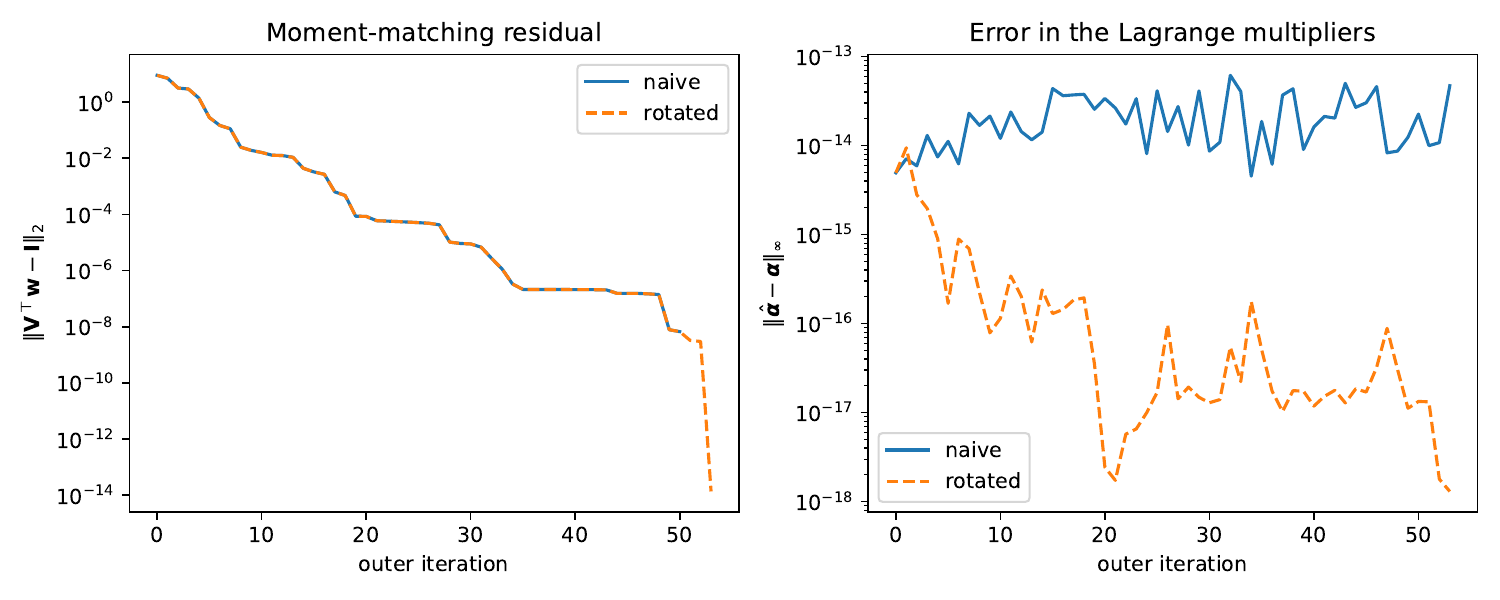}
  \caption{Two runs of \Cref{alg:lawson-hanson-nnls} on the same consistent moment-matching problem: $\Omega = [-1,1]$, the monomial basis up to degree $N = 20$, a uniform grid of $n = 60$ candidate nodes, and a moment vector generated from a nonnegative rule supported on twelve of those nodes. The two runs share an identical active set, least squares solve and termination criterion, and differ only in how the Lagrange multipliers are formed. \emph{Left}: the moment-matching residual against outer iteration. The two trajectories agree bit for bit through iteration $50$, after which the na\"{i}ve iteration stalls at $6.7 \times 10^{-9}$ while the stabilized iteration descends to $1.4 \times 10^{-14}$. \emph{Right}: the corresponding error in the Lagrange multipliers, measured against their values at the exact passive set least squares solution, computed in extended precision. The error in the direct expression \eqref{eq:lawson-hanson-lagrange-multipliers-1} never falls below roughly $10^{-14}$, while that in the stabilized expression \eqref{eq:lawson-hanson-lagrange-multipliers-2} decreases along with the residual; both are derived in \Cref{sec:lawson-hanson}.}\label{fig:nnls-dual-stall}
\end{figure}

\Cref{fig:nnls-dual-stall} isolates the effect on a problem small enough to inspect directly.
We take $\Omega = [-1, 1]$ with the monomial basis up to degree $N = 20$ and a uniform grid of $n = 60$ candidate nodes, and we choose the moment vector so that the problem is consistent: a nonnegative $\mw^\star$ supported on twelve of the grid points is drawn at random and $\mI = \mV^\top\mw^\star$ is formed from it.
An exact nonnegative solution therefore lies on the grid and the attainable residual is zero, which is what permits the iteration to descend far enough for the cancellation in the direct expression for the multipliers---\eqref{eq:lawson-hanson-lagrange-multipliers-1} of \Cref{sec:lawson-hanson}---to become visible at all.
We then run two implementations which share an identical active set, an identical least squares solve on the passive set, and an identical termination criterion, and which differ only in whether the multipliers are taken from that expression or from Lawson and Hanson's stabilized form, \eqref{eq:lawson-hanson-lagrange-multipliers-2}.
The two iterations agree bit for bit through iteration $50$.
At that point the na\"{i}ve multipliers can no longer resolve a descent direction and the iteration halts with nineteen nodes at a residual of $6.7 \times 10^{-9}$, while the stabilized iteration continues for three further steps, admits two more nodes, and terminates at $1.4 \times 10^{-14}$.
The right panel of \Cref{fig:nnls-dual-stall} shows the mechanism.
The error in \eqref{eq:lawson-hanson-lagrange-multipliers-1} stays near $10^{-14}$ for the whole run and never decreases, being governed by the cancellation incurred in forming $\mI - \mV^\top\mw$ rather than by the size of the residual; the error in \eqref{eq:lawson-hanson-lagrange-multipliers-2} falls with the residual itself, to roughly $10^{-17}$.
We stress that the size of the gap is strongly problem-dependent: over a sweep of degrees, grid sizes and random supports we observed everything from no appreciable difference to the gap shown here.

\subsection{Basis pursuit quadrature}\label{sec:basis-pursuit-quadrature}

\begin{figure}
  \centering
  \label{fig:nnls-vs-bp}
  \includegraphics[width=0.45\textwidth]{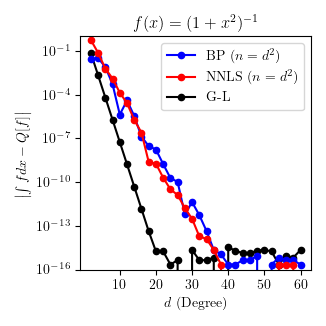}%
  \includegraphics[width=0.45\textwidth]{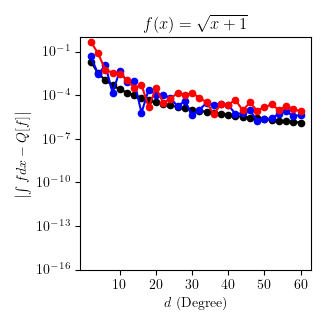}
  \includegraphics[width=0.45\textwidth]{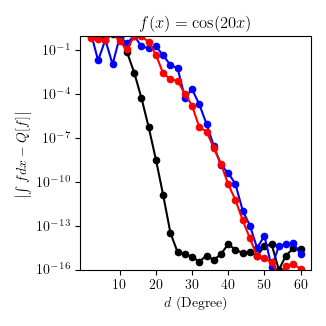}%
  \includegraphics[width=0.45\textwidth]{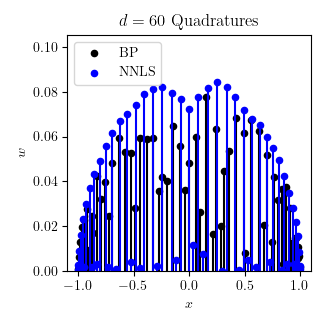}
  \caption{We compare the basis pursuit (BP) quadratures (\Cref{sec:basis-pursuit-quadrature}) and nonnegative least squares (NNLS) quadratures (\Cref{sec:nonnegative-least-squares-quadrature}) with Gauss--Legendre quadratures on the interval for several integrands up to degree 60. To avoid ill-conditioning, we formulate the problem in the Chebyshev basis. We observe that while the NNLS and BP quadratures are comparable in terms of their order of accuracy (being roughly one half of what is achieved by the equivalent Gauss--Legendre rule, as expected), we find that the NNLS rules are more symmetric, with a smooth decay of quadrature weights towards the boundary, as one might expect.}
\end{figure}

A classic model from the compressed sensing literature not too far removed from nonnegative least squares is the \emph{basis pursuit} problem~\cite{chenAtomicDecompositionBasis2001}:
\begin{equation}
  \begin{split}
    \label{eq:basis-pursuit}
    \mbox{minimize}\quad & \|\mw\|_1 \\
    \mbox{subject to}\quad & \mV^\top \mw = \mI \\
    & \mw \geq 0.
  \end{split}
\end{equation}
In the original formulation of basis pursuit, the nonnegativity constraints are optional. We include them in our formulation to present a form of basis pursuit pertinent to the design of positive quadratures.
This model was explored for the purposes of computing quadrature rules by Glaubitz~\cite{glaubitzStableHighorderCubature2021}, and still earlier by Mehrotra and Papp~\cite{mehrotraGeneratingMomentMatching2013}.

\begin{figure}
  \includegraphics[width=0.24\linewidth]{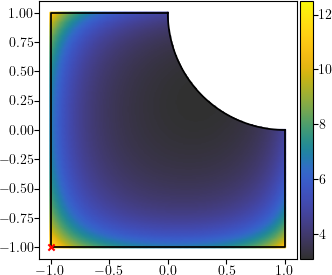}%
  \hspace{0.0066\linewidth}%
  \includegraphics[width=0.24\linewidth]{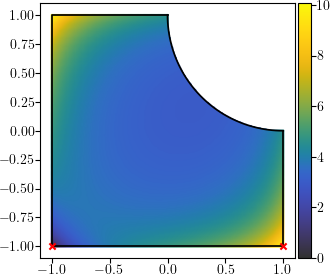}%
  \hspace{0.0066\linewidth}%
  \includegraphics[width=0.235\linewidth]{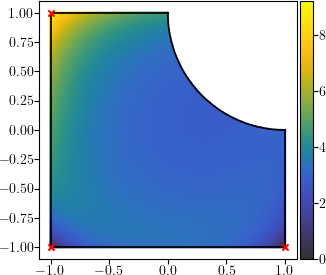}%
  \hspace{0.0066\linewidth}%
  \includegraphics[width=0.245\linewidth]{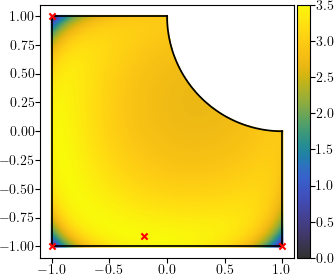}
  \includegraphics[width=0.24\linewidth]{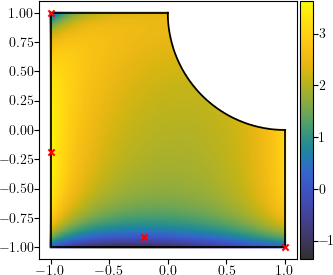}%
  \hspace{0.008\linewidth}%
  \includegraphics[width=0.241\linewidth]{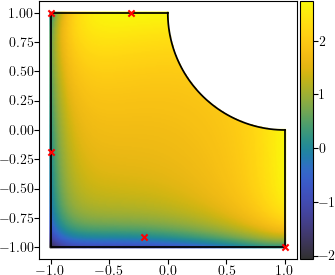}%
  \hspace{0.0066\linewidth}%
  \includegraphics[width=0.24\linewidth]{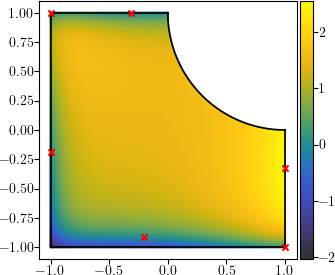}%
  \hspace{0.001\linewidth}%
  \includegraphics[width=0.249\linewidth]{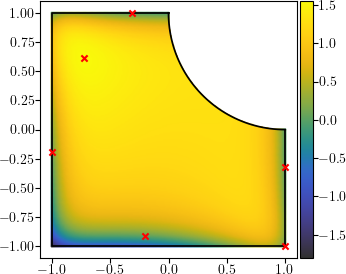}
  \caption{\emph{Left to right, top to bottom}: the Lagrange multiplier $\malpha$ for the first iterations of \Cref{alg:lawson-hanson-nnls}, the Lawson--Hanson iteration. Red crosses show the tentative locations of the quadrature nodes at each iteration. The grid resolution is deliberately chosen to be very fine to illustrate the smooth nature of the underlying Lagrange multiplier field.}\label{fig:nnls-LMs}
\end{figure}

In the presence of positivity and moment constraints, basis pursuit collapses to a feasibility problem.
Since $\mw \ge 0$, we have $\|\mw\|_1 = \m{1}^\top \mw$.
Because the constant function lies in $\mathbb{P}_N^d$, the vector $\mathbf{1}$ lies in the column space of $\mV$; hence, $\m{1}^\top \mw = \mathcal{Q}[1] = \operatorname{vol}(\Omega)$ for all feasible $\mw$.
The LP's cost is therefore constant over the entirety of the feasible set, reducing it to a feasibility testing problem.
In \Cref{sec:ryu-boyd-quadrature}, we will see that Ryu and Boyd's idea of introducing a nontrivial objective restores selectivity while preserving exactness.

Solving \eqref{eq:basis-pursuit} with an off-the-shelf implementation of the simplex method will produce an order $m$ quadrature, matching Tchakaloff (\Cref{thm:tchakaloff}).
On the other hand, if we solve it using an interior point method (see \Cref{sec:fast-primal-heuristic}), we may easily produce an interior point which is optimal but dense, resulting in an order $n$ quadrature.
Since all feasible $\mw$ have the same cost, the particular solution is governed by the algorithmic details of a black box solver.

A direct comparison between basis pursuit and NNLS on a uniform grid over $[-1,1]$ is shown in \Cref{fig:nnls-vs-bp}.
Both methods produce rules of comparable accuracy.
However, the rules produced by NNLS more naturally follow the expected limiting density of $\rho(x) = (1 - x^2)^{-1/2}$ for a quadrature supported on $[-1, 1]$.
An explanation for this phenomenon is provided by \Cref{fig:nnls-LMs}, where we plot the Lagrange multipliers for the first several iterations of the Lawson--Hanson algorithm.
The Lagrange multipliers for \eqref{eq:nnls} are given by $\malpha = \mV(\mI - \mV^\top\mw)$, so they can be interpreted as samples of a polynomial in $\mathbb{P}_N^d$.
Greedily reducing the maximum value of $\malpha$ naturally guides the Lawson--Hanson iteration toward well-placed nodes.
On the other hand, there is no reason to expect that the node distribution provided by the simplex algorithm will be anything other than biased by different pivoting strategies.
For this reason, we regard Lawson--Hanson (\Cref{alg:lawson-hanson-nnls}) as the preferred method.

\subsection{Gauss--Newton solution of the moment-matching equations}\label{sec:gauss-newton}

The remaining methods in the quadrature toolbox use the uniform grid to produce an initial quadrature whose moment-matching residual is already small.
Each method then refines that initial rule by applying Newton’s method to the moment-matching equations.
We allow both $\mw$ and $\mX$ to vary and seek an $m$-dimensional zero of the $(d+1)n$-dimensional nonlinear system $\mF(\mw, \mX) = \m{0}$ (see \Cref{sec:moment-matched-quadratures} for definitions).
Our assumptions is that $\eff(\mathcal{Q}) = \tfrac{m}{(d+1)n} \leq 1$ for any feasible quadrature, it only makes sense to consider cases with $m \leq (d+1)n$.
Thus, the moment-matching system is either square or underdetermined.
In the underdetermined case, the solutions of $\mF(\mw, \mX)=0$ form a smooth nonlinear manifold of dimension $(d+1)n - m$, unless some degeneracy occurs.
Enhancements to the basic Newton's method exist for underdetermined systems~\cite{polyakSolvingUnderdeterminedNonlinear}, but we do not pursue them here.

As noted in \Cref{sec:moment-matched-quadratures}, we may apply Newton’s method directly to the rootfinding problem (using a generalized pseudoinverse step), or Gauss--Newton to the nonlinear residual.
In either case the resulting iteration is identical; for completeness, we derive the Newton form.
Let $\mw_*$ and $\mX_*$ satisfy $\mF(\mw_*, \mX_*) = \m{0}$, and assume that $\mw_n$ and $\mX_n$ denote the $n$th approximation of $\mw_*$ and $\mX_*$.
Then, if $\mw_* = \mw_n + \delta\mw_n$ and $\mX_* = \mX_n + \delta\mX_n$, we Taylor expand:
\begin{equation}
  \begin{split}
    \mF(\mw_*, \mX_*) &= \mF(\mw_n + \delta\mw_n, \mX_n + \delta\mX_n) \\
                      &= \mF(\mw_n, \mX_n) + \mathsf{D}_{\mw}\mF(\mw_n, \mX_n) \delta \mw_n \\
    &\hspace{4em} + \mathsf{D}_{\mX}\mF(\mw_n, \mX_n) \delta \mX_n + O(\|\delta\mw_n\|^2) + O(\|\delta\mX_n\|^2),
  \end{split}
\end{equation}
suggesting the iteration:
\begin{equation}
  \begin{bmatrix}
    \mw_{n+1} \\
    \mX_{n+1}
  \end{bmatrix} = \begin{bmatrix}
    \mw_{n} \\
    \mX_{n}
  \end{bmatrix} - \begin{bmatrix}
    \mathsf{D}_{\mw}\mF(\mw_n, \mX_n) & \mathsf{D}_{\mX}\mF(\mw_n, \mX_n)
  \end{bmatrix}^\dagger \mF(\mw_n, \mX_n).
\end{equation}
Kantorovich's theorem gives conditions for quadratic convergence~\cite{stoerIntroductionNumericalAnalysis2002a}.
When Kantorovich's theorem holds, convergence is cheap and rapid.

\begin{theorem}\label{thm:gauss-newton-complexity}
  If $m \sim (d+1)n$ and $(\mw_0, \mX_0)$ is already in the quadratic convergence region, Newton's method applied to the moment-matching equations attains $O(\epsilon)$ accuracy in $O(m^3 \log \log \tfrac{1}{\epsilon})$ operations.
\end{theorem}
  
\begin{proof}
  We assume $m \sim (d+1)n$ since we will only apply Newton's method when $m$ is within a modest constant factor of $(d+1)n$ (e.g., it will always be that case that $\eff(\mathcal{Q}) \geq (d+1)^{-1}$).
  Then, each step of Newton's method takes $O(m^3)$ to solve a least squares problem to apply the pseudoinverse.
  Let $\epsilon_n$ be the error in the $\ell_2$ norm at the $n$th step.
  Since quadratic convergence is immediately attained, there is some constant $C > 0$ such that $\epsilon_{n+1} \leq C \epsilon_n^2$ for $n \geq 0$.
  This recursion can be solved to give $\epsilon_n \leq C^{2^n - 1} \epsilon_0^{2^n}$.
  Taking logs shows that $n = O(\log \log \tfrac{1}{\epsilon})$, from which the result can be concluded.
\end{proof}

In general, Newton’s method need not converge for an arbitrary initial guess $(\mw_0,\mX_0)$.
However, the methods we introduce below reliably produce initial quadratures with small moment-matching residuals, and these almost always lie inside the basin of quadratic convergence.
Stabilized variants such as Levenberg--Marquardt guarantee convergence to a local minimizer~\cite{nocedalNumericalOptimization2006}, but we view these as unnecessary for quadrature optimization: if Newton diverges, decreasing the grid size $h$ or subdividing the cut cell is typically a far more effective strategy.
Forcing convergence with a stabilized method often yields a lower-quality rule at significantly greater computational cost.

\subsection{Ryu--Boyd quadrature}\label{sec:ryu-boyd-quadrature}

The only difference between the basis pursuit quadratures (\cref{sec:basis-pursuit-quadrature}) and the ``Gauss LP'' quadratures due to Ryu and Boyd is the choice of cost function~\cite{ryuExtensionsGaussQuadrature2015}.
Let $r(\mx) : \Omega \to \bbR$, and let $\mr = r(\mX)$.
Ryu and Boyd refer to this as a \emph{sensitivity function}.
The LP for Ryu--Boyd quadrature is:
\begin{equation}
  \begin{split}
    \label{eq:ryu-boyd-LP}
    \mbox{minimize}\quad & \mr^\top \mw, \\
    \mbox{subject to}\quad & \mV^\top \mw = \mI, \\
    & \mw \geq 0.
  \end{split}
\end{equation}
From this, we can see that basis pursuit is an instance of Ryu--Boyd with $r(\mx) = 1$.
As we saw in \Cref{sec:basis-pursuit-quadrature}, basis pursuit is totally insensitive to the choice of quadrature nodes, since we always consider polynomial spaces that contain constant functions.

\begin{figure}
  \centering
  \includegraphics[width=0.48\linewidth]{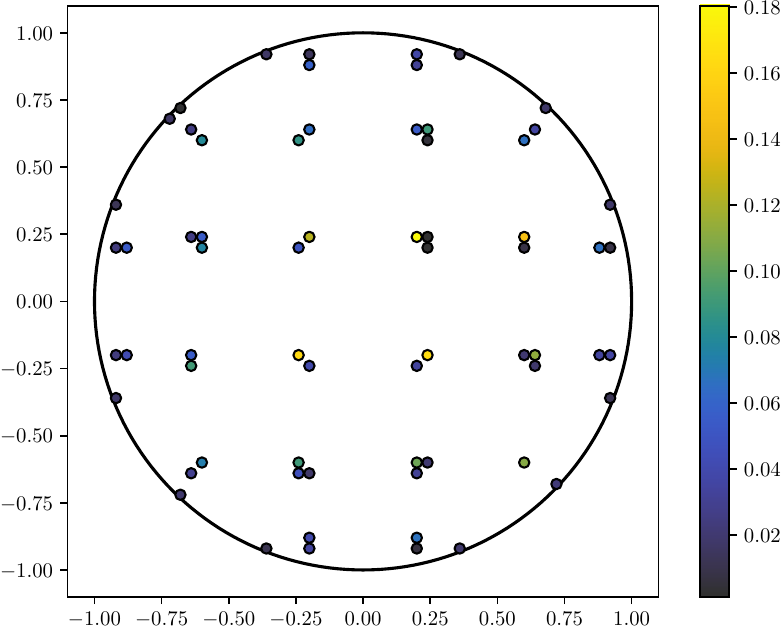}%
  \hspace{0.025\linewidth}%
  \includegraphics[width=0.48\linewidth]{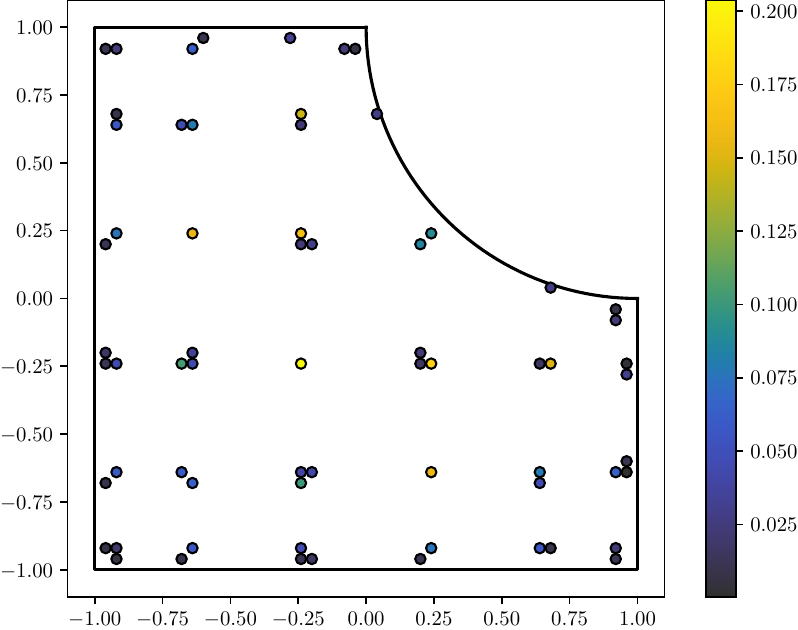}%
  \vspace{0.025\linewidth}
  \includegraphics[width=0.48\linewidth]{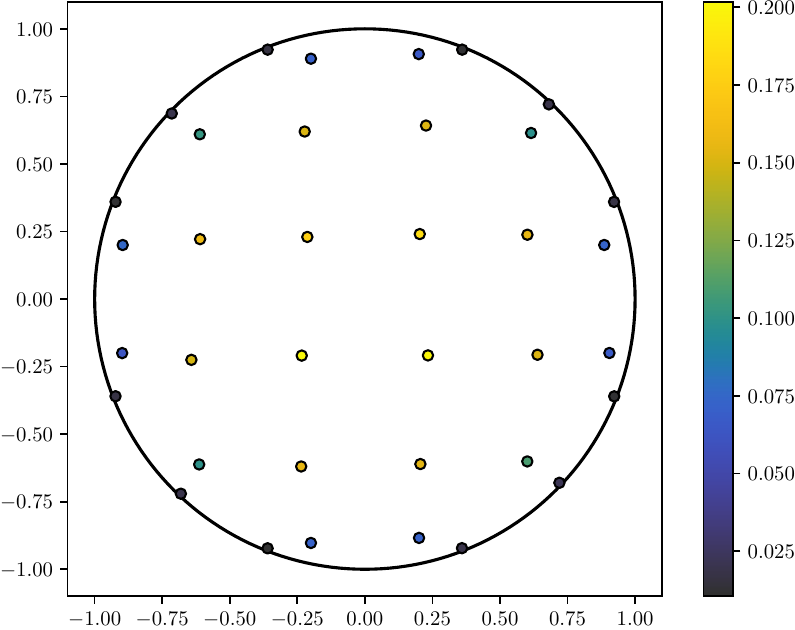}%
  \hspace{0.025\linewidth}%
  \includegraphics[width=0.48\linewidth]{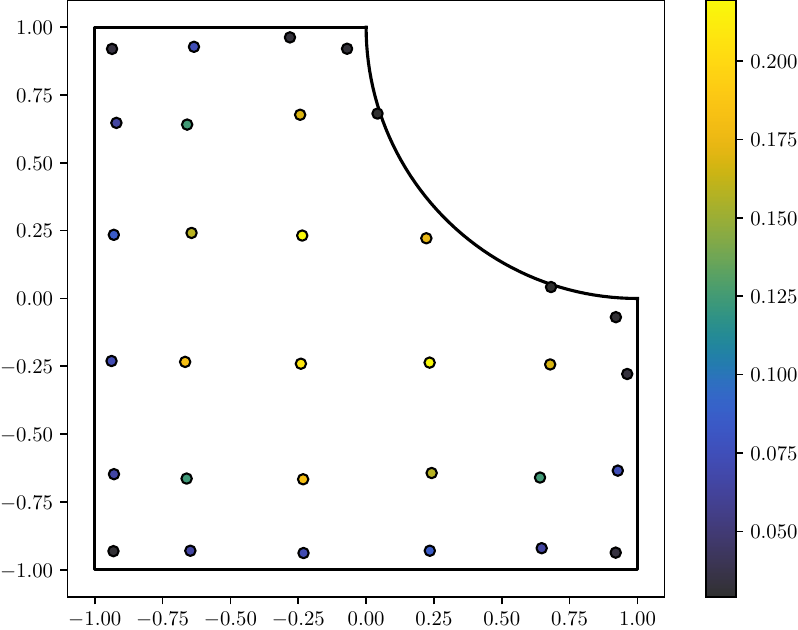}
  \caption{The original Ryu--Boyd algorithm (as described in \Cref{sec:ryu-boyd-quadrature}) applied to the two test problems shown in \Cref{fig:simple-nnls-2d} (the left and right columns in this figure correspond to the left and right columns in \Cref{fig:simple-nnls-2d}). The same spatial grid and moment-matching system are used. The residual function $r(x, y) = x^{20} + y^{20}$ is used. \emph{Top}: the LP \eqref{eq:ryu-boyd-LP} is solved using the simplex method. The resulting quadrature is shown, exhibiting the effect originally observed by Ryu and Boyd, as small clusters form near ideal node placements. These quadratures integrate $\mathbb{P}_{10}^2$ exactly and have exactly $\dim\mathbb{P}_{10}^2 = 66$ nodes. \emph{Bottom}: The nodes are clustered following Ryu and Boyd's prescription and subsequently optimized using Gauss--Newton. After running Gauss--Newton, all nodes remain in the interior of $\Omega$ and all weights remain positive. The left and right column quadratures integrate $\mathbb{P}_{10}^2$ within numerical roundoff and have 36 and 32 nodes, respectively. This is a significant increase in efficiency compared to the NNLS quadratures. Note that the minimum number of nodes possible for these test problems is 22; so, the quadratures have an efficiency of $0.6\overline{1}$ and 0.6875, respectively, while the NNLS quadratures each have an efficiency of $1/3$.}\label{fig:ryu-boyd-2d}
\end{figure}

What makes a good sensitivity function is open-ended---several choices are discussed in the original paper.
Clearly, we should choose a function $r$ that doesn't lie in the polynomial space being integrated so that $\mr$ isn't contained in the column space of $\mV$.
Ryu and Boyd prove that if the infinite-dimensional version of \eqref{eq:ryu-boyd-LP} is solved on $[-1, 1] \subseteq \bbR$ with $r(x) = x^{2m}$, then its solution recovers the corresponding Gaussian quadrature.
They did not prove but observed an analogous result in higher dimensions, but they observe experimentally that a choice like $r(x, y) = x^{2m} + y^{2m}$ for the polynomial space $\mathbb{P}_m^2$ leads to good results.
They also observe that while choosing $r(x, y) = x^{2m} + y^{2m}$ tends to produce ``Gauss'' quadratures (i.e., efficient quadratures with strictly interior points), the choice $r(x, y) = -x^{2m} - y^{2m}$ yields ``Lobatto'' quadratures, where many points from the boundary are selected.

These observations can be contextualized somewhat by considering the convex geometry of the feasible set.
Namely, for a fine enough uniform grid sampling $\Omega$, we observe that:
\begin{itemize}[nosep]
\item ordinary least squares produces a fully interior point (\Cref{sec:least-squares-quadrature}),
\item nonnegative least squares using the Lawson--Hanson algorithm produces an order $m$ quadrature corresponding to a point on the boundary of the feasible set, lying on a vertex or some facet of the polytope (\Cref{sec:nonnegative-least-squares-quadrature}),
\item basis pursuit does likewise, although the difference in algorithmic details biases the result to a different order $m$ quadrature (\Cref{sec:basis-pursuit-quadrature}).
\end{itemize}
Consequently, it can be seen that introducing the sensitivity function ``steers'' the optimum of the LP to different quadratures which are at most order $m$.

In fact, this steering effect frequently produces quadratures with nodes organized into small clusters centered around an off-grid point.
Ryu and Boyd presented a heuristic method in which these clusters are replaced with their weighted average and used as a warm start for Gauss--Newton (\Cref{sec:gauss-newton}).
See \Cref{fig:ryu-boyd-2d} for an overview.
A sketch of the Ryu--Boyd algorithm as presented in the original paper follows:
\begin{sfpalgorithm}[Ryu and Boyd's quadrature optimization algorithm.]
  \begin{enumerate}[nosep]
  \item Following \Cref{ssec:discretizing-the-cut-cell}, let $\mX$ be a uniform grid sampling $\Omega$ (let $h > 0$ be the grid spacing).
  \item Solve \eqref{eq:ryu-boyd-LP}, letting $\mw$ denote the optimum.
  \item Let $\mX_{\operatorname{LP}} = \{\mx \in \mX : \mw_{\mx} > 0\}$.
  \item Let $\mw_{\operatorname{LP}} = \{\mw_{\mx} : \mx \in \mX_{\operatorname{LP}}\}$.
  \item Cluster $\mX_{\operatorname{LP}}$ to obtain $\mX_{\operatorname{clust}}$ and $\mw_{\operatorname{clust}}$. This is heuristic. As an example, we could call two points $\mx, \my \in \mX$ ``neighbors'' if $\|\mx - \my\|_{\infty} \leq h$ (giving $3^d - 1$ point neighborhoods). Finding all connected components in the resulting graph gives a connected component for each cluster: to find $\mx \in \mX_{\operatorname{clust}}$, compute their mean. To obtain $\mw_{\operatorname{clust}}$, sum together the weights of each node in the connected component (this maintains $\sum_{\mx \in \mX_{\operatorname{clust}}} \mw_{\mx} = \operatorname{vol}(\Omega)$).
  \item Use Gauss--Newton to optimize the moment-matching equations with $(\mw_{\operatorname{clust}}, \mX_{\operatorname{clust}})$ as a warm start---denote the result by $(\mw_{\operatorname{opt}}, \mX_{\operatorname{opt}})$.
  \item The quadrature determined by $(\mw_{\operatorname{opt}}, \mX_{\operatorname{opt}})$ is the result.
  \end{enumerate}
  \label{algo:ryu-boyd-quadrature}
\end{sfpalgorithm}
\noindent If $h$ is too large, this algorithm can fail in two ways.
First, clusters might begin to overlap or connect, preventing us from determining the location of all quadrature nodes.
Second, the Gauss--Newton iteration may diverge if $(\mw_{\operatorname{clust}}, \mX_{\operatorname{clust}})$ is an insufficiently good warm start. 
Both of these problems are mitigated by letting $h \to 0$, giving us some confidence that this algorithm is consistent.

\subsection{Xiao--Gimbutas elimination}\label{sec:xiao-gimbutas-elimination}

\begin{figure}
  \includegraphics[width=0.48\linewidth]{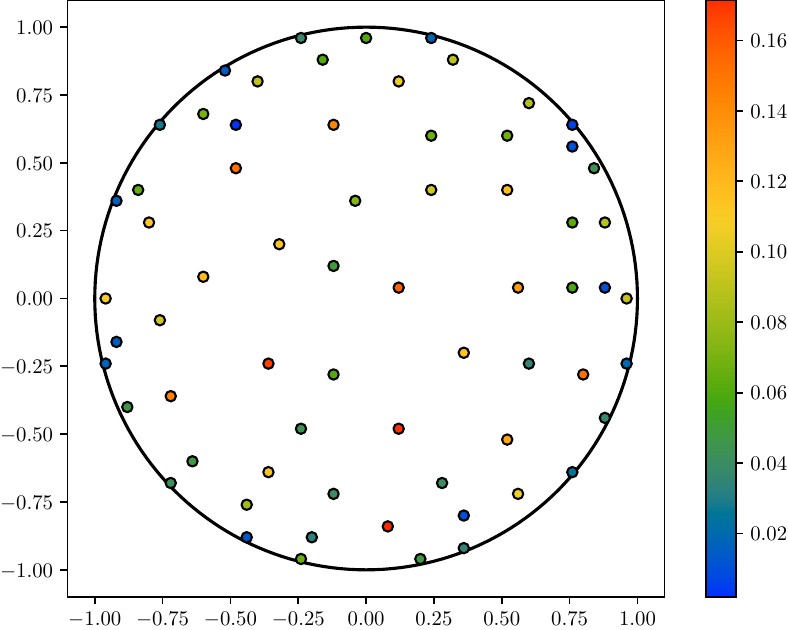}%
  \hspace{0.025\linewidth}%
  \includegraphics[width=0.48\linewidth]{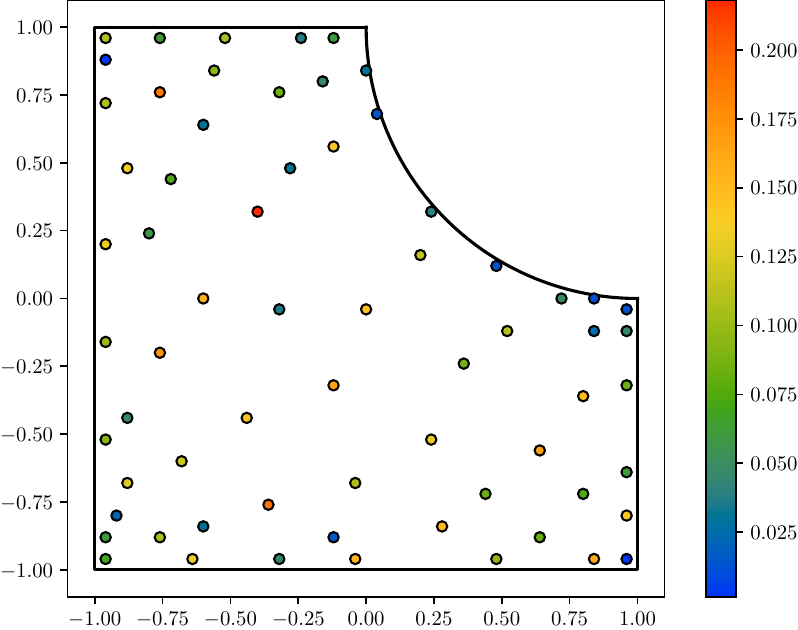}
  \caption{Plots of node significances $w_j \|\mphi(\mx_j)\|_2$ (defined in \Cref{sec:xiao-gimbutas-elimination}) for the NNLS quadratures computed in \Cref{fig:simple-nnls-2d}.}
\end{figure}

\begin{figure}
  \includegraphics[width=0.48\linewidth]{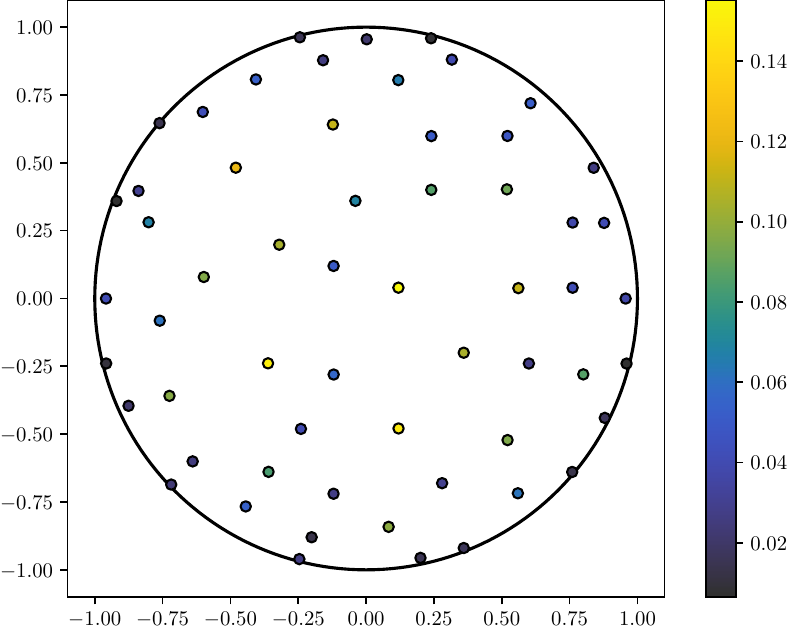}%
  \hspace{0.025\linewidth}%
  \includegraphics[width=0.48\linewidth]{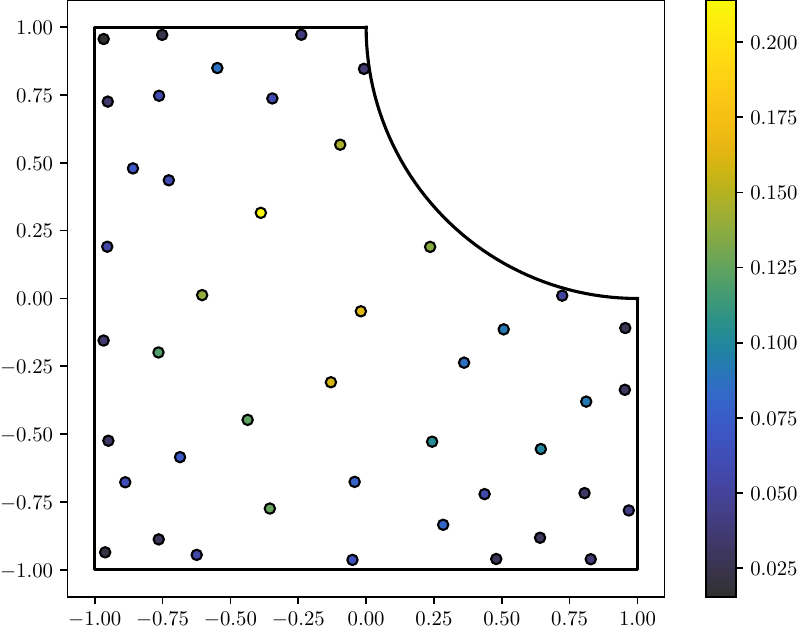}
  \caption{NNLS quadratures from \Cref{fig:simple-nnls-2d} after using Xiao--Gimbutas elimination to remove one node at a time until Gauss--Newton fails to converge. In each case, we can see that this process is able to remove a significant number of nodes.}\label{fig:xiao-gimbutas-nnls-2d}
\end{figure}

Steinitz elimination (\Cref{sec:steinitz-elimination}) allows us to reduce an order $n > m$ positive quadrature to an order $n = m$ positive quadrature.
If $\calQ$ is the Steinitz-eliminated quadrature, then $\eff(\calQ) = (d+1)^{-1}$.
As we saw in the previous section, examples of positive quadratures in $d$ dimensions with $\eff(\mathcal{Q}) > (d + 1)^{-1}$ abound.
It is natural to ask whether we can further economize the quadrature by eliminating nodes in some other way.
Sadly, little is known in general about optimal multivariate quadratures, so it is unclear what our notion of optimality should be, other than driving $\eff(\calQ) \leq 1$ as high as possible.
Nevertheless, Xiao and Gimbutas provide just such a tool~\cite{xiaoNumericalAlgorithmConstruction2010}.
We summarize their approach to node elimination here.

Xiao and Gimbutas define the \emph{significance} of a quadrature node to be either:
\begin{equation}
  \sigma_j^{\operatorname{XG1}} = \sum_{i=1}^m \phi_i(\mx_j)^2 \qquad\mbox{or}\qquad \sigma_j^{\operatorname{XG2}} = w_j \sum_{i=1}^m \phi_i(\mx_j)^2. \label{eq:xiao-gimbutas-significance}
\end{equation}
We can see that $\sigma_j$ will be large if it contributes a significant amount to the moments  and small otherwise.
Nodes are eliminated one at a time by repeatedly eliminating the node with the smallest significance and re-optimizing using Gauss--Newton.
Ultimately, their goal is to design quadratures for canonical domains which are as close to optimal ($\eff(\calQ) = 1$) as possible.
Consequently, they recursively backtrack over possible node elimination sequences, keeping the most efficient quadrature found.
This recursion terminates if Gauss--Newton diverges or if it converges to a quadrature with negative weights or nodes outside the integration region.
For canonical domains, their quadratures are likely the best that exist in many cases~\cite{kloecknerModepy2024}.

Xiao and Gimbutas' measure of significance is simply related to the moment-matching equations.
Let $\calQ$ be a quadrature with weights $w_1, \hdots, w_n$ and nodes $\mx_1, \hdots, \mx_n$ which integrates a polynomial space $\bbP = \Span\,\{\phi_1, \hdots, \phi_m\}$, and let $\tilde{\calQ}$ denote the quadrature obtained by removing the $j$th node from $\calQ$.
Let $\mV$ and $\mw$ be the basis matrix and quadrature weight vector associated with $\calQ$; likewise, define $\tilde{\mV}$ and $\tilde{\mw}$ to go along with $\tilde\calQ$.
As usual, let $\mI$ be the moment vector.
Since $\calQ$ integrates $\bbP$, the moment-matching equations are satisfied.
Define $\mphi(\mx) = (\phi_1(\mx), \hdots, \phi_m(\mx))$.
With $\mphi(\mx)$ so defined, the Xiao--Gimbutas significances become $\sigma_j^{\operatorname{XG1}} = \|\mphi(\mx_j)\|_2^2$ and $\sigma_j^{\operatorname{XG2}} = w_j \|\mphi(\mx_j)\|_2^2$, respectively.
It is easy to see that $\mI - \tilde{\mV}{}^\top\tilde{\mw} = w_j\mphi(\mx_j)$.
That is, $w_j\mphi(\mx_j)$ is exactly the residual of the moment-matching equations for $\tilde{\calQ}$.
Clearly, the closer either of the two Xiao--Gimbutas significances are to zero, the closer the moment-matching equations for $\tilde{\calQ}$ are to being satisfied.

The full Xiao--Gimbutas algorithm (including recursive backtracking) is too costly to run on the fly.
However, to take advantage of their method of eliminating nodes in the context of quickly generating quadrature rules for many different unstructured cut cells, we can try speculatively eliminating nodes.
Clearly, the sequence of eliminated nodes depends on the definition of $\sigma_j$.
\begin{sfpalgorithm}[Simple greedy algorithm for node elimination.]
  \label{alg:greedy-node-elimination}
  \begin{enumerate}[nosep]
  \item Compute $\m{\sigma} = (\sigma_1, \hdots, \sigma_N)$ for some definition of $\sigma_j$.
  \item Remove the $j$th node, where $j = \Arg\min_j \sigma_j$.
  \item Use Gauss--Newton to reoptimize the quadrature by solving the moment-matching equations.
  \item Repeat from (1) while all weights are positive and all nodes lie within $\Omega$.
  \end{enumerate}
\end{sfpalgorithm}

We collect some numerical experiments showing the result of pruning different NNLS quadratures in \Cref{table:pruning}.
We note that the Xiao--Gimbutas significances are not \emph{exactly} the desired residual.
To this end, we include another node significance in our experiment which is exactly the norm of the moment-matching residual contribution of the $j$th node:
\begin{equation}
  \sigma_j^{\operatorname{Res}} = w_j \|\mphi(\mx_j)\|_2.
\end{equation}
We found that $\sigma_j^{\operatorname{Res}}$ performs similarly to XG2.
However, In practice, we found that $\sigma_j^{\operatorname{XG1}}$ fares quite poorly, typically failing to eliminate very many nodes at all.

\begin{table}
  \centering
  \setlength\tabcolsep{4.5pt}
  {\footnotesize
  % \rotatebox{90}{
    \begin{tabular}{ccccccc|ccccccc}
  \multicolumn{7}{c}{Unit disk} & \multicolumn{7}{c}{Square with circular bite} \\
  \hline
  $N$ & $\phi_{\max}$ & $\sigma_j$ & $\dim\mathbb{P}_N^2$ & $n$ & \mbox{Eff} & GN &$N$ & $\phi_{\max}$ & $\sigma_j$ & $\dim\mathbb{P}_N^2$ & $n$ & \mbox{Eff} & GN \\
  \hline
  8 & -0.003 & Res & 45 & 37 & \emph{40.5\%} & 5 & 8 & -0.003 & Res & 45 & 43 & 34.9\% & 5 \\
  8 & -0.003 & XG1 & 45 & 44 & 34.1\% & 8 & 5 & -0.003 & XG1 & 45 & 42 & \textbf{35.7\%} & 7 \\
  8 & -0.003 & XG2 & 45 & 37 & \emph{40.5\%} & 5 & 8 & -0.003 & XG2 & 45 & 43 & 34.9\% & 5 \\
  \hline
  8 & -0.01 & Res & 45 & 35 & 42.9\% & 5 & 8 & -0.01 & Res & 45 & 35 & 42.9\% & 6 \\
  8 & -0.01 & XG1 & 45 & 44 & 34.1\% & 5 & 8 & -0.01 & XG1 & 45 & 42 & 35.7\% & 7 \\
  8 & -0.01 & XG2 & 45 & 34 & \textbf{44.1\%} & 5 & 8 & -0.01 & XG2 & 45 & 34 & \textbf{44.1\%} & 6 \\
  \hline
  8 & -0.03 & Res & 45 & 32 & 46.9\% & 6 & 8 & -0.03 & Res & 45 & 29 & \emph{51.7\%} & 10 \\
  8 & -0.03 & XG1 & 45 & 44 & 34.1\% & 5 & 8 & -0.03 & XG1 & 45 & 42 & 35.7\% & 7 \\
  8 & -0.03 & XG2 & 45 & 30 & \textbf{50.0\%} & 5 & 8 & -0.03 & XG2 & 45 & 29 & \emph{51.7\%} & 10 \\
  \hline
  10 & -0.003 & Res & 66 & 56 & \emph{39.3\%} & 5 & 10 & -0.003 & Res & 66 & 58 & \textbf{37.9\%} & 5 \\
  10 & -0.003 & XG1 & 66 & 63 & 34.9\% & 8 & 10 & -0.003 & XG1 & 66 & 61 & 36.1\% & 8 \\
  10 & -0.003 & XG2 & 66 & 56 & \emph{39.3\%} & 5 & 10 & -0.003 & XG2 & 66 & 60 & 36.7\% & 5 \\
  \hline
  10 & -0.01 & Res & 66 & 54 & 40.7\% & 5 & 10 & -0.01 & Res & 66 & 56 & 39.3\% & 6 \\
  10 & -0.01 & XG1 & 66 & 63 & 34.9\% & 7 & 10 & -0.01 & XG1 & 66 & 64 & 34.4\% & 6 \\
  10 & -0.01 & XG2 & 66 & 52 & \textbf{42.3\%} & 5 & 10 & -0.01 & XG2 & 66 & 54 & \textbf{40.7\%} & 6 \\
  \hline
  10 & -0.03 & Res & 66 & 48 & \emph{45.8\%} & 7 & 10 & -0.03 & Res & 66 & 47 & 46.8\% & 8 \\
  10 & -0.03 & XG1 & 66 & 63 & 34.9\% & 6 & 10 & -0.03 & XG1 & 66 & 62 & 35.5\% & 8 \\
  10 & -0.03 & XG2 & 66 & 48 & \emph{45.8\%} & 7 & 10 & -0.03 & XG2 & 66 & 45 & \textbf{48.9\%} & 9 \\
  \hline
  12 & -0.003 & Res & 91 & 85 & 36.5\% & 4 & 12 & -0.003 & Res & 91 & 79 & \textbf{39.2\%} & 7 \\
  12 & -0.003 & XG1 & 91 & 90 & 34.4\% & 5 & 12 & -0.003 & XG1 & 91 & 90 & 34.4\% & 6 \\
  12 & -0.003 & XG2 & 91 & 78 & 39.7\% & 6 & 12 & -0.003 & XG2 & 91 & 84 & 36.9\% & 6 \\
  \hline
  12 & -0.01 & Res & 91 & 80 & 38.8\% & 7 & 12 & -0.01 & Res & 91 & 76 & \emph{40.8\%} & 7 \\
  12 & -0.01 & XG1 & 91 & 90 & 34.4\% & 5 & 12 & -0.01 & XG1 & 91 & 88 & 35.2\% & 7 \\
  12 & -0.01 & XG2 & 91 & 78 & \textbf{39.7\%} & 5 & 12 & -0.01 & XG2 & 91 & 76 & \emph{40.8\%} & 5 \\
  \hline
  12 & -0.03 & Res & 91 & 63 & 49.2\% & 7 & 12 & -0.03 & Res & 91 & 65 & 47.7\% & 11 \\
  12 & -0.03 & XG1 & 91 & 89 & 34.8\% & 6 & 12 & -0.03 & XG1 & 91 & 88 & 35.2\% & 9 \\
  12 & -0.03 & XG2 & 91 & 62 & \textbf{50.0\%} & 7 & 12 & -0.03 & XG2 & 91 & 59 & \textbf{52.5\%} & 13
\end{tabular}

    % }
    }
  \caption{For our two test domains, and for a variety of degrees, we run Xiao--Gimbutas elimination until failure for each node significance. The column labels, from left to right: the total degree ($N$), the maximum level set value used by inside-outside test ($\phi_{\max}$), the node significance used ($\sigma_j$), the number of basis functions ($\dim\bbP_N^2$), the number of quadrature nodes ($n$), the quadrature efficiency ($\eff$), and the maximum number of Gauss--Newton iterations observed during Xiao--Gimbutas elimination (GN). We find that the $\ell_2$-based significance and XG2 perform the best and are routinely capable of boosting the efficiency of each quadrature by 10-15\%.}\label{table:pruning}
\end{table}

Reoptimizing the quadrature rule using Gauss--Newton after eliminating a node is quite cheap compared to solving the original NNLS problem.
At the beginning of running the greedy algorithm, the eliminated quadratures will tend to be very good warm starts for the Gauss--Newton iteration.
The conditions of \Cref{thm:gauss-newton-complexity} typically hold, with Gauss--Newton converging to numerical roundoff in 5 or 6 iterations (see \Cref{table:pruning}).
The initial NNLS quadrature will have $m$ nodes, and the minimum possible is $\lceil m/(d+1) \rceil$.
Conservatively, if we assume the cost of each elimination is $O(m^3)$, then the overall cost of running the simple greedy algorithm for elimination is $O(m^4)$.

\paragraph{Adapting the NNLS discretization for Xiao--Gimbutas elimination}
If the bounding box of $\Omega$ is tight and we have many grid points on $\Omega$'s boundary, the Lawson--Hanson minimizer will tend to include points along the boundary (the Lawson--Hanson quadratures tend to be Lobatto quadratures).
Deleting a node from one of these Lobatto quadratures and reoptimizing will spread the remaining points out slightly, pushing some boundary points outside of $\Omega$, violating one of our constraints.
To avoid this problem, when sampling our grid, we can choose a value $\phi_{\max}$ slightly less than zero and only keep grid points $\mx$ which satisfy $\phi(\mx) \leq \phi_{\max}$ to avoid an initial grid with boundary points.
With this modification, \Cref{alg:greedy-node-elimination} will often run for a significant number of steps before any nodes exit the domain, allowing for a significantly improvement in $\eff(\mathcal{Q})$.

\subsection{Pruning inefficient rules}\label{sec:pruning-inefficient-rules}

\begin{figure}
  \centering
  \begin{subfigure}[b]{0.345\textwidth}
    \includegraphics[width=\textwidth]{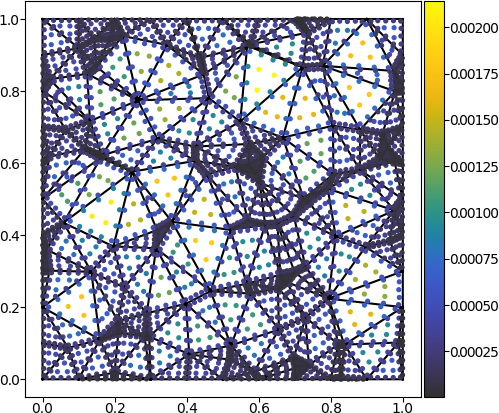}
    \caption{Unpruned}
    \label{fig:first}
  \end{subfigure}%
  \begin{subfigure}[b]{0.3253\textwidth}
    \includegraphics[width=\textwidth]{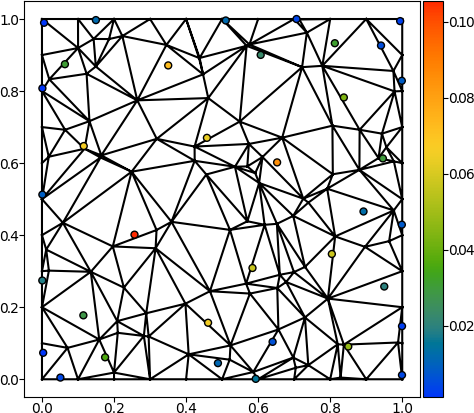}
    \caption{NNLS}
    \label{fig:second}
  \end{subfigure}%
  \begin{subfigure}[b]{0.3253\textwidth}
    \includegraphics[width=\textwidth]{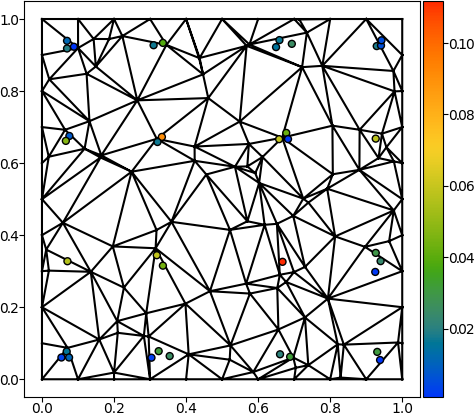}
    \caption{LP}
    \label{fig:third}
  \end{subfigure}
\caption{Pruning a composite quadrature on the square $\Omega=[0,1]^2$.
\emph{Left}: A composite rule obtained by meshing $\Omega$ and mapping Xiao--Gimbutas simplex rules~\cite{xiaoNumericalAlgorithmConstruction2010,kloecknerModepy2024} onto each triangle (see also \Cref{sec:xiao-gimbutas-elimination}).
Each elementwise rule integrates $\mathbb{P}_7^2$ exactly, yielding a global quadrature that is exact but highly inefficient.
\emph{Middle}: Pruning the composite rule using NNLS (\eqref{eq:nnls}) via the Lawson--Hanson algorithm (see \Cref{sec:lawson-hanson}) collapses the rule to a dramatically smaller set of nodes while preserving exactness.
\emph{Right}: A similar reduction obtained by solving the Ryu--Boyd LP~\eqref{eq:ryu-boyd-LP} with sensitivity function $r(x,y)=x^{14}+y^{14}$.
Both approaches illustrate how sparse optimization can be used to extract an efficient quadrature rule native to $\Omega$ from an oversampled composite rule.}\label{fig:pruning}
\end{figure}

Although our primary message is that an inside-outside test along with a moment vector is sufficient to compute a good quadrature rule, in many existing codes a convenient first step is to rely on meshing infrastructure to obtain candidate nodes and moments.
In such settings, composite elementwise quadratures serve purely as a moment generator, not as the final rule.
Our sparse optimization methods then extract a global quadrature with far fewer nodes.
This section illustrates how to reinterpret traditional element-based quadratures through the moment-matching lens.

The simple observation is that if we have a quadrature rule with more than $m$ nodes integrating an $m$-dimensional space, we can simply invoke \Cref{thm:caratheodory-quadrature} and run any of the sparse optimization methods discussed in this section to eliminate unnecessary nodes.
As an illustration, consider the simple setup shown in \cref{fig:pruning}, in which a square domain is meshed and simple Gauss rules on triangles are used to compute the desired moments.
The extension of this to polygonal or polyhedral domains produces an exact composite rule; for curved domains, adaptive refinement can ensure the mapped rule computes the moments to within a prescribed error tolerance.
Either way, the resulting composite quadrature gives us exactly what we need: a set of candidate points and the correct moment vector for the target polynomial space.
Once these moments are known, the elementwise structure becomes irrelevant.
It is important to note that while the original composite rule also integrates a larger space of piecewise polynomials on the mesh, this additional exactness is often unnecessary: the pruned rule integrates the intended global space on $\Omega$ and is far more efficient.
In practice, this provides a simple and robust pathway from legacy mesh-based workflows to compact, high-quality moment-matched quadratures.

% \begin{table}
%   \centering
%   \input{tables/table1.tex}
%   \caption{For the same setup as \Cref{fig:gauss-lp-3d}, some statistics for the primal-dual iterations. These statistics particularly concern the selected warm start (the ``clustered'' quadrature) and the subsequent Gauss-Newton iteration for solving the moment-matching equations. The statistics show the ``quality'' of the warm start is only marginally improved with additional primal-dual iterations. The efficiency of the quadrature improves somewhat after several iterations (see \textbf{bold}), but begins to degrade beyond that point. \emph{From left to right}: the primal-dual iteration (``Iter.''), the condition number of the Jacobian of moment-matching residuals evaluated at the warm start (``$\kappa(\mathsf{D}F)$"), the max norm of the moment-matching residual for the warm start (``$\|R_0\|_\infty$''), the max norm of the final moment-matching residual (``$\|R_{\operatorname{final}}\|_\infty$''), the number of quadrature points (``$|X|$''), the efficiency of the quadrature (``Eff.''), and the number of Gauss-Newton iterations needed for convergence (``\#GN'').}
% \end{table}

\section{Fast construction of efficient cut cell quadratures}\label{sec:fast-construction}

\subsection{A fast primal heuristic}\label{sec:fast-primal-heuristic}

% \begin{figure}
%   \centering
%   \includegraphics[width=\textwidth]{figures/ipmstats.png}
%   \caption{\emph{Left}: TODO. \emph{Right}: the ``clustered'' quadrature ($\times$'s), and the optimized quadrature ($\circ$'s). All nodes in both the clustered quadrature and subsequent optimized quadrature are strictly interior to $\Omega$, and the weights are positive (lying in a range of about 0.014 to 0.239, the exact values being unimportant).}
% \end{figure}

\begin{figure}
  \centering
  \includegraphics[width=0.33\linewidth]{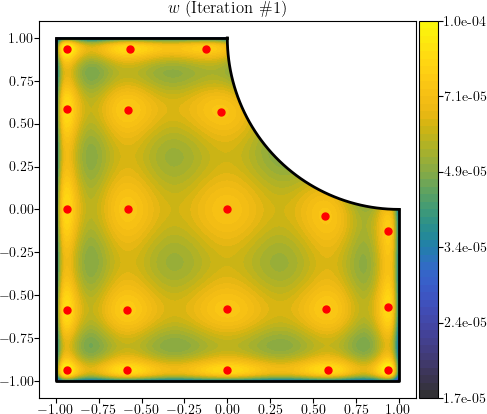}\vspace{0.01\linewidth}%
  \includegraphics[width=0.33\linewidth]{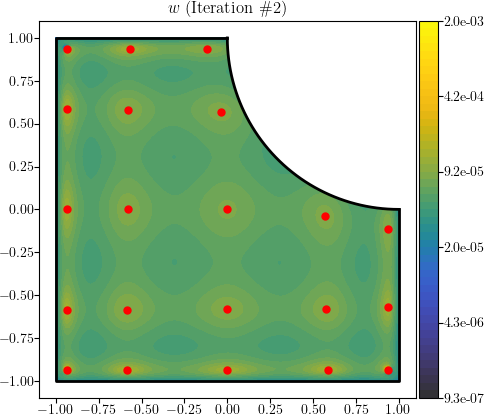}\hspace{0.01\linewidth}%
  \includegraphics[width=0.33\linewidth]{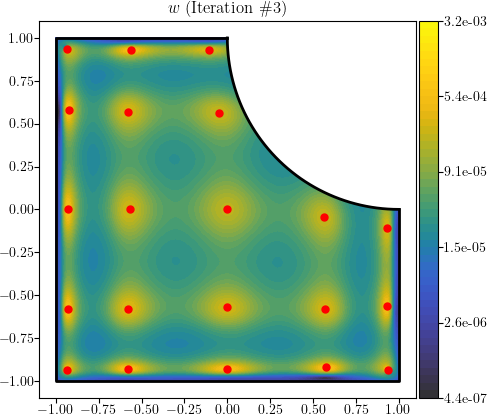}\vspace{0.025\linewidth}
  \includegraphics[width=0.33\linewidth]{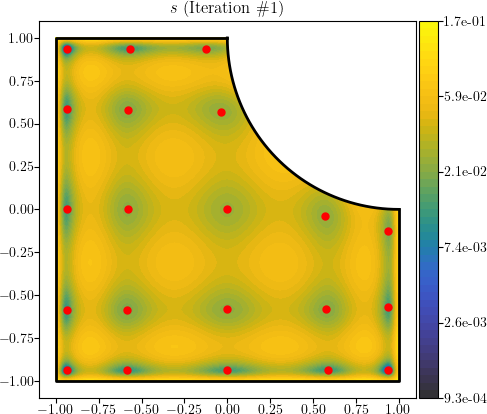}\vspace{0.01\linewidth}%
  \includegraphics[width=0.33\linewidth]{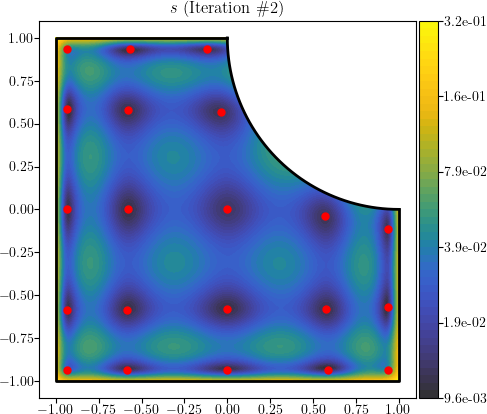}\hspace{0.01\linewidth}%
  \includegraphics[width=0.33\linewidth]{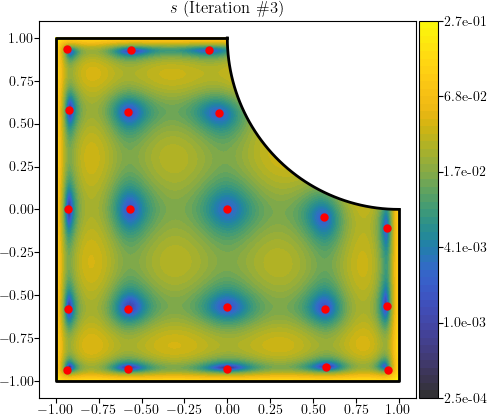}
  \caption{For the same test problem as in the left column of \Cref{fig:simple-nnls-2d} and \Cref{fig:ryu-boyd-2d}, we plot the first three iterations' primal and dual variables for the predictor-corrector interior point method applied to the Ryu--Boyd quadrature LP. \emph{Left to right}: the iterations of the interior point method. \emph{Top}: the primal variables. \emph{Bottom}: the dual variables. For each iteration, the local maximum clustering is applied to produce a warm start for Gauss--Newton used to solve the moment-matching equations (all Gauss--Newton iterations immediately achieve quadratic converge and produce quadratures with a residual on the order of machine precision). The red dots denote locations of local maxima.}\label{fig:primal-dual-2d-viz}
\end{figure}

Solving the Ryu--Boyd LP from \Cref{sec:ryu-boyd-quadrature} to optimality with a black-box LP solver is accurate but often expensive.
By examining the behavior of different LP solvers, we are led naturally to a fast primal heuristic.
While the simplex method reveals information about individual quadrature points as it proceeds, interior point methods follow the central path in the interior of the feasible set and converge to the optimal solution (which is unique if the residual function is chosen appropriately; see \Cref{sec:basis-pursuit-quadrature}).
This suggests that interior point iterates might contain global information about the quadrature rule, and indeed they do.
As shown in \Cref{fig:primal-dual-2d-viz}, the primal weights become increasingly concentrated around the final quadrature nodes, while the dual variables develop local extrema at precisely the same locations.
Somewhat surprisingly, these local extrema become discernible after a single iteration in the case of the primal variable, and immediately in case of the dual (slack) variables.
This observation motivates the fast heuristic introduced below.

We start from a standard primal-dual interior point method (``Algorithm MPC'') presented by Wright~\cite{wrightPrimaldualInteriorpointMethods1997}, which is a particular instance of Mehrotra's predictor-corrector scheme~\cite{mehrotraImplementationPrimalDualInterior1992}.
This method iteratively computes affine-scaling and centering directions for the KKT system of \eqref{eq:ryu-boyd-LP}, and then takes damped steps along these directions to remain in the strictly feasible region as it tracks the central path.
A lightly specialized version of the algorithm follows:
\begin{sfpalgorithm}[Primal-dual predictor-corrector method for \eqref{eq:ryu-boyd-LP}.]
  \begin{enumerate}[nosep]
  \item Form the normal matrix $\mV^\top\mV$ and compute its LU decomposition.
  \item Initialize $(\mw_0,\mlam_0,\ms_0)$ by solving the least squares systems
    $\mw_0 = \mV{(\mV^\top\mV)}^{-1}\mI$ and
    $\mlam_0 = {(\mV^\top\mV)}^{-1}\mV^\top\mr$,
    then set the slack variable $\ms_0=\mr-\mV\mlam_0$ and shift so that $\mw_0>0$, $\ms_0>0$.
  \item For $k=0,1,\ldots$ until convergence:
    \begin{enumerate}[nosep]
    \item Compute the affine-scaling search direction by solving the Newton system for the KKT residual.
    \item Compute the centering direction using Mehrotra’s correction term.
    \item Combine both directions and choose primal/dual step lengths that maintain strict positivity.
    \item Update $(\mw_k,\mlam_k,\ms_k)$ and repeat.
    \end{enumerate}
  \end{enumerate}
  \label{alg:ryu-boyd-MPC}
\end{sfpalgorithm}
\noindent The initialization step produces the least squares quadrature (\cref{fig:least-squares-quadrature}), and subsequent predictor-corrector iterations move the primal weights and dual variables along the central path toward the sparse, optimal Ryu--Boyd solution, which is sparse.

We plot several iterates of \Cref{alg:ryu-boyd-MPC} in \Cref{fig:primal-dual-2d-viz}.
From these plots, we can see straightaway that the clusters found in \Cref{fig:ryu-boyd-2d} coincide with the local extrema of the Lagrange multiplier fields.
We can run \Cref{alg:ryu-boyd-MPC} until it converges, but since the algorithm converges linearly, this can easily take 10-20 steps.
Instead, we find the local extrema of the primal or slack variable and use those nodes as a warm start for Newton's method (\Cref{sec:gauss-newton}).
An algorithm based on finding local maxima of some iterate $\mw_k$ might look like:
\begin{sfpalgorithm}[Heuristic primal solver for computing Ryu--Boyd quadratures.]
  \begin{enumerate}[nosep]
  \item Run \Cref{alg:ryu-boyd-MPC} to compute $\mw_k$ for $k > 0$ small; say, $k = 1$ or $2$.
  \item Find all strict local maxima of $\mw_k$.
  \item For the grid nodes corresponding to the local maxima, run \Cref{alg:least-squares} to compute quadrature weights (they will generally be positive).
  \item Run Gauss--Newton to optimize the result.
  \end{enumerate}
  \label{alg:ryu-boyd-MPC-primal-heuristic}
\end{sfpalgorithm}

\begin{figure}
  \centering
  \includegraphics[scale=0.75]{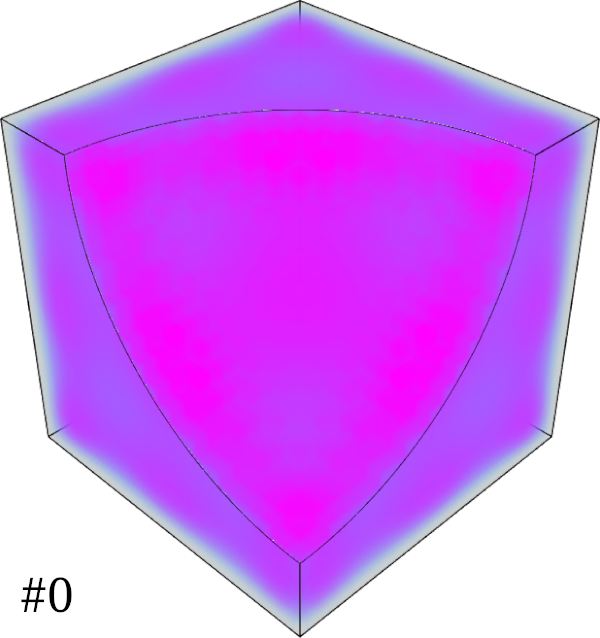}\hspace{0.015\linewidth}%
  \includegraphics[scale=0.75]{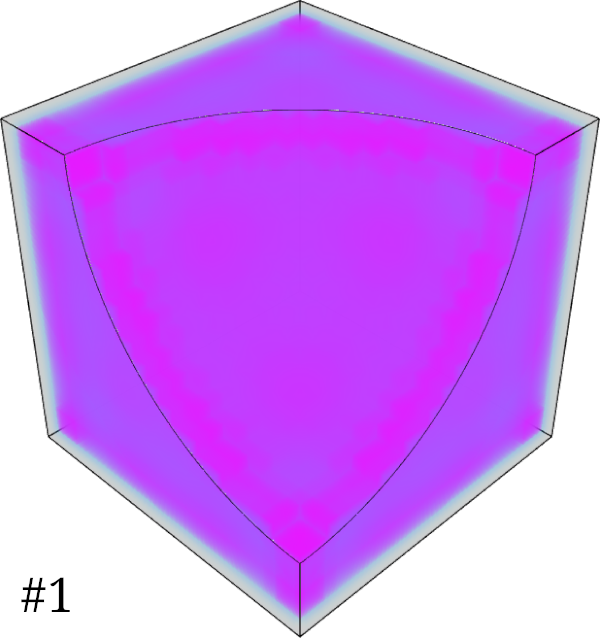}\hspace{0.015\linewidth}%
  \includegraphics[scale=0.75]{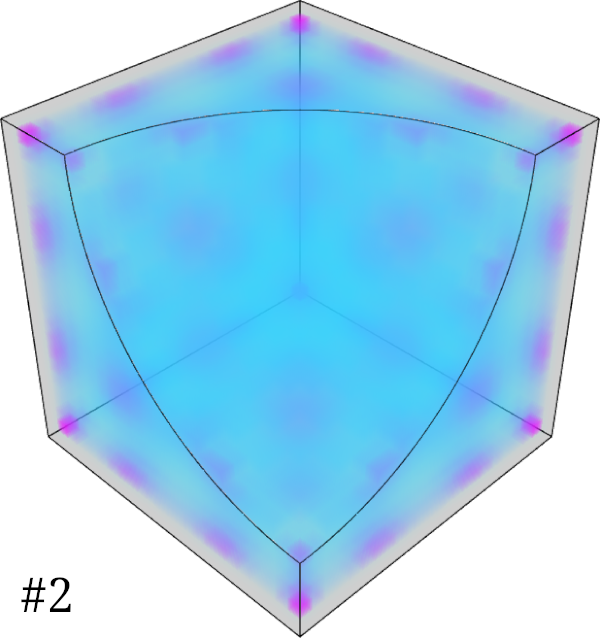}\hspace{0.015\linewidth}
  \includegraphics[scale=0.75]{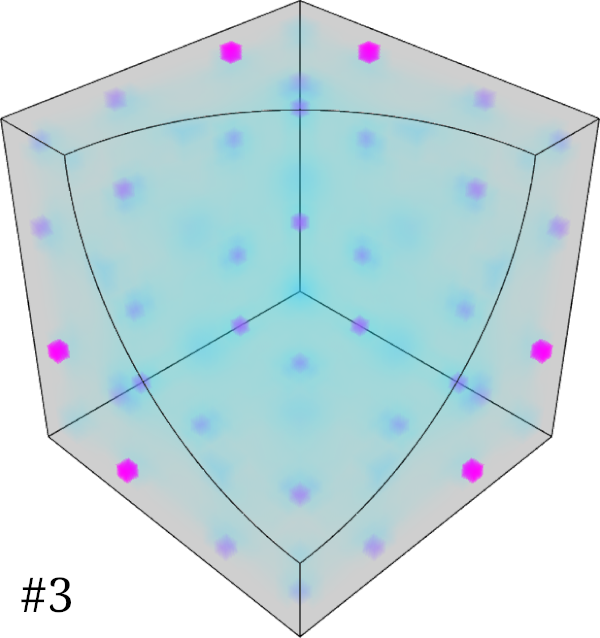}\hspace{0.015\linewidth}%
  \includegraphics[scale=0.75]{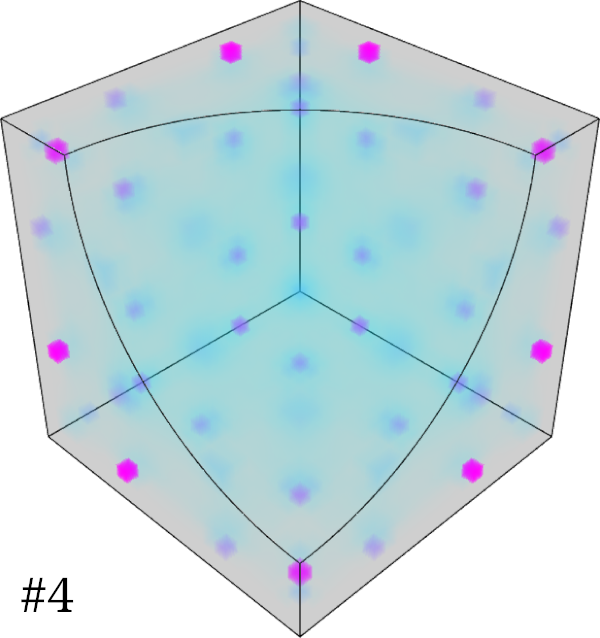}\hspace{0.015\linewidth}%
  \includegraphics[scale=0.75]{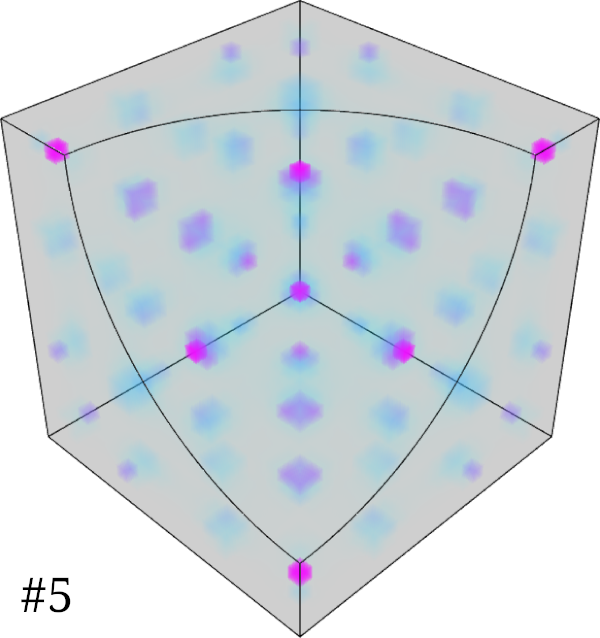}\hspace{0.015\linewidth}
  \caption{\emph{From left to right, top to bottom}: primal variables for the interior point method (the least squares quadrature for the grid), followed by the first several interior point primal iterations. The domain is $\Omega = [-1, 1]^3 - B_{3/2}((1,1,1))$. Despite appearances, after the first iteration, the local maxima of the primal variable give an initial iterate for Gauss--Newton which immediately leads to quadratic convergence when solving the moment matching equations. We neglect to plot the scale of the primal variables here (the color is scaled to the range $[0, \max(w_k)]$ in each subplot, where $w_k$ is the primal variable field at the $k$th iteration, with $0 \leq k < 8$). The purpose of the plot is to drive home the point that each successive interior point iteration serves to concentrate the measure nearer and nearer to the true quadrature points (which are almost surely not grid points). The key observation is that we can already deduce the location of each quadrature point from early iterations by looking at local maxima of the primal variable field.}\label{fig:gauss-lp-3d}
\end{figure}

\subsection{A fast dual heuristic}\label{sec:fast-dual-heuristic}

\begin{figure}
  \centering
  \includegraphics[width=0.435\linewidth]{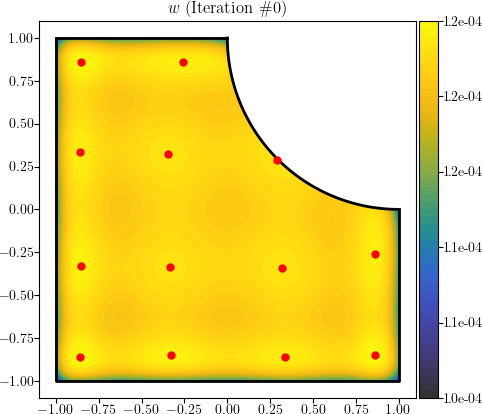}\hspace{0.04\linewidth}%
  \includegraphics[width=0.44\linewidth]{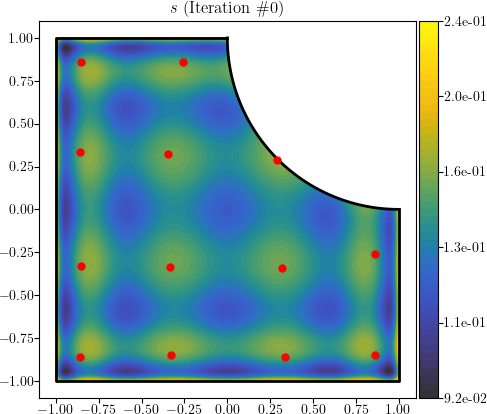}
  \caption{The leftmost column (``Iteration \#1'') is the warm start, which is essentially a least squares quadrature (see \Cref{fig:least-squares-quadrature}). The quadrature computed this way is clearly incorrect for Iteration \#0. However, it is interesting to observe that the interior local minima of the dual variable at Iteration \#0 correspond to the correct quadrature.}\label{fig:ryu-boyd-MPC-initialization}
\end{figure}

The least-squares initialization $\mw_0$ in \Cref{alg:ryu-boyd-MPC} does not, by itself, reveal the locations of the final quadrature nodes: as seen in \Cref{fig:ryu-boyd-MPC-initialization}, its local maxima do not coincide with the Ryu--Boyd support.  
However, an unexpected and striking phenomenon is already present at initialization: the local minima of the dual variable $\ms_0$ lie very close to the true quadrature nodes.
In fact, $\ms_0$ has a simple and transparent interpretation.

\paragraph{Interpretation of $\ms_0$}
Since the columns of $\mV$ are the basis functions sampled at grid points, the matrix $h^d \mV^\top \mV$ converges to the $L^2$ Gram matrix on $\Omega$ as $h \to 0$.  
Likewise, multiplication by $h^d \mV^\top$ approximates $L^2$ inner products:
\begin{equation}
  h^d (\mV^\top \mr)_i \;\to\; (\phi_i, r)_{L^2(\Omega)} \quad \mbox{as} \quad h \to 0.
\end{equation}
Hence, initializing $\mlam_0 \gets (\mV^\top \mV)^{-1} \mV^\top \mr$ computes the coefficients of the $L^2$ projection of $r$ onto $\mathbb{P}_N^d$ with respect to the discrete uniform measure determined by the nodes.
The corresponding slack variable $\ms_0 = \mr - \mV \mlam_0$ is then the residual of this discrete projection; i.e., the component of $r$ orthogonal to $\mathbb{P}_N^d$.
This explains why $\ms_0$ gives so much information about the structure of the optimal quadrature.
This motivates the following simple heuristic.
\begin{sfpalgorithm}[Heuristic dual solver for computing Ryu--Boyd quadratures]
  \begin{enumerate}[nosep]
  \item Compute the projection residual $\ms = \mr - \mV(\mV^\top \mV)^{-1} \mV^\top \mr$.
  \item Find all strict local minima of $\ms$.
  \item Compute weights at these grid points using least squares and refine with Gauss--Newton (as in \Cref{alg:ryu-boyd-MPC-primal-heuristic}).
  \end{enumerate}
  \label{alg:ryu-boyd-MPC-dual-heuristic}
\end{sfpalgorithm}

\paragraph{Continuous interpretation and justification}
We can consider the continuous analog of this construction, with $\{\phi_i\}$ a basis for $\mathbb{P}_N^d$ and $r\in L^2(\Omega)$.
Let $\m{G} \in \bbR^{m\times m}$ and $b \in \bbR^m$ be defined by:
\begin{equation}
  G_{ij} = (\phi_i, \phi_j)_{L^2(\Omega)}, \qquad b_i = (\phi_i, r)_{L^2(\Omega)}.
\end{equation}
let $\mc = \mG^{-1} \mb$, and define:
\begin{equation}
  s(\mx) = r(\mx) - \sum_{i=1}^m c_i \phi_i(\mx).
\end{equation}
As $h\to 0$, the vector $\ms$ computed above converges to $s$.
Using this formulation, we can understand why the heuristic is effective.

\begin{example}
  Let $\Omega = [-1, 1]$ and $r(x) = x^{2n}$, and let $\phi_i(x) = x^i$ for $0 \leq i < 2n$.  
  Then, $s(x) = P_{2n}(x)$, the Legendre polynomial of degree $2n$.  
  The strict local minima of $P_{2n}$ occur at every other zero of $P'_{2n}$, and these points are asymptotic to the zeros of $P_n$---precisely the nodes of the order-$n$ Gauss--Legendre quadrature rule.  
  In numerical experiments (for $1 \leq n \leq 200$), computing least-squares weights at these minima always produced positive weights, and the subsequent Gauss--Newton iteration achieved quadratic convergence immediately.
\end{example}
\noindent The success of this approach depends strongly on the choice of residual function $r$; it is unclear whether a universal choice exists that guarantees correctness for either of the heuristics in \Cref{alg:ryu-boyd-MPC-primal-heuristic,alg:ryu-boyd-MPC-dual-heuristic}.  
Nevertheless, for many practical choices of $r$, the dual heuristic performs remarkably well.

In $\bbR^d$, the local minima of the function $s(\mx)$ are simultaneous zeros of $d$ polynomials: the partial derivatives $s_{\mx_1}, \hdots, s_{\mx_d}$.
There is a known connection between simultaneous zeros of polynomials and quadrature rules which extends the classical results in one dimension to higher dimensions.
A summary of these results provided by Mysovskikh surveys the field as of 1980~\cite{mysovskikhApproximationMultipleIntegrals1980}.
A more recent monograph is also available~\cite{xuCommonZerosPolynomials2020}.
Understanding the relationship of \Cref{alg:ryu-boyd-MPC-dual-heuristic} to this literature could help put it on firmer footing.

\section{Comparison of quadrature algorithms}\label{sec:comparison}

\begin{figure}
  \includegraphics[width=0.5\linewidth]{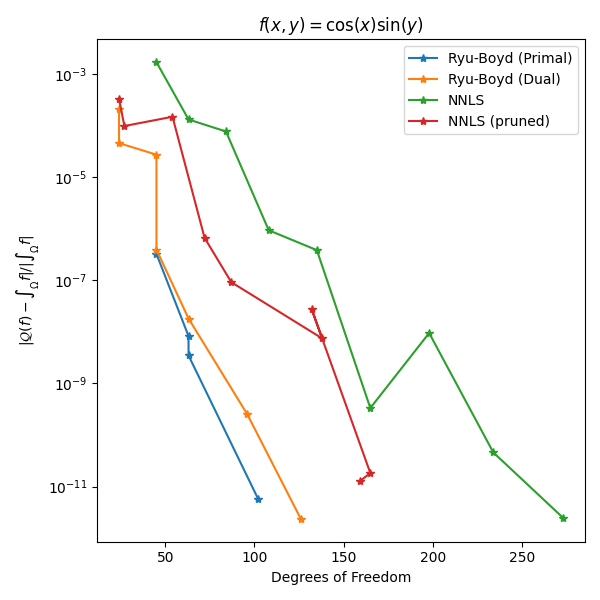}%
  \includegraphics[width=0.5\linewidth]{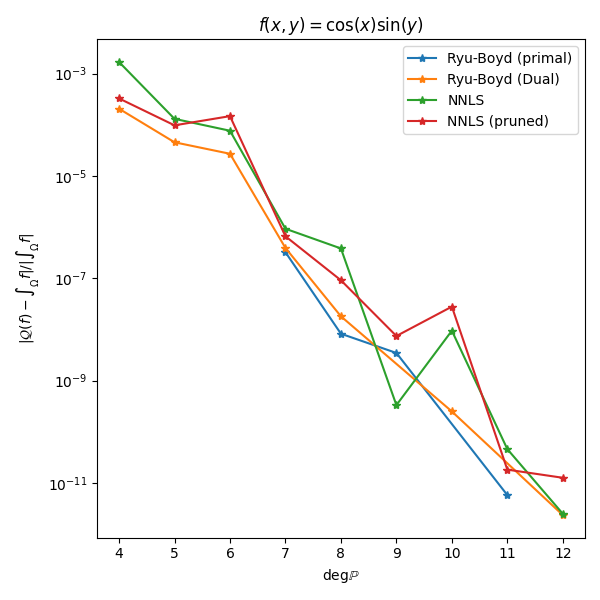}
  \caption{Relative integration error using each quadrature rule to integrate $f(x, y) = \cos(x)\sin(y)$ over $\Omega$. Here, $\Omega$ is the cut cell with a disk-shaped bite removed, as shown in e.g.\ \Cref{fig:ryu-boyd-MPC-initialization}. \emph{Left}: plotted vs.\ $\mbox{dof}(\mathcal{Q})$. \emph{Right}: plotted vs. the degree of the polynomial space being integrated.}\label{fig:errors}
\end{figure}

\begin{figure}
  \includegraphics[width=0.5\linewidth]{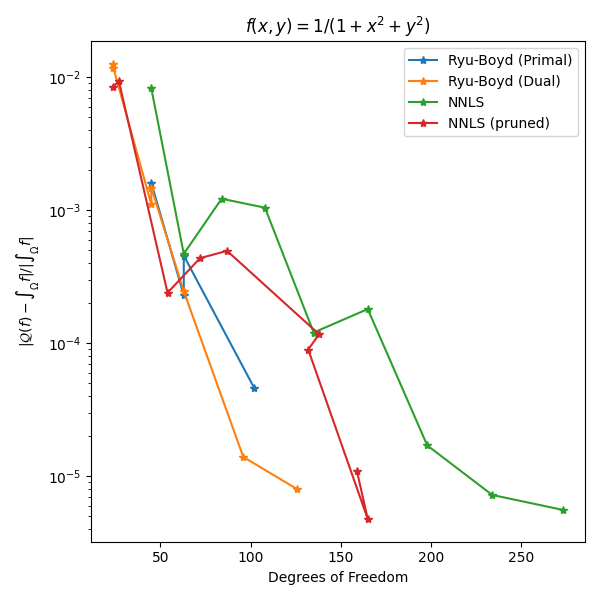}%
  \includegraphics[width=0.5\linewidth]{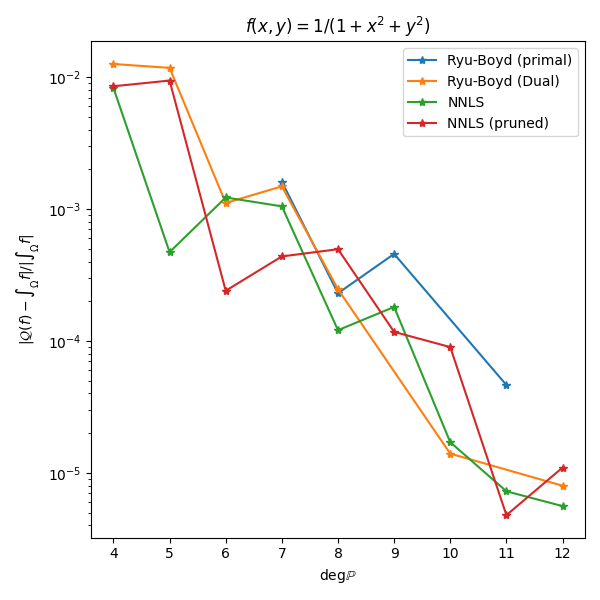}
  \caption{Relative integration error using each quadrature rule to integrate $f(x, y) = (1 + x^2 + y^2)^{-1}$ over $\Omega$. \emph{Left}: plotted vs.\ $\mbox{dof}(\mathcal{Q})$. \emph{Right}: plotted vs. the degree of the polynomial space being integrated.}\label{fig:errors-hard}
\end{figure}

\begin{figure}
  \includegraphics[width=0.5\linewidth]{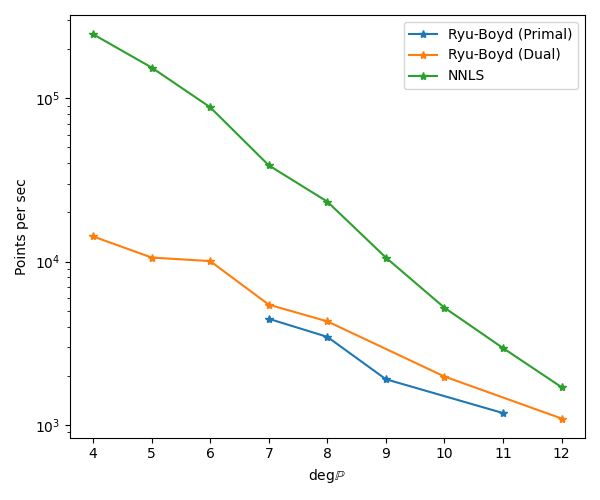}%
  \includegraphics[width=0.5\linewidth]{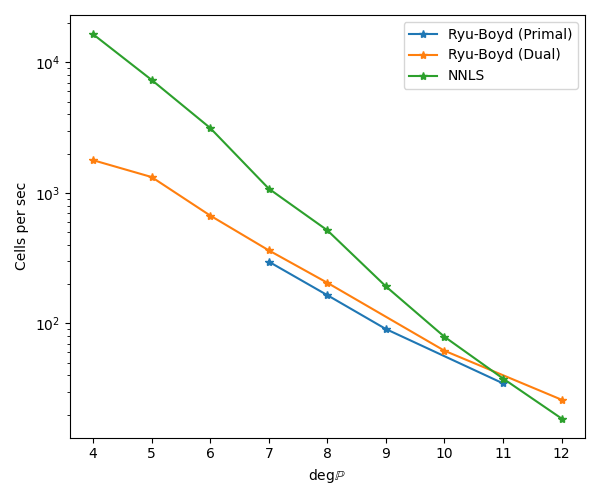}
  \caption{Timings result for building different quadrature rules on the cut cell shown in \Cref{fig:ryu-boyd-MPC-initialization}. \emph{Left}: numbers of quadrature points computed per second on a single core. \emph{Right}: number of quadrature rules (i.e. cut cells processed) per second per core.}\label{fig:timings}
\end{figure}

The purpose of this section is not to perform an exhaustive benchmark, but to illustrate representative accuracy and performance characteristics of the proposed quadrature construction methods on a typical cut-cell geometry.
We evaluate the accuracy, efficiency, and robustness of the quadrature algorithms developed in the preceding sections.
Each of these algorithms follows the same steps, at least conceptually:
\begin{enumerate}[nosep]
\item Discretize the bounding box of $\Omega$ into a uniform grid and discard all points lying outside the domain.
\item Compute the moments of each basis function in the polynomial space being integrated.
\item Solve some version of the moment-matching equations using an algorithm described in the preceding sections.
\end{enumerate}
To avoid overcomplicating things, we assume that $\Omega$ is described by a level set function and that we know its exact bounding box.
In particular, we consider the domain:
\begin{equation}
  \Omega = [-1, 1]^2 \backslash \Set{(x, y) \in \bbR^2 : (x - 1)^2 + (y - 1)^2 < 1}.
\end{equation}
In terms of the level set functions:
\begin{align}
  \psi_{\operatorname{box}}(x, y) &= \max(|x|, |y|) - 1, \\
  \psi_{\operatorname{disk}}(x, y) &= (x - 1)^2 + (y - 1)^2 - 1, \\
  \psi_{\operatorname{diff}}(x, y) &= \max(\psi_{\operatorname{box}}(x, y), -\psi_{\operatorname{disk}}(x, y)),
\end{align}
this can is equivalently given by:
\begin{equation}
  \Omega = \Set{(x, y) \in \bbR^2 : \psi_{\operatorname{diff}}(x, y) \leq 0}.
\end{equation}
For these examples, we assume the moments have been computed accurately to numerical roundoff in double precision using an adaptive Gauss quadrature on the domain and the divergence theorem.
For small $N$ (up to, say, $N = 15$ to 20), the monomial basis is accurate enough for our use---the exact conditions under which the poor conditioning of the monomial basis begins to cause problems have been studied~\cite{shenPolynomialInterpolationMonomial2023}.
Volume integration of multivariate monomials can be expedited by the divergence theorem.
One approach is given by the formula:
\begin{equation}
  \int_{\Omega} f(\mx) d\mx = \frac{1}{d+q} \int_{\partial\Omega} \mn(\mx) \cdot \mx f(\mx) d\mx.
\end{equation}
One proof directly applies the divergence theorem; it can also be shown using Euler's theorem on homogeneous functions~\cite{chinNumericalIntegrationHomogeneous2015}.

We first examine accuracy integrating the smooth function $f(x,y) = \cos(x)\sin(y)$.
Figure~\ref{fig:errors} shows the relative integration error as a function of the degrees of freedom $\mathrm{dof}(\mathcal{Q})$ and total polynomial degree $N$.
All methods exhibit the expected decay of the integration error as $N$ increases.
For a fixed number of degrees of freedom, the Ryu--Boyd quadratures and their heuristic variants consistently achieve the lowest errors, reflecting their ability to identify near-minimal node sets with well-conditioned weights.
The NNLS-based rules show somewhat larger errors but maintain the expected convergence trends.

Figure~\ref{fig:errors-hard} presents the corresponding results for the more challenging integrand $f(x,y) = (1 + x^2 + y^2)^{-1}$.
Errors decrease more slowly, as expected, but the relative performance of the methods is unchanged: the LP-based quadratures remain the most accurate, and the primal and dual heuristics track the full Ryu--Boyd rules closely.
These experiments provide evidence that the heuristics successfully capture the essential structure of the black-box Ryu--Boyd quadratures.

We next examine the construction time for the various rules.
Figure~\ref{fig:timings} reports throughput measurements on a single CPU core: the number of quadrature points computed per second (left), and the number of entire quadrature rules constructed per second (right).
The NNLS solver is fastest for lowest degrees, with the primal and dual heuristics further behind.
The Ryu--Boyd rules computed using a black box LP solver were computed too slowly to be of much use, so they were left out of the comparison.
Interestingly, the heuristic Ryu--Boyd methods appear to scale better with polynomial degree, eventually outpacing construction of NNLS rules.
These results indicate that the primal and dual heuristics inherit the accuracy of the LP formulation while achieving construction costs comparable to or better than NNLS at high degree.
As such, they offer a practical path to computing high-order cut-cell quadratures in production codes.

% --- section withheld: needs the subdivision discussion written ---
% \subsection{Composite rules in two dimensions}\label{sec:composite-rules}
%
% \todo{address subdivision here... if any of the quadrature algorithms fail, just subdivide until success, and then ``on the way back up'' run pruning algorithms to enforce the correct number of quadrature nodes...}
%
% \begin{figure}
%   \centering
%   \includegraphics[width=0.41\linewidth]{figures/nnls-disk.pdf}
%   \includegraphics[width=0.57\linewidth]{figures/disk-nnls-N8-d6.png}
%   \caption{Convergence of nonnegative least squares rules on the unit disk. \emph{Left}: absolute error computing the integral of $f(x, y) = \exp(x^2 + y^2)$ for different total degree $N$ (each different graph: $N = 2, 4, \hdots, 12$), and for varying numbers of background cells along the $x$ and $y$ axes, ranging from a single cell to 128 cells in each direction. \emph{Right}: a depiction of the NNLS rules for total degree $N = 6$ with 8 cells along each axis. Each NNLS rule has 28 points in it, matching the Tchakaloff bound (i.e., equal to the dimension of $\bbP_8^2$), and was solved to zero residual.}\label{fig:composite-nnls-disk}
% \end{figure}

\section{Conclusion}

In this work we have attempted to systematize, unify, and extend a surprisingly scattered body of ideas surrounding the construction of efficient quadrature rules on arbitrary domains.
The central theme is that once one can evaluate which candidate points lie in $\Omega$ and can compute the moment vector for a chosen polynomial space, the full complexity of meshing and handcrafted elementwise quadratures can, in principle, be bypassed.
Classical results such as Tchakaloff's theorem and tools from sparse optimization turn the problem into one of solving underdetermined moment-matching equations with positivity constraints.
By combining these ideas with fast moment computation and an efficient and robust inside-outside test, one can obtain practical algorithms that quickly deliver near-optimal rules.

Before concluding, we present a few open questions based on the work collected in this survey:
\begin{enumerate}
\item The structure of \Cref{alg:lawson-hanson-nnls} suggests that the complexity of $O(m^2 n)$ determined by \Cref{thm:lawson-hanson-nnls-complexity} is pessimistic.
  Is it possible to design a practical algorithm which computes a Tchakaloff efficient (i.e., with not more than $m$ nodes) quadrature in $O(m^3)$ time (leaving aside the cost of the computation of the moment vector)?
\item Is there some ``ideal'' Ryu--Boyd sensitivity function which is to be preferred among all others?
\item In \Cref{sec:fast-dual-heuristic}, we mention that there is a well-developed theory describing the connection between simultaneous zeros of multivariate orthogonal polynomials and minimal quadratures~\cite{xuCommonZerosPolynomials2020}.
  Can this theory be used to understand when \Cref{alg:ryu-boyd-MPC-dual-heuristic} is likely to fail and put it on a more secure foundation?
\item One issue touched on here only briefly is the use of quadrature nodes as collocation points (important in mass-lumped finite element simulations~\cite{voetMathematicalTheoryMass2023}).
  For a polynomial interpolant to be stable at high order, its Lebesgue constant should grow slowly.
  The number of nodes in a Tchakaloff efficient quadrature matches the dimension of the space being integrated---is there any way to adapt the algorithms presented here to ensure the computation of stable sets of interpolation nodes?
\end{enumerate}

\section{Acknowledgments}

Preliminary work for this research was carried out while employed at Coreform, LLC.
Many thanks to the employees and management of Coreform for providing a fruitful and collaborative environment and for motivating this work in the first place.
Special thanks are due to Drs.\ Derek Thomas, Michael Scott, Greg Vernon, David Kamensky, and Christopher Whetten for useful conversations and for providing helpful orientation within the world of computational mechanics.

\bibliographystyle{plain}
\bibliography{../sam,local}

\appendix
% number appendix equations A.1, A.2, ... (and reset at each appendix section)
\numberwithin{equation}{section}
% make \Cref say "Appendix A" rather than "Section A"
\crefalias{section}{appendix}

\section{Accuracy of the Lawson--Hanson algorithm}\label{sec:lawson-hanson}

In this section, we explain a detail of Lawson and Hanson's Fortran implementation of the Lawson--Hanson algorithm which appears to have been omitted from their book~\cite{lawsonSolvingLeastSquares1995} and is not explained in a comment in the original Fortran code~\cite{lawsonHansonFortranCode}.

For bookkeeping, we use ``numpy'' notation (i.e., 0-indexed ``MATLAB notation''), where, for example, $\mx_{i:j} = (x_i, x_{i+1}, \hdots, x_j-1)$.
Expressions like $\mA_{:,:k}$ are to be interpreted following this convention.
We also consider expressions with index set subscripts; e.g., given $I = \{i_1, ..., i_n\}$, we write $\mx_I = (x_{i_1}, \hdots, x_{i_n})$.
Of course, a set is an unordered collection of unique items---when we index a matrix's rows or columns with an index set, we always assume that the indexing happens in sorted order.
This avoids an abuse of notation and keeps things conceptually simple.
For an index vector $I \subseteq [n]$, we write $I^c \subseteq [n]$ for its complement.
Let:
\begin{equation}
  \mR_k = \parens{\mQ_k^\top \mV^\top}_{:|J_k|,J_k}, \qquad \mb_k= \parens{\mQ_k^\top \mI}_{:|J_k|}
\end{equation}
Note that:
\begin{equation}
  \parens{\mQ_k^\top \mV^\top}_{|J_k|:,J_k} = \m{0}_{(n - |J_k|) \times |J_k|}.
\end{equation}
At the $k$th step, we let $J_k \subseteq [n]$ be the set of inactive indices and compute a new tentative value of $\mw_k$ by solving the least squares problem on the set of inactive indices $J_k$:
\begin{equation}
  \begin{split}
    \mR_k \parens{\mw_k}_{J_k} &= \mb_k, \label{eq:lawson-hanson-least-squares-problem} \\
    \parens{\mw_k}_{J_k^c} &= \m{0}_{|J_k^c|}.
  \end{split}
\end{equation}
If any components of $\mw_k$ are negative, backtracking is used to find the largest step such that $\mw_k \geq 0$.
Afterwards, the Lagrange multiplier is recomputed from:
\begin{equation}
  \malpha_k = \mV \parens{\mV^\top \mw_k - \mI} \label{eq:lawson-hanson-lagrange-multipliers-1}
\end{equation}
If $n > m$, the least squares problem \eqref{eq:lawson-hanson-least-squares-problem} will be overdetermined until the last step, when it becomes square.

Lawson and Hanson elect to solve \eqref{eq:lawson-hanson-least-squares-problem} using a QR decomposition.
This QR decomposition is updated and downdated on the fly, with columns being added and removed from the orthogonal basis which spans some subset of the columns of $\mV^\top$.
For clarity, note the following:
\begin{enumerate}[nosep]
\item An index becoming \emph{active} (being added to the active set) corresponds to a quadrature node being removed from the rule; equivalently, a column of $\mV^\top$ being removed from the basis, and an entry of $\mw$ being reset to zero.
\item An index becoming \emph{inactive} is the same as a new quadrature node being added and an unused column of $\mV$ being incorporated into the QR decomposition's running orthogonal basis.
\end{enumerate}
Columns are added to the orthogonalized basis using a Householder reflector; they are removed using a Givens rotation.

Both $\mQ_k^\top \mV^\top$ and $\mQ_k^\top \mI$ are maintained throughout the iteration, where $\mQ_k \in \mathbb{R}^{m \times m}$ is the product of all Householder reflectors and Givens rotations applied throughout.
The matrix $\mQ_k$ is not stored explicitly.
Consequently, the least squares problem \eqref{eq:lawson-hanson-least-squares-problem} reduces to an upper triangular system.
No additional pivoting is applied for numerical stability other than what could be regarded as the problem-dependent pivoting induced by the choice of component to add or remove from the active set at each step.

The key detail is this.
At each step, the new vector of Lagrange multipliers is computed from:
\begin{equation}
  \begin{split}
    \label{eq:lawson-hanson-lagrange-multipliers-2}
    \parens{\malpha_k}_{J_k} &= \m{0}_{|J_k|}, \\
    \parens{\malpha_k}_{J_k^c} &= \squareb{\parens{\mQ_k^\top \mV^\top}^\top}_{J_k^c} \mQ_k^\top \mI,
  \end{split}
\end{equation}
The form of \eqref{eq:lawson-hanson-lagrange-multipliers-2} allows us to solve for the components of $\malpha_k$ as follows:
\begin{enumerate}[nosep]
\item Set $\left(\malpha_k\right)_J := \m{0}_n$.
\item Set $\left(\malpha_k\right)_{J^c} := \left(\mQ_k^\top \mV_k^\top\right)_{(|J|+1):}^\top \left(\mQ_k^\top \mI\right)_{(|J|+1):}$.
\end{enumerate}
We can see that na\"{i}vely computing $\malpha_k$ from \eqref{eq:lawson-hanson-lagrange-multipliers-1} requires us to first compute the residual $\mI - \mV^\top \mw_k$.
As $k$ increases and the residual reduces in magnitude, this residual computation can sometimes introduce an artificially large ``error floor'' as we begin ``approximating zero'' using the differences of moderately large numbers.
On the other hand, with \eqref{eq:lawson-hanson-lagrange-multipliers-2}, components of $\malpha_k$ are directly set to zero if they are active, or computed using a triangular solve if they are not.
In the latter case, the triangular system has been transformed from the original using a sequence of orthogonal transforms, preventing spurious numerical errors from accumulating.

We now justify our claim that Lawson and Hanson's stabilized expression for the Lagrange multipliers given by \eqref{eq:lawson-hanson-lagrange-multipliers-2} can be used instead of \eqref{eq:lawson-hanson-lagrange-multipliers-1}.

\begin{theorem}
  Expressions \eqref{eq:lawson-hanson-lagrange-multipliers-1} and \eqref{eq:lawson-hanson-lagrange-multipliers-2} are equivalent.
\end{theorem}

\begin{proof}
  First, note that since $\mQ_k$ is a (square) orthogonal matrix, we trivially have:
  \begin{equation}
    \malpha_k = \mV\mQ_k \parens{\mQ_k^\top\mV^\top \mw_k - \mQ_k^\top \mI}.
  \end{equation}
  We expand the first term in the parenthesis into parts corresponding to $J_k$ and $J_k^c$ and simplify to get:
  \begin{equation}
    \mQ_k^\top\mV^\top \mw_k = \parens{\mQ_k^\top \mV^\top}_{:,J_k}\parens{\mw_k}_{J_k} + \parens{\mQ_k^\top \mV^\top}_{:,J_k^c}\underbrace{\parens{\mw_k}_{J_k^c}}_{=0} = \parens{\mQ_k^\top \mV^\top}_{:,J_k}\parens{\mw_k}_{J_k}.
  \end{equation}
  We can further rewrite the matrix in this last expression as a block matrix to find that:
  \begin{equation}
    \parens{\mQ_k^\top \mV^\top}_{:,J_k} = \begin{bmatrix}
      \mR_k \parens{\mw_k}_{J_k} \\ \m{0}_{(n - |J_k|)\times |J_k|}
    \end{bmatrix}.
  \end{equation}
  Consequently:
  \begin{equation}
    \mQ_k^\top\mV^\top \mw_k - \mQ_k^\top \mI = \begin{bmatrix} \mR_k \parens{\mw_k}_{J_k} - \mb_k \\ - \parens{\mQ_k^\top \mI}_{|J_k|:} \end{bmatrix} = \begin{bmatrix} \m{0}_{|J_k|} \\ - \parens{\mQ_k^\top \mI}_{|J_k|:} \end{bmatrix}.
  \end{equation}
  Now, consider $\parens{\malpha_k}_{J_k}$.
  We have:
  \begin{equation}
    \begin{split}
      &\parens{\malpha_k}_{J_k} = \parens{\mV\mQ_k}_{J_k,:} \begin{bmatrix} \m{0}_{|J_k|} \\ - \parens{\mQ_k^\top \mI}_{|J_k|:} \end{bmatrix} \\
      &\hspace{3em} = \begin{bmatrix} \mR_k^\top & \m{0}_{|J_k| \times (n - |J_k|)} \end{bmatrix} \begin{bmatrix} \m{0}_{|J_k|} \\ - \parens{\mQ_k^\top \mI}_{|J_k|:} \end{bmatrix} = \m{0}_{|J_k|}.
    \end{split}
  \end{equation}
  For $\parens{\malpha_k}_{J_k^c}$, we can write:
  \begin{equation}
    \parens{\malpha_k}_{J_k^c} = \parens{\mV\mQ_k}_{J_k^c,:} \begin{bmatrix} \m{0}_{|J_k|} \\ - \parens{\mQ_k^\top \mI}_{|J_k|:} \end{bmatrix} = -\parens{\mV \mQ_k}_{J_k^c, |J_k|:} \parens{\mQ_k^\top \mI}_{|J_k|:}.
  \end{equation}
  This proves the claim.
\end{proof}

\end{document}